\documentclass[msom,nonblindrev]{informs3} 

\OneAndAHalfSpacedXII 

\usepackage{titletoc}
\usepackage{natbib}
 \bibpunct[, ]{(}{)}{,}{a}{}{,}%
 \def\bibfont{\fontsize{8}{10}\selectfont} 
 \def\bibsep{\smallskipamount} 
 \def\bibhang{24pt}%

\usepackage{tablefootnote}
\usepackage{hyperref}
\hypersetup{
	colorlinks=true,
	linkcolor=black,
	filecolor=black,
	urlcolor=black,
	citecolor=black,
}
\makeatletter
\def\UrlAlphabet{%
      \do\a\do\b\do\c\do\d\do\e\do\f\do\g\do\h\do\i\do\j%
      \do\k\do\l\do\m\do\n\do\o\do\p\do\q\do\r\do\s\do\t%
      \do\u\do\v\do\w\do\x\do\y\do\z\do\A\do\B\do\C\do\D%
      \do\E\do\F\do\G\do\H\do\I\do\J\do\K\do\L\do\M\do\N%
      \do\O\do\P\do\Q\do\R\do\S\do\T\do\U\do\V\do\W\do\X%
      \do\Y\do\Z}
\def\UrlDigits{\do\1\do\2\do\3\do\4\do\5\do\6\do\7\do\8\do\9\do\0}
\g@addto@macro{\UrlBreaks}{\UrlOrds}
\g@addto@macro{\UrlBreaks}{\UrlAlphabet}
\g@addto@macro{\UrlBreaks}{\UrlDigits}

\usepackage{endnotes}
\usepackage{makecell}
\usepackage{subfigure}
\usepackage{caption} 
\usepackage{multirow}
\usepackage{dcolumn}
\usepackage{booktabs}
\usepackage{amsmath}
\usepackage{bbm}
\usepackage{ulem}
\TheoremsNumberedThrough     
\ECRepeatTheorems

\def\qed{\hspace*{10pt}\hfill{$\square$}\hfilneg\par}
\EquationsNumberedThrough    

\newtheorem{thm}{Theorem}
\newtheorem{lem}{Lemma}
\newtheorem{cor}{Corollary}
\newtheorem{prop}{Proposition}

\newcommand{\ben}{\begin{equation}}   \newcommand{\een}{\end{equation}}
\newcommand{\bea}{\begin{eqnarray}}   \newcommand{\eea}{\end{eqnarray}}
\newcommand{\beq}{\begin{eqnarray*}}  \newcommand{\eeq}{\end{eqnarray*}}
\newcommand{\PaperTitle}{Service Deployment in the On-Demand Economy: Employees, Contractors, or Both?}

\begin{document}


\RUNAUTHOR{Lu, Weng, and Xiao}

\RUNTITLE{\PaperTitle}

\TITLE{\PaperTitle}

\ARTICLEAUTHORS{%
        \AUTHOR{Lijian Lu}
	\AFF{School of Business and Management, Hong Kong University of Science and Technology, HK,     \EMAIL{lijianlu@ust.hk}}
        \AUTHOR{Xin Weng}
	\AFF{Tsinghua-Berkeley Shenzhen Institute, Tsinghua Shenzhen International Graduate School, Tsinghua University, Shenzhen, China, \EMAIL{wengx20@mails.tsinghua.edu.cn}}
        \AUTHOR{Li Xiao}
	\AFF{Tsinghua-Berkeley Shenzhen Institute, Tsinghua Shenzhen International Graduate School, Tsinghua University, Shenzhen, China, \EMAIL{xiaoli@sz.tsinghua.edu.cn}}
} 

\ABSTRACT{%
The recent advancements in mobile/data technology have fostered a widespread adoption of on-demand or gig service platforms. The increasingly available data and independent contractors have enabled these platforms to design customized services and a cost-efficient workforce to effectively match demand and supply. In practice, a diverse landscape of the workforce has been observed: some rely solely on either employees or contractors, others use a blended workforce with both types of workers. 
In this paper, we consider a profit-maximizing service provider (SP) that decides to offer a single service or two differentiated services, along with the pricing and staffing of the workforce with employees and/or contractors, to price- and waiting-sensitive customers. Contractors independently determine whether or not to participate in the marketplace based on private reservation rates and per-service wage offered by the SP, while it controls the number of employees who receive per-hour wage. Under a single service, we show that the SP relies on either employees or contractors and identify sufficient and necessary conditions in which one workforce is better than the other. Under the optimal service deployment, we show that the SP offers either a single service relying solely on employees or contractors, or two differentiated services with a hybrid workforce depending on the service value and cost efficiencies of employees and contractors. 
Our analysis suggests that proliferating services with a blended workforce could improve the SP's profit significantly, and identifies conditions in which this value is significant. Our results provide an in-depth understanding and insightful guidance to on-demand platforms on the design of service differentiation and workforce models.
}%

\KEYWORDS{
sharing economy, on-demand workforce, staffing, pricing, incentive compatibility, service proliferation}

\maketitle

\section{Introduction} 
The on-demand or gig service platforms, facilitated by recent advancements in mobile/data technology, have been rapidly growing to be indispensable in the digital economy. These digital gig platforms have revolutionized the two sides (workers and customers) of many markets including the food delivery (DoorDash,  Meituan), ride-hailing (Uber,  Didi), hotel (Airbnb, Interhome), e-commerce and logistics (Amazon Fresh, SF Rush). On the demand side, the availability of massive customer data has enabled platforms to predict the heterogeneity of customers and then offer differentiated and customized services. On the supply side, the increasingly available independent contractors empower the platforms with the flexibility to optimize the traditional workforce of full-time employees. 
For example, the US independent workforce in 2022 represented $36\%$ of the employed population (\citealt{Mckinsey2022}). In contrast to the employees, or {\it firm-scheduled} agents, who receive hourly wages and could be directly controlled by the platforms, the contractors, or {\it self-scheduled} agents, are paid on a per-service wage and independently determine when and where to work\endnote{Throughout this paper, we use `contractors' and `self-scheduled agents', `employees' and `firm-scheduled agents', `on-demand service' and `self-scheduled service', `standard service' and `firm-scheduled service', interchangeably.}.

In practice, the gig platforms often dynamically adjust the offering of differentiated services with various forms of workforces.   
For example, Amazon introduced the on-demand grocery delivery service Amazon Fresh in 2007 as a complement to its existing standard service \citep{wells2018amazon}, and SF Express launched SF Rush, an on-demand intra-city instant delivery service, in 2016 to supplement SF Express and serve the impatient but high-value customers\endnote{\url{https://ir.sf-cityrush.com/en/investor-relations/}}.
Platforms providing food delivery service, on the other hand, had supplemented their on-demand services with firm-scheduled services to make the delivery more reliable. 
For instance, Meituan initially provided only on-demand service through crowdsourcing and later introduced a firm-scheduled service, Meituan Zhuansong, in rural areas by hiring full-time employees\endnote{\url{https://www.sixthtone.com/news/1004343}; \url{https://www.caixinglobal.com/2021-09-14/meituan-offers-glimpse-into-how-takeout-delivery-times-are-set-101773276.html}}, and DoorDash started hiring couriers as full-time employees in New York to provide reliable food delivery service\endnote{\url{https://www.ft.com/content/663359a2-f67a-4b86-938e-734145e452b0}}.
Besides, switching from one workforce to another entirely different model has also been witnessed.  For example, Uber Eats opted for a sub-contracting approach in Spain due to the new regulatory decree that reclassifies food delivery workers as employees by the Spanish government in 2021\endnote{\url{https://www.rfi.fr/en/business-and-tech/20211031-in-spain-delivery-riders-law-reshuffles-deck-for-take-away-market}.}, while on-demand startups Instacart and Sprig shifted away from independent contractors to rely on full-time employees\endnote{See Forbes article ``Yet Another On-Demand Service, Sprig, Switches Its Independent Contractors To Employees" \url{https://www.forbes.com/sites/ellenhuet/2015/08/06/on-demand-sprig-switches-independent-contractors-to-employe}.}.  

Despite the popularity of service proliferation and diverse landscape of different workforce models in the gig platforms, the optimal service deployment that effectively facilitates the match between proliferated services and the staffing of workforce has been limitedly understood. The service provider (SP) must decide whether or not to offer differentiated services via pricing and staffing. Specifically, to staff the workforce, the SP directly controls the number of employees while indirectly affecting the number of independent contractors via incentives with per-service wages, and decides on the portfolio of services together with corresponding prices and workforce to serve the price- and waiting-sensitive customers. 





\subsection{Our Approach and Contributions}
To understand this challenging problem, we study a profit-maximizing service provider's optimal service deployment in a stylized queueing model that consists of heterogeneous customers that are price- and waiting-sensitive, workers that are heterogeneous in earning sensitivities, and the service provider that determines whether or not to offer differentiated services, set the price, and manage the staffing of the workforce for each service to maximize the expected total profit. Specifically, we aim to address the following research questions:
\begin{itemize}
\item[1.]  What are the market conditions under which the SP is more profitable to operate a single service, on-demand service via independent contractors or standard service via employees?  How are the incentives of the SP, customers, and labor affected? Will the on-demand or standard service be beneficial to all of them, and when would this coordination happen? 
\item[2.] When does the SP offer proliferated services with both on-demand and standard service and how can these two services be coordinated? How would the supply market, such as the supply pool of contractors and hourly wage of employees, affect the pricing and staffing?
\item[3.] What are the impacts of service proliferation with both employees and contractors? How does adding another service affect the price and lead time of the existing service, as well as the SP profit, consumer surplus, labor welfare, and social welfare? Specifically, under which conditions would the service proliferation be significantly or marginally valuable?
\end{itemize}

To address the aforementioned research questions, we consider a single SP that offers a single service with a workforce of employees or contractors, or two differentiated services with a hybrid workforce, to serve price- and waiting-sensitive customers. 
The SP first decides the portfolio of services, along with the prices and staffing of the workforce by determining the number of employees and the per-service wage for contractors, to maximize the expected total profit. Then, the independent contractors, who are heterogeneous in the reservation rates, self-voluntarily decide whether or not to participate in the marketplace. Finally, the utility-maximizing customers, with heterogeneous sensitivity to waiting, select which service to join or balk at the service upon arrival. 

We start by examining the single service portfolio with employees or contractors. Our analysis reveals that the workforce plays a fundamentally different role in pricing and staffing. Specifically, under the workforce of employees, the SP may shut down the service when the service value is small and employees are expensive, and otherwise sets a sufficiently low price to serve the whole market. Under the workforce of contractors, on the other hand, the SP may be better off to reduce the lead time of patient customers by discouraging the impatient customers from joining the service. Besides, the structure of optimal staffing under each model is similar to the newsvendor solution with a safety capacity and the contractors are more efficient with a smaller safety capacity required. 
More importantly, we identify necessary and sufficient conditions under which the SP is better off with contractors or employees, and establish tight bounds between the performances of these two single systems.  Interestingly, we find that the incentives between the SP, labors, and consumers could be coordinated with employees or contractors, and provide market conditions in which this coordination occurs. 

We then investigate the optimal service deployment in a hybrid model. Specifically, the SP offers two heterogeneous services, that are differentiated in both prices and lead times, using a blended workforce with both employees and contractors. We find that, depending on the cost efficiency of these labor types, the optimal service deployment could be: providing a single service with employees (or contractors) if it has a significantly higher cost efficiency compared to the contractors (or employees), and two differentiated services with both employees and contractors otherwise when they have comparable cost efficiencies.

We further explore the impacts and values of proliferating one service with two differentiated services. First, we show that service proliferation may have two opposing effects, {\it competition effect} and {\it market expansion effect}, such that adding the on-demand service reduces the price and increases the lead time of the standard service due to the competition effect, while adding the standard service may either increase or reduce the price of the on-demand service depending on the net effect of these two effects. Second, we establish tight lower and upper bounds on the value of service proliferation to the SP, and show that adding another service with alternative labors could improve the profit of the SP significantly when the existing service has relatively lower cost efficiency. We also conduct extensive numerical simulations to illustrate the impacts of service proliferation on customers, labors, and social welfare. Finally, we demonstrate the robustness of our research findings by considering various extensions that relax the premises of the base model.

\subsection{Literature Review} \label{sec-lite}
Our paper is related to three streams of literature: on-demand service,  dual business models, and price/lead time differentiation.

{\bf On-demand service.} 
Our work is related to the growing stream of literature exploring operational issues in on-demand platforms (e.g., see \citealt{benjaafar2020operations,chen2020om,hu2021classics} for comprehensive reviews). Examples include the optimal price and wage \citep{cachon2017role,bai2019coordinating,bimpikis2019spatial,benjaafar2022labor,liu2023operating}, workforce scheduling \citep{besbes2024workforce}, supply and demand matching \citep{hu2020price,afeche2023ride,lobel2024detours,keskin2024order}, contract design \citep{feldman2023managing}, and  regulatory interventions \citep{yu2020balancing,bernstein2021competition}.  
\cite{taylor2018demand} shows that the delay sensitivity (or agent independence) could either increase or decrease the optimal price (or wage) in a queueing model under certain conditions regarding the uncertainty of service value and/or contractors' opportunity costs.  \cite{gurvich2019operations} find that the on-demand service might prompt the firm to opt for a lower staffing level and a reduced service level, resulting in lower profit, via a multi-period newsvendor model.
Fundamentally different from the above papers, our work aims to understand how the SP should deploy the differentiated services with employees and/or contractors, and to evaluate the impacts of service proliferation from various perspectives on the SP, consumers, labors and social welfare.


{\bf Price/lead time differentiation.}
Our paper also closely relates to the literature on service differentiation with price- and waiting-sensitive customers \citep{boyaci2003product,maglaras2018optimal,zhou2023benefit}.  
The first stream considers a single-line queueing system in which all heterogeneous customers join the same line and are served based on a priority rule (e.g., see \citealt{mendelson1990optimal,afeche2013incentive,maglaras2018optimal}). 
The second category studies two separate-line queueing systems, whereby heterogeneous customers could choose which line to join (e.g., see \citealt{guo2014downs, qian2017comparison, zhang2021design,zhou2023benefit}). \cite{boyaci2003product} show the impacts of the values of capacity costs on the degree of service differentiation. \cite{zhao2012lead} establish the efficacy of differentiation on customer valuation and service cost. 
\cite{zhou2023benefit} examine the impacts of privatization in a mixed duopoly service system with one private service provider and the other public service provider.
Our paper adopts a similar setting in which the SP may provide differentiated service to the price- and waiting-sensitive customers. However, our work significantly differs from the above papers on the endogenous decision on the workforce of labors and service offerings, i.e., the SP decides whether to offer a single service with either employees or independent contractors, or to provide two differentiated services with both labors, along with the pricing and staffing of each service. 

{\bf Dual business models.} This paper belongs to the recent emerging literature on queues with two distinct business models. 
Examples include pay-per-sale and pay-per-use pricing of information goods \citep{sundararajan2004nonlinear,ladas2022product},  
subscriptions and pay-per-use pricing in the service sector \citep{ladas2022product,wu2023bundling}, fee-for-service and bundled-price in  healthcare \citep{adida2017bundled,guo2019impact}. Our work is mostly related to the following papers that examine the optimal staffing of the workforce between full-time employees and independent contractors in the context of on-demand service. \cite{chakravarty2021blending} examines the impacts of two demand rationing strategies, preference for employees and equal opportunity for all drivers, on the blended workforce decision. Our model framework is fundamentally different from \cite{chakravarty2021blending} in the sense that the price- and waiting-sensitive customers select which service to join or balk at the service upon arrival, and services are differentiated in both prices and staffing of the workforce with two different types of labors.

\cite{dong2019impact} investigate the staffing problem for a cost-minimizing SP providing a single service with a workforce of flexible contractors and/or full-time employees. They consider a random number of servers and solve the staffing problem via a fluid model and a stochastic fluid model asymptotically, and find that the blended workforce could significantly reduce the staffing cost. Similar to their work, we also adopt a queueing-theoretic framework to examine the service deployment with employees and/or contractors. However, our paper differs fundamentally from their paper in the following aspects. First, we model the number of independent contractors by using an incentive-compatible approach, in which each contractor self-voluntarily decides whether or not to participate in the marketplace based on the private reservation rate and per-service wage offered by the SP.
Second, in addition to staffing, the SP also decides the portfolio of services along with the pricing and staffing for each service, to represent the popular practices of service differentiation. 
Furthermore, we evaluate the impacts of service proliferation via a blended workforce from various different perspectives including the SP, consumers, and labors, and social welfare, and provide tight bounds on the value of the blended workforce to the SP. 


\cite{lobel2024frontiers} consider a newsvendor setting with a single service that could be served by employees, or contractors, or a hybrid model and show that the value of the hybrid model is marginal. In contrast, we adopt a queueing-theoretic framework with price- and waiting-sensitive customers to explore the impacts of service proliferation with two differentiated services, thus, our model formulation is fundamentally different from \cite{lobel2024frontiers}. We show that the value of the hybrid model, which offers two differentiated services with a blended workforce, could be significant and provide tight bounds on this value.

{\bf Outline.} The rest of this paper is organized as follows. We first present the modeling framework (service mechanism, customer selection, platform profit, labor welfare, consumer surplus and social welfare) in Section \ref{sec-model}. We then study the single service model with a workforce of employees or independent contractors in Section \ref{sec-SO}, and investigate the hybrid model with two differentiated services with both employees and contractors in Section \ref{sec:Hybrid}. 
Section \ref{sec-extension} illustrates the robustness of our findings by studying various extensions and Section \ref{sec-conclusion} concludes this paper with potential future directions.
  
\section{Problem Formulation} \label{sec-model}
We consider a service provider (SP), or a firm, that offers customers access to two service channels, the {\it standard service} and the {\it on-demand service}, which differ in both the prices and the lead times\endnote{Throughout this paper, we use lead time, time of delay, and expected  sojourn time interchangeably, representing the total time (expected waiting time + expected service time) that a customer spends in the system.}. We use the subscript `$s$' to denote the standard service channel and `$o$' to denote the on-demand service channel. 

{\bf Customer behavior.}  Customers arrive randomly and sequentially according to a Poisson process with arrival rate $\Lambda$, which is referred to as the market size. Customers are sensitive to the price and lead time, with heterogeneous sensitivities to waiting. 
Specifically, by selecting a service with price $p$ and lead time $W$, a type-$\theta$ customer derives the following expected utility:
\[
U(p,W \;|\; \theta) = V - p - \theta W,  
\] 
where $V$ represents the value of service, $p$ and $W$ are the price and lead time of the chosen service, respectively, and $\theta \in  [0,1]$ measures customers' heterogeneity in the waiting. For the convenience of analytical tractability, we assume that $\theta$ follows a uniform distribution. The above utility model is widely used in the literature (see, e.g., \citealt{afeche2016optimal,zhang2021design}) capturing the heterogeneous sensitivities of waiting among the price- and waiting-sensitive customers. Upon arrival, each customer selects which service to join by choosing the service that generates the highest non-negative utility, and balks if this utility is negative\endnote{We assume that the customers do not observe  real-time information about the queue and can rationally anticipate the average waiting time of the queue. Thus, the expected utility of a type-$\theta$ customer selecting a service with price $p$ and expected lead time $W$ is given by ${\sf E}\left[U(p,\tilde{W}|\theta)\right] =  V - p - \theta {\sf E}\tilde{W} = V-p-\theta W$. We also assume that customers select the service only at the time of arrival and cannot change the service afterward.}. For notation convenience, we use $U_i$ as a brief for the expected utility of selecting service $i$, i.e., $U(p_i,W_i|\theta)$, for each service $i=o,s$.

{\bf Service mechanism.} The SP can employ  two different service mechanisms to achieve the desired lead time target, which includes both the service time and waiting time. 
We evaluate the lead time under each service mechanism via a multi-server queue, in which customers are served on a first come first served (FCFS) rule. For analytical tractability, we follow the common approximation in the literature (for example, see \citealt{wang2019offering,armony2021pooling,zhang2021design,hu2024regulation}) that approximates the lead time in a multi-server queue using an $M/M/1$ queue with the same service rate.  
Specifically, for each service $i\in\{s,o\}$, the total service rates are $k_i\mu_i$, where $k_i$ is the number of servers (or agents) and $1/\mu_i$ is the average service time per customer, and the lead time from $M/M/1$ queue is given by: 
\begin{eqnarray}
W_i=1/(k_i\mu_i -\lambda_i),\label{W}
\end{eqnarray}
where $\lambda_i$ is the arrival rate of service $i$. These two service mechanisms for the standard service and on-demand service differ in the supply cost functions due to distinct contracts between the service-providing agents and the SP, with details to be discussed below. 

\begin{itemize}
    \item[{\it Firm-scheduled agents (or employees).}] The SP hires regularly-contracted employees, or {\it firm-scheduled agents}, to provide the standard service. These agents are employed on a full-time basis with a specific {\it hourly wage} $w_s$. By adjusting the number of employees $k_s$, as well as the operational schedules of these employees, the SP can achieve a desired target of lead time, at a total cost of $C_s(k_s) = w_s k_s  $ per hour.

    \item[{\it Self-scheduled agents (or contractors).}] The SP outsources the on-demand service to {\it self-scheduled agents}, or contractors, who are paid by a {\it per-service wage} $w_o$. 
    Unlike the aforementioned employees, each contractor independently and voluntarily decides whether to accept or reject providing service, based on the per-service wage offered by the SP and her outside options. 
    Specifically, the contractors are {\it incentive-compatible} and participate in providing on-demand service if and only if doing so is more profitable. That is, each contractor evaluates the expected earning per hour of providing on-demand service, i.e., $\lambda_o w_o/k_o$ where $\lambda_o  /k_o$ is the average number of on-demand services per hour for each contractor, and participates if and only if this earning rate exceeds her reservation rate $r\in [0,1]$. We assume that contractors are earning-sensitive with heterogeneous reservation rates that are independently drawn from a uniform distribution. 
Therefore, the equilibrium number of participating contractors $k_o$ satisfies:
\begin{eqnarray}
k_o = K {\sf Pr} \left( \frac{\lambda_o w_o}{k_o} \ge r \right),\label{ko}
\end{eqnarray}
where $K$ is the number of all potential self-scheduled agents, i.e., the supply pool. We remark that, to achieve a lead time target $W_o$ for the on-demand service, the SP can only indirectly influence self-scheduled agents to participate by adjusting the per-service wage $w_o$, at the cost rate of $C_o(w_o) = \lambda_o w_o $ per hour. 
\end{itemize}

{\bf Service deployment and pricing.} The SP determines the service portfolio, i.e., by selecting which service to operate, and the corresponding prices, the number of employees and/or per-service wage for contractors, aiming to maximize his total profit. In particular, based on the service mechanisms included, there are three service portfolios: a single model with standard or on-demand service, and a hybrid model with both services.  Throughout this paper, we use superscript `$S$', `$O$', `$*$' to denote the service portfolio with standard service only, on-demand service only, and hybrid model, respectively.  Notations are summarized in Table \ref{tab-notation}.
 
{\it Standard service only.} The system is reduced to a one-tier service system with only the standard service channel. Upon arrival, each customer has two options, to either join the service or balk, based on her expected utility derived from the service. The SP decides the price of standard service and the number of employees to maximize his expected profit, i.e., by solving the following optimization problem: 
\bea
\pi^S = \left\{\begin{array}{rl}
\max\limits_{p_s,k_s} &  \pi(p_s,k_s) = p_s \lambda_s - C_s(k_s) \\
\mbox{s.t.} 
& \lambda_s = \Lambda {\sf Pr}(U_s \geq 0)  ,  \ W_s=1/(k_s\mu_s -\lambda_s).
\end{array}\right. \label{model-S}
\eea

{\it On-demand service only.} This is a one-tier service system with only the on-demand service channel. Upon arrival, customers decide whether to join the service or balk, based on their utilities derived from the service. The SP decides the price of on-demand service and the per-service wage to maximize his expected profit, i.e., by solving the following optimization problem:
\bea
\pi^O = \left\{\begin{array}{rl}
\max\limits_{p_o,w_o} &  \pi(p_o,w_o) = p_o \lambda_o - C_o(w_o) \\
\mbox{s.t.} 
& \lambda_o = \Lambda {\sf Pr}(U_o \geq 0)  , \ W_o=1/(k_o\mu_o -\lambda_o), \  k_o = K {\sf Pr} \left( \frac{\lambda_o w_o}{k_o} \ge r \right),
\end{array}\right. \label{model-O}
\eea
where the last equality in the constraints reflects the incentive-compatible constraint for the contractors to voluntarily participate in providing on-demand services.

{\it Hybrid model.} The SP provides both the standard and on-demand services, and each arriving customer decides which service to join or balk at the services based on the expected utility associated with each service. The SP determines the prices of these two services, the number of employees, and the per-service wage for contractors to maximize his expected total profit.  The optimization problem for the SP is formulated as: 
\begin{align}\label{opt}
&\pi^* =   \left\{\begin{array}{rl}
 \max\limits_{p_s,p_o,k_s,w_o} & \pi(p_s,p_o,k_s,w_o)= p_s \lambda_s - C_s(k_s) + p_o\lambda_o-C_o(w_o)   \\
\mbox{s.t.} 
& \lambda_s = \Lambda {\sf Pr}(U_s \ge 0, U_s \ge U_o), \  \lambda_o= \Lambda {\sf Pr} (U_o \ge 0, U_o \ge U_s), \\
& W_i=1/(k_i\mu_i -\lambda_i), i\in\{s,o\}, \  k_o = K {\sf Pr} \left( \frac{\lambda_o w_o}{k_o} \ge r \right).
\end{array}\right.
\end{align}  

The main objectives of this paper are to provide a comprehensive understanding of the value of service proliferation with different types of service providers. In addition to the SP profit, we also aim to evaluate the influence from different perspectives that include the service-providing agents, the customers, and the entire society. To this end, we formally define the consumer surplus as $CS =  \Lambda  \int_{0}^{1} \max_{i\in\mathcal{P}}\{U_i\} d \theta$, the labor welfare, and the total social welfare (i.e., the sum of SP profit, labor welfare, and consumer surplus) as: 
\beq
LW  & = &  \left( k_s w_s - k_s \int_0^1 r \mathbb{I}_{ \{r \leq w_s \}} dr \right)\mathbb{I}_{\{s\in \mathcal{P}\}}  + \left(K \int_0^1 w_o \mathbb{I}_{ \{r \leq \frac{\lambda_o w_o}{k_o}  \}} dr - K \int_0^1 r \mathbb{I}_{ \{r \leq \frac{\lambda_o w_o}{k_o}  \}} dr \right)\mathbb{I}_{\{o\in \mathcal{P}\}}, \\
SW  &=&  \left( \Lambda  \int_{0}^{1} (V - \theta W_s) \mathbb{I}_{ \{U_s = \max_{i\in\mathcal{P}}\{U_i\} \} } d \theta - k_s \int_0^1 r \mathbb{I}_{ \{r \leq w_s \}} dr \right)\mathbb{I}_{\{s\in \mathcal{P}\}} \nonumber \\
&&  + \left( \Lambda  \int_{0}^{1} (V - \theta W_o) \mathbb{I}_{ \{U_o = \max_{i\in\mathcal{P}}\{U_i\}\} } d \theta  - K \int_0^1 r \mathbb{I}_{ \{r \leq \frac{\lambda_o w_o}{k_o}  \}} dr \right)\mathbb{I}_{\{o \in \mathcal{P}\}}, 
\eeq
where $\mathbb{I}_{A}$ is an indicator function that equals one if condition $A$ is satisfied and zero otherwise, $\mathcal{P}$ denotes the service portfolio with $\mathcal{P} = \{s,\varnothing\},\{o,\varnothing\},\{s,o,\varnothing\}$ when the SP operates the standard service only, the on-demand service only, the hybrid model, respectively, and $\varnothing$  represents the balking decision of a customer with $U_\varnothing = 0$. 

\begin{table}[h]
\vspace{-0.1in}
\caption{Summary of Notations}
    \centering
    \begin{tabular}{c|l} \hline
      Notation &  meaning \\
       \hline 
       $V $ & value of service to each customer \\
       $\Lambda$ & customers' arrival rate (per hour), i.e., market size \\
       $K$ & number of all potential self-scheduled contractors, i.e., supply pool \\ 
       $\theta$ &  customer's sensitivity of waiting, $\theta\sim U[0,1]$ \\
       $r$ & hourly reservation rate of contractors, $r\sim U[0,1]$ \\
       $w_s$ & hourly wage of firm-scheduled employees \\
       $w_o$ & per-service wage of self-scheduled contractors \\
       $\mu_s, \mu_o$ & service rate (per hour) of an employee, a contractor  \\
       $k_s, k_o$ & number of participating employees, contractors \\
       $p_s,p_o$  & price of standard service, on-demand service \\
       $W_s,W_o$ & lead time (hours) of standard service, on-demand service \\  \hline
    \end{tabular}
    \label{tab-notation}
\vspace{-0.2in}
\end{table}

We close this section with a few discussions of our model assumptions. First, to derive a thorough understanding, we approximate the lead time using an $M/M/1$ queue for the convenience of tractability. We relax this assumption and consider an $M/M/k$ queue in Subsection \ref{subsec-multi} to show the robustness of the research findings.  
Second, our model assumes that the customers are heterogeneous in waiting and agents are heterogeneous in the reservation rate, both are independently drawn from uniform distributions. We show in Subsection \ref{subsec-distribution} that our research findings and key managerial insights continue to hold with a different formulation of the heterogeneity. Third, we assume that the two services have the same valuation, which reflects various business practices when the services have roughly the same quality and customers do not care about how they are served (e.g., package delivery service). In Subsection \ref{subsec:quality}, we extend our analysis to a different setting in which the services also differentiate in the quality so that customers derive a higher valuation from the on-demand service. 
Finally, our base model assumes that each type of agents provide only one specific service, that is, employees only provide standard services and contractors only serve on-demand services. In Online Appendix Section \ref{ec-sec-flexible}, we relax this assumption and consider a {\it flexible} setting in which the SP could utilize on-demand contractors for standard services and/or employees for on-demand services (see Figure \ref{fig-dedicated-flexible-servers} for an illustration), and we show that these two settings are equivalent for the SP, yielding the same expected profit when he endogenously sets service parameters. 

\section{Portfolio of Single Service: Standard vs On-demand} \label{sec-SO} 
This section considers the service portfolio with standard service or on-demand service. As discussed above, these two services differ in the prices, the lead times, and also the supply mechanisms. Specifically, the standard service relies on full-time employees with a per-hour wage, while the on-demand service is outsourced to part-time contractors with a per-service wage. Compared to the employees, the contractors could either benefit the SP because they will not be paid during their idle time thus saving labor cost for the SP, or hurt the SP due to the uncertain supply of labor and the high cost to attract 
self-scheduled independent contractors to participate. 

Our objectives in this section are to provide a thorough understanding of the impacts of the service mechanism on the SP profit, consumer surplus, labor welfare, and total social welfare. We are particularly interested in providing market conditions under which standard or on-demand service mechanism is more valuable from different perspectives. To this end, we first derive the  optimal solutions under the standard service only system \eqref{model-S} and the on-demand service only system \eqref{model-O}, and then apply these solutions to examine the impacts of service mechanisms on different performance metrics. For the sake of notational convenience, we refer to the one-tier system with the standard service as {\it System S}, and the one-tier system with the on-demand service as {\it System O}.  Define\endnote{We make a few comments on the monotonicity of thresholds defined in \eqref{eq-barK-barLo}. First, for any $x\ge 0$, the function $L_o^3 - x(1+L_o)$ first decreases and then increases in $L_0\ge 0$, with a unique positive root given by $\bar L_o(x) = \frac{2 \sqrt[3]{3} x + 
 \sqrt[3]{2} (9 x + \sqrt{3(27 - 4 x)} x)^{\frac{2}{3}}}{6^{\frac{2}{3}} (9 x + \sqrt{3(27 - 4 x)} x)^{\frac{1}{3}}}$ that increases in $x$. Thus, 
 $\bar L_o$ defined \eqref{eq-barK-barLo} can be written as $\bar L_o = \bar L_o\left(\frac{2\Lambda^2}{K\mu_o^2}\right)$. Second, it is obvious that $\frac{x}{(x-1)^3}$ decreases in $x\ge 1$, thus, $\bar K$ decreases in the value of service $V$ and $\overline{C}_s$ increases in the hourly wage of employees $w_s$.}
\bea 
\begin{array}{c}
\bar L_o =   \max\left\{L_o>0: \frac{L_o^3}{1+L_o}\le \frac{2\Lambda^2}{K\mu_o^2} \right\}, \ 
\bar K = \frac{2\Lambda^2 \sqrt{V \Lambda +1}}{\mu_o^2(\sqrt{V\Lambda + 1}-1)^3}, \ 
\overline{C}_s = \frac{w_s}{\mu_s} + 2\sqrt{\frac{w_s}{\Lambda \mu_s}},  \ 
\overline{C}_o = \frac{\bar L_o (3+\bar L_o)}{2\Lambda},
\end{array}
\label{eq-barK-barLo}
\eea 
where, as we show later, $\overline{C}_s$ and $\overline{C}_o$ reflect the per-service cost that includes the markdown of the price and staffing cost under the standard and on-demand mechanisms, respectively, when all customers join the service.  

\subsection{System S: Standard Service}  
We start with the single model with standard service and provide a full characterization of the optimal prices and number of employees that solves \eqref{model-S}. In the next proposition, we show the optimal strategy under the standard only service system is to either serve all customers when the value of service is sufficiently large or shut down the service without providing service to any customers. 

 
\begin{prop}\label{S-only} 
{\sc (Optimal solution of System S)} Under the single service model with the standard service, the optimal policy for the SP is characterized by:
\beq
\left\{\begin{array}{ll}
    \lambda_s^S =  \Lambda, \quad W_s^S = \sqrt{\frac{w_s}{\Lambda \mu_s}} , \quad \pi^S = \Lambda (V - \overline{C}_s), & \mbox{ if $V \ge \overline{C}_s$,} \\
    \lambda_s^S =  0, \quad W_s^S  = 0 , \quad \pi^S =0, & \mbox{ otherwise.}
\end{array}
\right.
\eeq
\end{prop}  

We prove the above proposition in two steps. First, we remark that the profit function in \eqref{model-S} is generally not concave, and we perform a one-to-one variable transformation from $(p_s,k_s)$ to $(\lambda_s,W_s)$, where $\lambda_s$ is the effective arrival rate and $W_s$ is  the lead time of the standard service. The transformed profit function is concave in the lead time $W_s$, thus, has a unique maximum. Second, the profit function at this optimal lead time is a linear  function of the effective arrival rate $\lambda_s$. Thus, the optimal effective arrival rate of System S, denoted by $\lambda_s^S$, equals one of the two boundaries, i.e., $0$ or the arrival rate $\Lambda$. We show that the profit function is increasing, thus, serving all customers is optimal if and only if the value of service is larger than the threshold $\overline{C}_s$, which represents the total cost per customer with employees (i.e., cost of serving $\frac{w_s}{\mu_s}$ plus cost of waiting $2\sqrt{\frac{w_s}{\Lambda \mu_s}}$) when all customers join the service. Specifically, as summarized in Table \ref{ec-table-opt-solutions}, the price of the standard service is set at the highest value to attract the most impatient customer joining the service (i.e., $p_s = V-W_s$), and the number of employees follows a `newsvendor' formula that includes both the mean and `safety' capacity (i.e., 
$k_s = \frac{\Lambda}{\mu_s}   + \sqrt{\frac{\Lambda}{\mu_s w_s}}$, where the square term represents the safety capacity).
The above results, illustrated in Figure \ref{fig-sol-S-O}, imply that, depending on the potential profit margin per service, the SP jointly adjusts the price and number of employees to either attract and serve all customers, or avoid serving any customers by setting a sufficiently high price (e.g., the value of service).


\subsection{System O: On-demand Service} 
In this subsection, we consider a single service model with on-demand service and provide a full characterization of the optimal price and per-service wage that maximize the SP profit. Analogous to the single model with standard service, 
the SP profit function in \eqref{model-O} fails to be concave; thus, we perform a one-to-one variable transformation with  decision variables $(\lambda_o,L_o)$, where $\lambda_o$ is the effective arrival rate and $L_o$ is the average number of customers in on-demand service. Using the transformation, we show that the profit function is concave in the effective arrival rate $\lambda_o$; thus, there exists a unique optimal arrival rate. We further show that the profit function at this optimal effective arrival rate is quasi-concave in $L_o$, therefore, it admits a unique optimal queue length.  
The optimal decisions under the single service model with on-demand service are then fully characterized in the following proposition. 


\begin{prop}\label{O-only}{\sc (Optimal solution of System O)} Under the single service model with on-demand service, the optimal policy for the SP is characterized by:
\beq
\left\{\begin{array}{ll}
    \lambda_o^O =  \Lambda, \quad L_o^O = \bar L_o , \quad \pi^O = \Lambda ( V  - \overline{C}_o), & \mbox{ if $K \ge \bar K$,} \\
    \lambda_o^O =  \frac{K}{\bar K} \Lambda, \quad L_o^O = \sqrt{V\Lambda + 1}-1, \quad \pi^O = \frac{K\mu_o^2 (\sqrt{V\Lambda +1}-1)^4}{4\Lambda^2}, & \mbox{ otherwise.}
\end{array}
\right. 
\eeq
\end{prop}

\begin{figure}[htb]
	\vspace{-0.1in}
	\FIGURE
	{
            \subfigure[~System S]{
                \begin{minipage}{0.5\textwidth}
                    \includegraphics[width=0.48\textwidth]{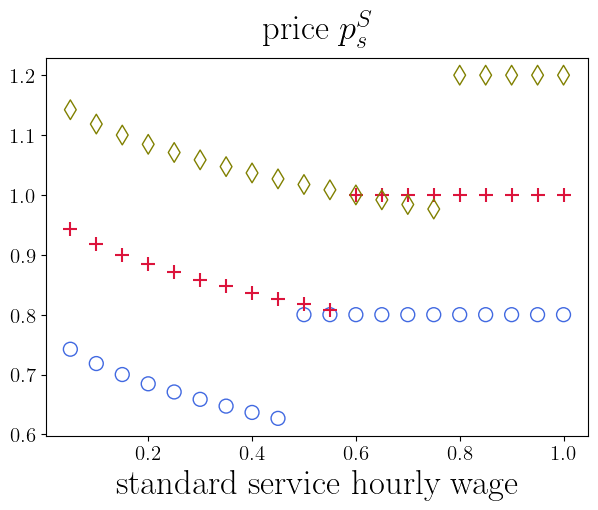} 
                    \includegraphics[width=0.48\textwidth]{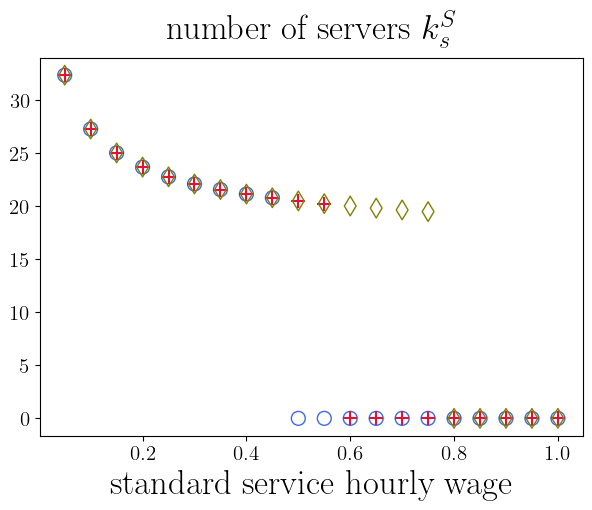} 
                    \vfill
                    \includegraphics[width=0.48\textwidth]{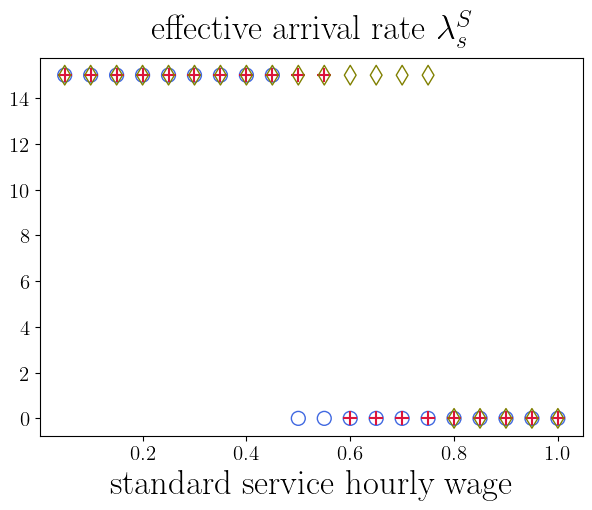} 
                    \includegraphics[width=0.48\textwidth]{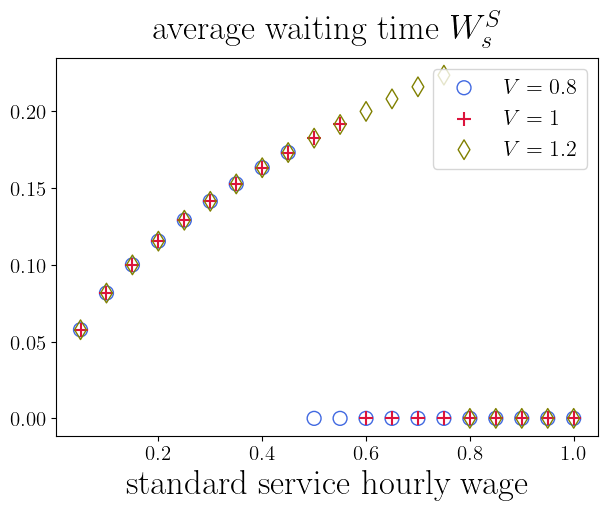} 
                \end{minipage}
            }
            \subfigure[~System O]{
                \begin{minipage}{0.5\textwidth}
                    \includegraphics[width=0.48\textwidth]{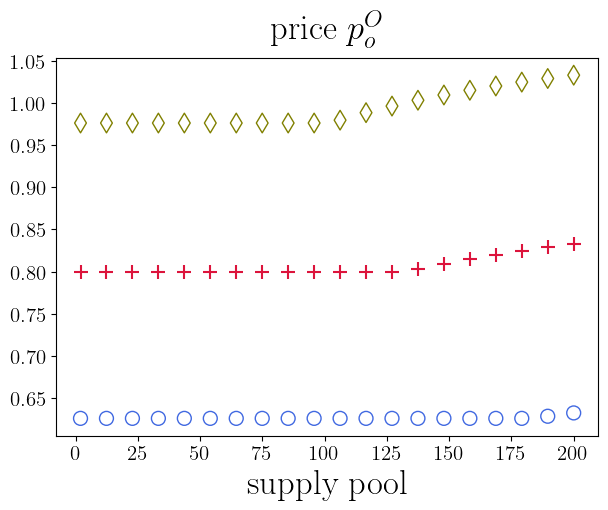} 
                    \includegraphics[width=0.48\textwidth]{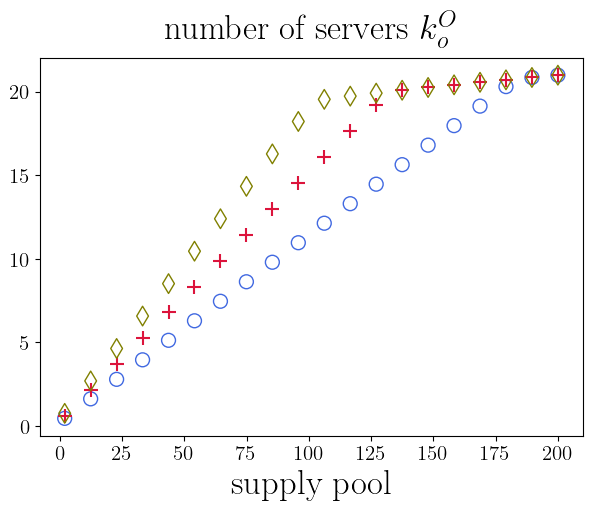} 
                    \vfill
                    \includegraphics[width=0.48\textwidth]{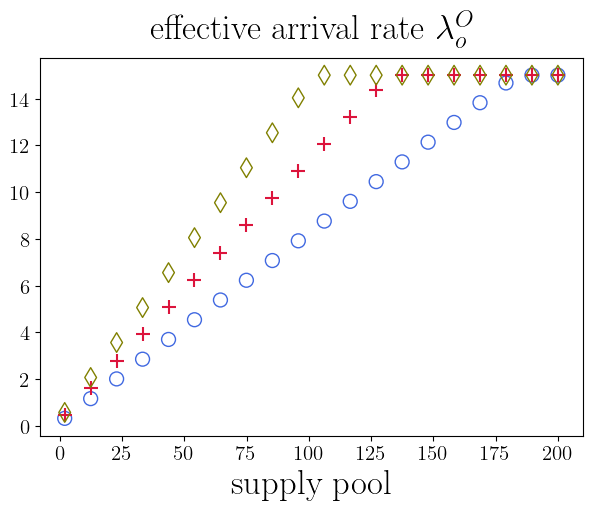} 
                    \includegraphics[width=0.48\textwidth]{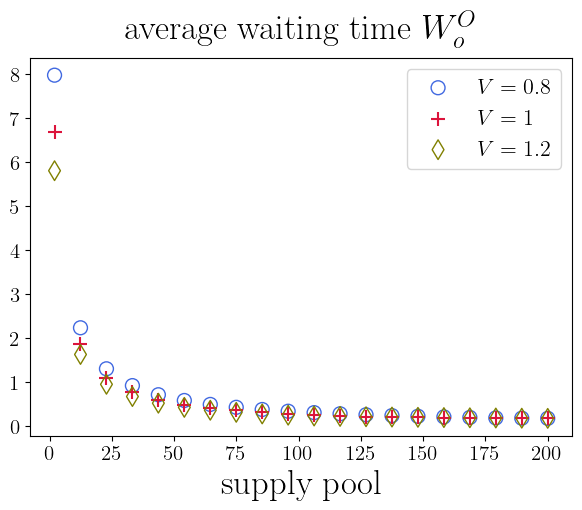} 
                \end{minipage}
            }
	}
	{\bf An visualization of optimal pricing and staffing: Standard vs On-demand \label{fig-sol-S-O}}
        {This figure demonstrates the optimal price and number of servers under the single service portfolio with the standard service (left panel) and the on-demand service (right panel). In each figure, we use a blue square, red plus, and green diamond to represent the service value at $0.8$, $1$, and $1.2$, respectively. 
        Parameters are $\Lambda = 15, \mu_o=1, \mu_s=1$.}
       \vspace{-0.2in}
\end{figure}

Our analysis reveals an important interaction between the supply pool and the value of service on the optimal pricing and per-service wage decisions of the SP. First, when the supply pool of potential self-scheduled agents or the value of service is large, the SP lowers the per-service wage and thus has a higher potential profit margin to provide service; therefore, it is optimal for the SP to attract and serve the entire market by setting a low price. Second, when the supply pool of contractors is small, it is more expensive for the SP to attract contractors to provide services, which, in turn, increases the lead time and thus demotivates customers to join the service. The above two propositions on the characterization of optimal pricing and staffing are visualized in Figure \ref{fig-sol-S-O}. 

Interestingly, we find that the optimal pricing and staffing are qualitatively consistent between these two different service mechanisms. 
First, the impacts of the supply-side and demand-side characteristics are consistent in the pricing and staffing. In particular, a higher value of service increases the price of service and the attractiveness of the service by reducing waiting time through increasing the number of servers, and a lower cost of labor, either because of a lower hourly wage of employees or a larger supply pool of contractors.
Second, the price is set to reach a targeted fraction of customers under each mechanism, and the optimal number of servers (employees or contractors) adopts a similar newsvendor structure that includes both the mean and safety capacity. For example, when $V\ge \overline{C}_s$ or $K\ge \bar K$, the price is set at $V-W_s$ under the standard mechanism or $V-W_o$ under the on-demand mechanism so that the most impatient customer joins the service under each mechanism, and the square term in the optimal staffing, given in \eqref{equ-staffing-single} in Corollary \ref{cor-staff}, reflects the safety capacity.  
Thus, the cost per service (i.e., $\overline{C}_s$ or $\overline{C}_o$) includes the markdown of the price and the staffing cost of serving in each mechanism. 

\begin{cor} {\sc (Optimal Staffing: Employees vs Contractors)} \label{cor-staff}
    Under the single service model, the optimal number of employees or contractors satisfies:
    \bea
    k_s^S = \frac{\lambda_s^S}{\mu_s} + \sqrt{\frac{\lambda_s^S}{\mu_s w_s}}, \quad \mbox{ or, } \quad k_o^O = \frac{\lambda_o^O}{\mu_o} + \sqrt{\frac{\lambda_o^O}{\mu_o \cdot \frac{w_o^O\lambda_o^O}{k_o^O} \cdot \max\left\{ \frac{2\bar K}{K},2\right\} }}, \label{equ-staffing-single}
    \eea
    where $\frac{w_o^O\lambda_o^O}{k_o^O}$ represents the hourly wage per contractor.
\end{cor}

Our results also indicate that the per-service payment mechanism with contractors may differ fundamentally from the per-time payment mechanism with employees in terms of customer targeting. Specifically, the SP may shut down the standard service by setting a sufficiently high price when it is not profitable to operate, while he always operates the on-demand service targeting a strictly positive fraction of customers that is proportional to the supply pool since contractors are not paid when they are not working. In this case, the SP sets a price to target only patient customers (i.e., those with a small cost of waiting) who join the service while the impatient customers balk at the service. In addition, under the same hourly wage and the same demand, according to \eqref{equ-staffing-single}, since $\max\left\{ \frac{2\bar K}{K},2\right\} \ge 2$, the service mechanism with contractors requires a smaller safety capacity and thus is more efficient than the mechanism with employees.





\begin{table}[bht] \small
    \vspace{-0.1in}
    \begin{center}
    {
    \caption{Optimal prices and staffing policy: A summary \label{ec-table-opt-solutions}}
        \begin{tabular}{l  c} \hline \addlinespace[0.5em]
        Service portfolio & service design $(p_i, k_i, \lambda_i, W_i)$ \\ \hline  \addlinespace[0.5em]
        \begin{tabular}[c]{@{}c@{}}Standard:\end{tabular} 
         & \begin{tabular}[c]{@{}c@{}} $\left\{\begin{array}{lllll}
         p_s^S = V -\sqrt{\frac{w_s}{\Lambda \mu_s}}, & k_s^S = \frac{\Lambda}{\mu_s} + \sqrt{\frac{\Lambda}{\mu_s w_s}}, & \lambda_s^S =  \Lambda, & W_s^S = \sqrt{\frac{w_s}{\Lambda \mu_s}}, & \mbox{if $V\ge \overline{C}_s$} \\ 
         p_s^S = V, &k_s^S = 0, &\lambda_s^S =  0, &W_s^S = 0, & \mbox{if $V< \overline{C}_s$}
         \end{array}\right.$
         \end{tabular}  \\ \addlinespace[0.5em]  
        \begin{tabular}[c]{@{}c@{}}On-demand:\end{tabular} & 
        \begin{tabular}[c]{@{}c@{}} 
        $\left\{\begin{array}{lllll}
         p_o^O = V - \frac{\bar L_o}{\Lambda}, &k_o^O = \frac{\Lambda (1 +\bar L_o)}{\mu_o \bar L_o}, &\lambda_o^O = \Lambda, &W_o^O = \frac{\bar L_o}{\Lambda}, & \mbox{if $K \ge \bar K$} \\ 
         p_o^O = V - \frac{\sqrt{V\Lambda+1}-1}{\Lambda}, &k_o^O = \frac{K\Lambda \sqrt{V\Lambda+1}}{\bar K \mu_o (\sqrt{V\Lambda+1}-1)}, &\lambda_o^O =  \frac{K}{\bar K} \Lambda, &W_o^O = \frac{\bar K (\sqrt{V\Lambda+1}-1)}{K \Lambda}, & \mbox{if $K<\bar K$}
         \end{array}\right.$
        \end{tabular} \\ \addlinespace[0.5em]  
        \hline 
        \multicolumn{2}{l}{{\it Note.} The wage per-service of contractors is $w_o^O = \frac{\left(k_o^O\right)^2}{K\lambda_o^O}$ by \eqref{ko}. }
        \end{tabular}
    }
    \end{center}
    \vspace{-0.2in}
\end{table}

\subsection{Standard vs On-demand}
In this section, we examine the implications of service mechanisms from various perspectives, including the SP profit, labor welfare, consumer surplus, and total social welfare. The optimal solutions derived in the previous two subsections are summarized in Table \ref{ec-table-opt-solutions}. To avoid the trivial case, in what follows, we assume $V\ge \overline{C}_s$ so that the SP will not shut down the standard service. 


We first examine the impacts of service mechanism on the SP profit. Denote $R(w_s , K)$ as the performance ratio of the SP profit under the System S relative to System O, i.e., $R(w_s , K) = \pi^S/\pi^O$, the following theorem provides tight bounds on this performance ratio. 

\begin{thm}\label{thm-compare-SO}
{\sc (SP Profit)}   Compared to the standard service, the on-demand service improves the profit of the SP if and only if the supply pool is large, i.e., $R(w_s,K)\le 1$ if and only if $ K \geq K_F(w_s) $, where $K_F(w_s)$ decreases in $w_s$.  Besides, the performance ratio is bounded by: 
\[ 1 - \frac{\overline{C}_s}{V}
\leq R(w_s, K) \leq   \frac{V}{V - \overline{C}_o}  .\]
\end{thm}

The above theorem reveals that the preference of the SP over employees and contractors is not universal. The SP is better off with contractors if and only if the supply pool of self-scheduled agents or the hourly wage of employees is large. Intuitively speaking, a larger supply pool of contractors or a higher hourly wage of employees increases the attractiveness of the service mechanism with contractors, due to a lower per-service wage of contractors or higher labor cost of employees. To be more precise, the key impacts of these two parameters on the SP are fundamentally different: the supply pool affects the SP mostly through the effective arrival rate under the on-demand model, while the fundamental driver of the hourly wage of employees is the profit margin per service under the standard model. This is because, as illustrated in Figure \ref{fig-sol-S-O}, the profit margin per service is relatively stable under the on-demand model and the demand is stable (i.e., all customers) under the standard model. Therefore, an increase in the supply pool benefits the SP profit via a higher demand under the on-demand model, and a larger hourly wage of employees hurts the SP due to a lower profit margin per-service under the standard model; both of these favor the on-demand model.

\begin{figure}[!hbt]
	\vspace{-0.1in}
	\FIGURE
	{   
             \begin{minipage}{0.8\textwidth}
                    \includegraphics[width=0.48\textwidth]{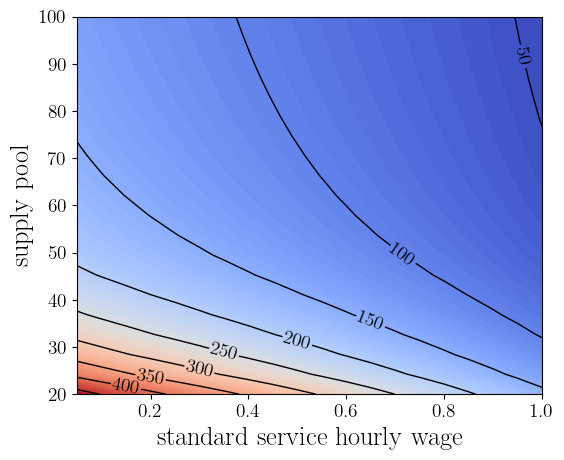}
				\includegraphics[width=0.45\textwidth]{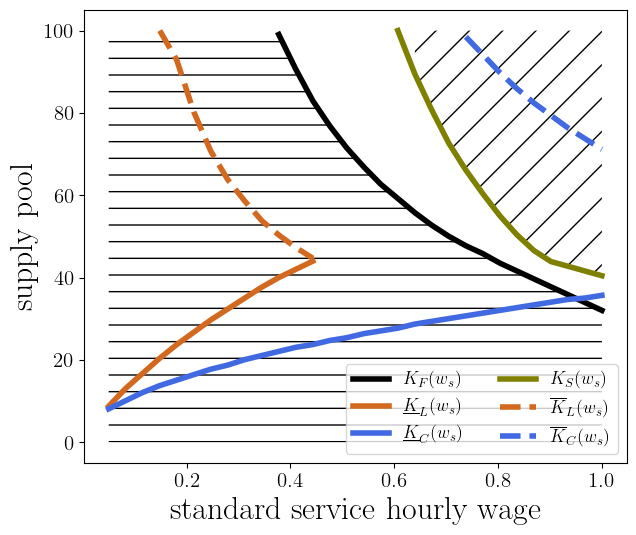}
		\end{minipage}
	}
	{The interplay between hourly wage and supply pool on the effectiveness of service mechanism \label{fig-ratio}}
	{This figure illustrates the impacts of the hourly wage of employees and the supply pool of contractors on the effectiveness of standard mechanism relative to on-demand mechanism. The left panel plots a heat map for the profit ratio of the SP (i.e., $R(w_s,K)\times 100\%$), in which a value of $<(>) 100$ means that the standard mechanism reduces (improves) the SP profit. The right panel depicts the regions in which the standard mechanism outperforms the on-demand mechanism from various perspectives, e.g., the region where the SP profit and the social welfare are both higher under the standard (or on-demand) mechanism is highlighted in horizontal (or diagonal) lines. 
}
	\vspace{-0.2in}
\end{figure}

We also present a tight lower bound and upper bound on the effectiveness of the supply mechanism with employees relative to contractors. These two bounds are tight in the sense that, for $\epsilon>0$, there exist scenarios $(w_s^\epsilon, K)$ and  $(w_s, K^\epsilon)$ such that the performance ratio $R(w_s^\epsilon, K) <  1 - \frac{\overline{C}_s^\epsilon}{V} +\epsilon$ and $R(w_s, K^\epsilon) > \frac{V}{V - \overline{C}_o^\epsilon}-\epsilon$, where $\overline{C}_s^\epsilon$ and $ \overline{C}_o^\epsilon$ are defined in \eqref{eq-barK-barLo}.
The lower bound is constructed by evaluating at the maximum SP profit over all possible supply pools under the on-demand model; thus, it only depends on the hourly wage of employees. In contrast, the upper bound is constructed via the maximum value over feasible hourly wages under the standard model; thus, it only depends on the supply pool of contractors via the queue length. As we explained below Proposition \ref{S-only}, $\overline{C}_s$ measures the cost of staffing per service under the standard model; thus, the lower bound represents the potential profit margin per service under the standard model, implying that the standard model with a higher profit margin is more effective. 
We remark that these two service mechanisms are approximately the same for the SP with a large service value, as these two bounds converge to $1$ when the service value becomes large.  The interplay between hourly wage of the employees and supply pool on the effectiveness is visualized in a heat map plot in Figure \ref{fig-ratio}, in which the parameters are $\Lambda = 30, V = 2, \mu_s = 1, \mu_o = 1$, and we use the same parameters in the following examples unless explicitly stated.

Next, we turn to investigate the impacts of service mechanisms on labor welfare, consumer surplus and total social welfare. We show in the following theorem that there is no universal dominance between the standard and on-demand mechanisms from the perspective of labor welfare, consumer surplus, and total social welfare, and provide sufficient and necessary conditions under which the standard dominates the on-demand mechanism.  

\begin{thm} \label{thm-compare-LW-CS}
{\sc (Labor Welfare, Consumer Surplus, and Social Welfare)} Compared to the standard mechanism with employees, the on-demand mechanism with contractors strictly improves  labor welfare, consumer surplus, and total social welfare if and only if $ \underline{K}_L(w_s)< K < \overline{K}_L(w_s)$, $\underline{K}_C(w_s) < K< \overline{K}_C(w_s) $, and $K> K_S(w_s)$, respectively.  Besides, $\underline{K}_L(w_s) \le \overline{K}_L(w_s)$ and $\underline{K}_C(w_s) \le \overline{K}_C(w_s)$ for any $w_s>0$, with the following monotone properties: $\underline{K}_L(w_s)$ is quasi-concave, $\overline{K}_L(w_s)$ is quasi-convex, $\underline{K}_C(w_s) $ increases, $ \overline{K}_C(w_s)$ and $K_S(w_s)$ decrease in $w_s$; see Figure \ref{fig-ratio} for an illustration.
\end{thm}

Our results above provide several intriguing implications of the service mechanism from various contexts. First, from the perspective of labor welfare, we find that the on-demand mechanism improves labor welfare if and only if the supply pool of self-scheduled agents is moderate. This is because, under the on-demand mechanism, the supply pool affects the total labor welfare in two ways: the number of contractors and the welfare per contractor. On one hand, a large on-demand supply pool intensifies competition in the labor market thus reducing the welfare of each contractor, via a lower per-service wage, in providing the service. On the other hand, a small supply pool of contractors reduces total labor welfare even if the welfare per contractor is improved. Therefore, as the supply pool of contractors increases, the total labor welfare under on-demand mechanism first increases and then decreases. This implies that, compared to the labor welfare under the standard mechanism (which does not depend on the supply pool), the on-demand mechanism either reduces labor welfare for any supply pool, which happens with a large hourly wage of employees, or improves the labor welfare with a moderate supply pool. 


Second, like labor welfare, consumers may prefer either the standard or the on-demand mechanism depending on the supply pool of contractors and the hourly wage of employees. Particularly,  consumers are better off under the on-demand mechanism if and only if the supply pool of contractors is moderate. The underlying intuition is that a small supply pool hurts the consumers due to a longer waiting time resulting from a smaller number of participating contractors, while a large supply pool also hurts consumers due to the increased price of services. Intuitively, one would expect that a large hourly wage of employees should reduce the number of employees and increase the waiting time, thus, hurting the consumer surplus under the standard service, which makes the on-demand mechanism more attractive from the perspective of consumer surplus. Counter-intuitively, our results reveal an opposite result: a large hourly wage makes the on-demand mechanism less attractive to the consumers. This is because, with a larger hourly wage of employees under the standard mechanism, the SP lowers the price of standard service to subsidize consumers for the longer waiting time, which improves consumer surplus.

Our analysis also reveals the impacts of the service mechanism on total social welfare that include SP profit, labor welfare, and consumer surplus.  Consistent with SP profit, we find that the on-demand mechanism improves total social welfare if and only if the supply pool of contractors is large. That is, though a large supply pool hurts self-scheduled agents and consumers, it increases the market demand, thus improving both SP profit and efficiency with a higher social welfare. Therefore, based on the value of the supply pool and two thresholds $K_F(w_s),K_S(w_s)$, compared to the on-demand mechanism with contractors, the standard mechanism improves, or reduces,  both SP profit and social welfare when the supply pool is small (i.e., $K\le \min\{K_F(w_s), K_S(w_s) \}$), or large (i.e. $ K \geq \max\{K_F(w_s),K_S(w_s) \}$), respectively, as visually highlighted in Figure  \ref{fig-ratio}.  

\begin{cor}\label{cor-joint-compare}
{\sc (Incentive Coordination)} The incentive coordination is achieved under the on-demand, or the standard mechanism, if and only if $ \max\{K_F(w_s),K_S(w_s), \underline{K}_L(w_s),\underline{K}_C(w_s)\} \leq K \leq \min\{ \overline{K}_L(w_s),\overline{K}_C(w_s)\}$, or $K \leq \min\{K_F(w_s),K_S(w_s), \underline{K}_L(w_s),\underline{K}_C(w_s)\}$, respectively. 
\end{cor}

Finally, we investigate the incentive coordination in which all parties, including the SP, labors, and consumers, are better off under one specific service mechanism. Our results above  indicate the interplay between the supply pool and hourly wage on coordinating incentives of a specific service mechanism. As visualized in Figure \ref{fig-compare-2}, we find that providing either the on-demand service or the standard service can align the incentives of all parties involved. Specifically, the standard service mechanism achieves coordination when the supply pool of contractors is small, while the on-demand service mechanism aligns the incentives of all parties when the supply pool is moderate; otherwise, coordination is impossible with a large supply pool. Thus, our finding highlights the key roles of the supply pool in determining whether or not coordination is possible with the standard or on-demand service mechanisms. On one hand, our findings suggest that a large supply pool will hurt labor welfare and consumer surplus under the on-demand mechanism with contractors; thus, providing theoretical support for the regulation of limiting the number of qualified and registered independent agents. On the other hand, our analysis also suggests that an insufficient supply pool of self-scheduled agents could result in inefficiency under the on-demand mechanism with contractors, reducing both SP profit and social welfare.  

\section{Portfolio of Differentiated Services: Hybrid Model}\label{sec:Hybrid}

In this section, we consider the hybrid model in which the SP opts to include both the standard and the on-demand services. Recall that, the standard service is provided by full-time employees and the on-demand service is served by part-time contractors. Our objectives in this section are to first provide a full characterization of the optimal service deployment that specifies when the SP should offer one or two services, and then thoroughly examine the impacts and values of proliferating a single service alongside both services. 

\subsection{Optimal Service Deployment}
Under a hybrid model, the SP could operate both the standard service with full-time employees and the on-demand service with part-time contractors. The SP selects the prices of the two services, the number of employees, and the per-service wage for contractors to maximize his expected total profit, see the formulation \eqref{opt}. Similar to the single service mechanism, to ensure the profit function is well-behaved (e.g., concavity), we first perform a one-to-one variable transformation that transfers the decision variables from $(p_s,p_o,k_s,w_o)$ to $(\lambda_s,W_s,\lambda_o,W_o)$. We then sequentially solve a single-variable optimization with decision variables $\lambda_s,W_s,W_o$ and $\lambda_o$. We show that the optimization problem \eqref{opt} can be simplified by finding the optimal value within three sub-cases corresponding to two single services and one both services, i.e., $\pi^* = \max\{\pi^S, \pi^O,  \pi^{H}    \}$, where $\pi^S,\pi^O,\pi^H$ represents the optimal profit under the standard only service, the on-demand only service, and both services. Finally, we show the monotonic properties of $\pi^S, \pi^O,  \pi^{H}$ on $(w_s,K)$ and divide the feasible space of $(w_s,K)$ into three different regimes (i.e., regime $\mathcal{S},\mathcal{O},\mathcal{H}$) to provide a full characterization of the optimal service deployment policy within each regime. 
We formally present and establish these results in the following Theorem \ref{opt-two}, which is visualized in Figure \ref{fig-two-tier}.

\begin{thm} \label{opt-two}{\sc (Optimal Service Deployment)} The feasible space of the hourly wage and the supply pool could be divided into three regimes, $\mathcal{S},\mathcal{O},\mathcal{H}$, such that the optimal service deployment \eqref{opt} for the SP is to: 
\begin{itemize}
\item[{\bf Regime $\mathcal{S}$:}] Provide the standard service only with employees if $(w_s,K)\in \mathcal{S}$.
\item[{\bf Regime $\mathcal{O}$:}] Provide the on-demand service only with contractors if $(w_s,K)\in \mathcal{O}$.
\item[{\bf Regime $\mathcal{H}$:}] Provide two services with both employees and contractors if $(w_s,K)\in \mathcal{H}$.  
\end{itemize}
\begin{figure}[!hbt]%
	\vspace{-0.1in}
	\FIGURE
        {
            \includegraphics[width=0.75\textwidth]{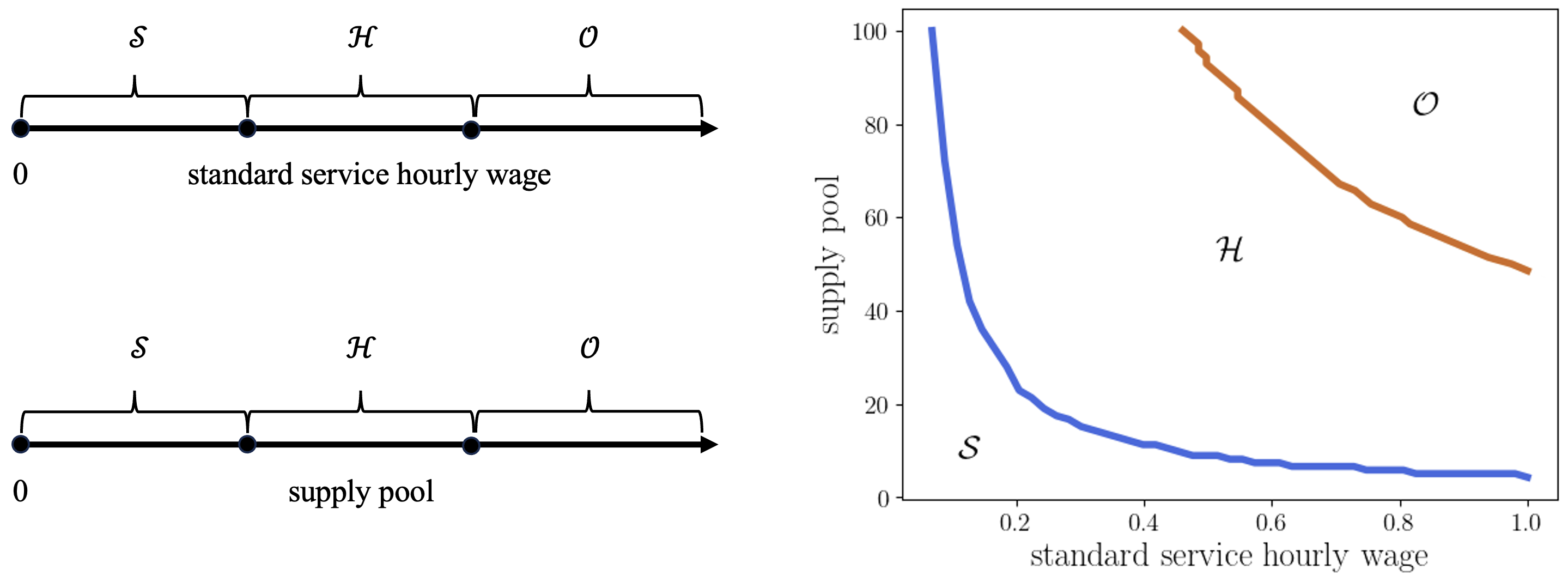} 	
        }
	{The visualization of the optimal service deployment \label{fig-two-tier} }
	{This figure visualizes the optimal service deployment characterized in Theorem \ref{opt-two}. The left panel illustrates the threshold structure for a given supply pool or fixed hourly wage, and the right panel demonstrates the partitions of the feasible space of $(w_s,K)$. 
 }
   \vspace{-0.2in}
\end{figure}
\end{thm}

Our results reveal that the SP may not be universally better off proliferating services that differentiate in prices, lead times, and service mechanisms. Instead, it is optimal for the SP to provide only one service when the labor cost of this service is significantly lower than that of the other one.  
Specifically, in Regime $\mathcal{S}$ with a small supply pool or a low hourly wage, the employees have a significant advantage over the contractors in terms of labor cost; thus, it is optimal for the SP to hire only employees to provide the standard service. On the other hand, when the contractors are much cheaper, due to either a sufficiently large supply pool or a high hourly wage of employees (i.e., Regime $\mathcal{O}$), it is more profitable for the SP to operate only the on-demand service, which is outsourced to self-scheduled agents with a per-service wage. 

We find that service proliferation benefits the SP only when the labor costs of these services are relatively comparable. That is, in Regime $\mathcal{H}$ with a moderate supply pool and a moderate hourly wage, a hybrid model with two differentiated services outperforms any single service deployment. Our analysis also provides two important implications under the hybrid deployment. First, the SP sets the prices to be attractive enough so that all customers will join one of these two services, i.e., the entire market is covered under the hybrid service model. Second, compared to the standard service, the on-demand service could be either more expensive with a shorter lead time, or cheaper with a longer lead time. This is because our base model considers a homogeneous quality between these two services without requiring a higher quality of service for the on-demand service. These two features indicate that the SP selects the pricing and staffing for these two complementary services so that the patient customers are served by the slower service and the impatient customers select the faster service. 

\subsection{Impacts of Service Proliferation on Prices and Lead times}
In this subsection, we explore the impacts of service proliferation that adds the other service to a one-tier system on the price and lead time. We first examine the impacts of adding the on-demand service to a single service system with only the standard service, and then study the impact of adding the standard service to the one-tier system that provides only the on-demand service. To that end, we focus on Regime $\mathcal{H}$ whereby the optimal service deployment is to operate both services.

 \begin{thm}\label{pro:compareTS}
 {\sc (Impacts on the Standard Service)}  
  Compared to System S, the hybrid model with two services reduces the price while increasing the lead time of the standard service.   
\end{thm}

The above results show how the price and lead time of the standard service are affected when the SP proliferates the standard service with the other service, which is outsourced to  independent contractors. This situation captures the practical observations in which many logistics companies, like SF Express, start to supplement their existing standard services with  on-demand services which are often outsourced. We find that, service proliferation reduces the price but simultaneously enlarges the lead time of the existing standard service, as introducing another service creates internal competition between the two services. Consequently, some impatient customers will switch from the standard service to the on-demand service when this service is available, while  patient customers will still remain in the standard service regardless of the availability of the other service.

Next, we consider another scenario whereby the SP is complementing the existing on-demand service with the standard service by hiring full-time employees. Recall that, as we discussed after Theorem \ref{opt-two}, in the hybrid model, the lead time under the on-demand service could be either shorter or longer than that under the standard service. Specifically, we could further divide the Regime $\mathcal{H}$ into two disjoint regimes,  $\mathcal{H}_1$ and $\mathcal{H}_2$ such that the on-demand service has a longer (or shorter) lead time when $(w_s,K)\in\mathcal{H}_1$ (or $\mathcal{H}_2$). Our next theorem shows that the impacts of service proliferation on the price and lead time of the existing on-demand service depend on whether or not the on-demand service has a shorter lead time in the two-service hybrid system. 

\begin{thm}\label{pro:compareTO}
{\sc (Impacts on the On-demand Service)} Compared to System O, the hybrid model with two services will:
{\rm (a)} reduce the price while enlarging the lead time of the on-demand service if it has a longer lead time (i.e., Regime $\mathcal{H}_1$) and  $(w_s, K) \in \mathcal{T}_1$;
{\rm (b)} increase the price but shorten the lead time of the on-demand service if it has a shorter lead time (i.e., Regime $\mathcal{H}_2$) and  $(w_s, K) \in \mathcal{T}_2$, where $\mathcal{T}_1$ and $\mathcal{T}_2$ are defined in the proof. 
\end{thm}

\begin{figure}[!hbt]%
	\vspace{-0.1in}
	\FIGURE
        {
            \begin{minipage}{0.9\textwidth}
                    \includegraphics[width=0.32\textwidth]{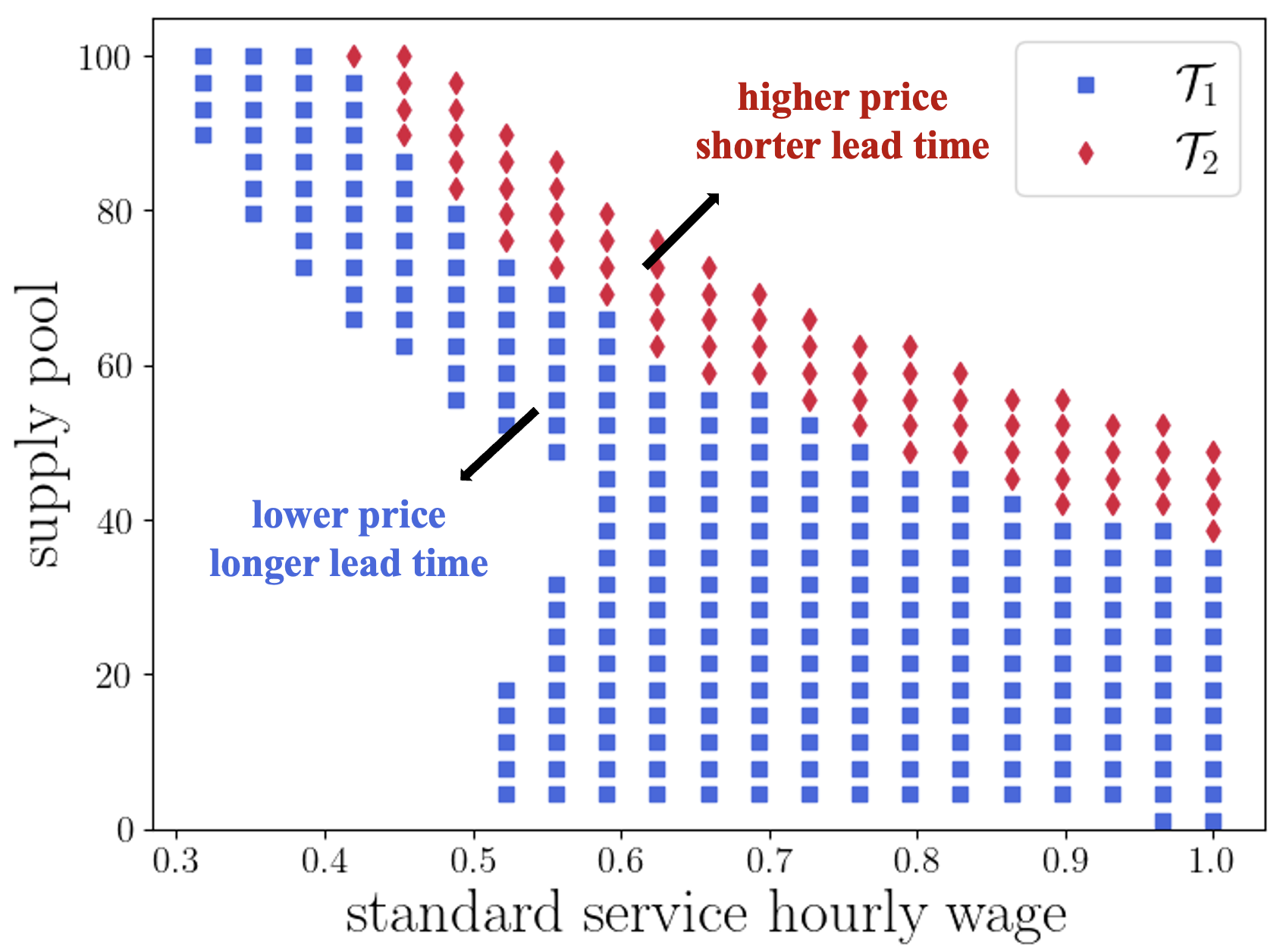}
				\includegraphics[width=0.64\textwidth]{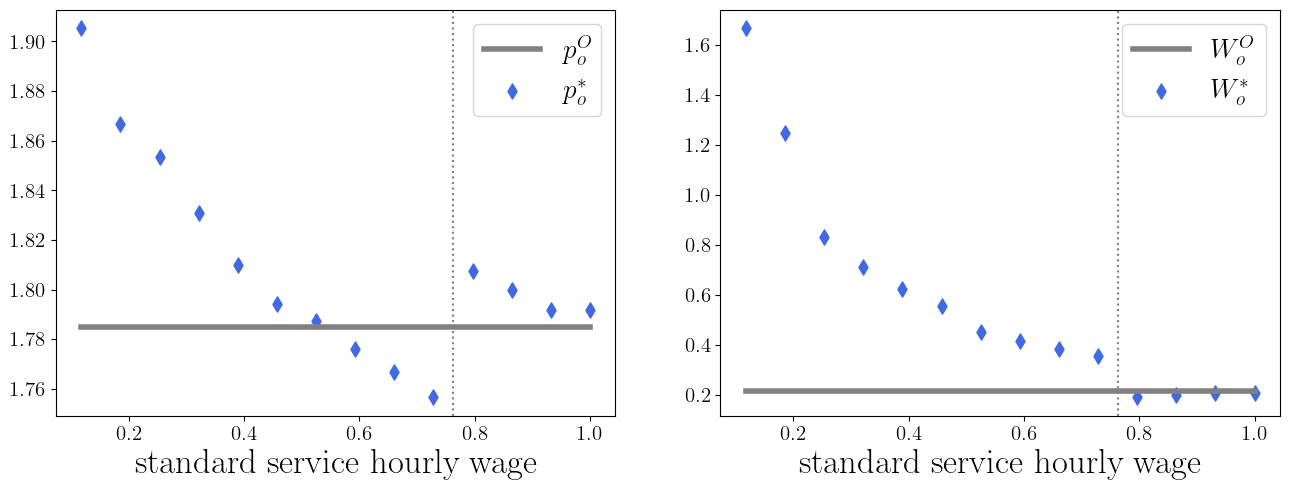}
		\end{minipage}
        }
	  {Impacts of service proliferation on the price and lead time of the on-demand service \label{fig-price-o-t}}
	{This figure shows that proliferating the on-demand service with the standard service may either shorten or lengthen the lead time of the on-demand service. The left panel illustrates the regimes under which proliferation enlarges ($\mathcal{T}_1$ with blue-square dots) or shortens ($\mathcal{T}_2$ with red-diamond dots) the lead time of on-demand service. The right panels use a specific example with $K=50$ to demonstrate the impacts of hourly wage on the prices and lead times, in which, the black-solid line, blue-diamond dots, and vertical black-dotted line represent metrics under system O, metrics under the hybrid model, and the boundary between the two regimes $\mathcal{H}_1$ and $\mathcal{H}_2$, respectively. 
 }
   \vspace{-0.2in}
\end{figure}



Unlike the standard service, the above results indicate a mixed effect of service proliferation on the price and lead time of the on-demand service. Specifically, as visualized in Figure \ref{fig-price-o-t}, adding the standard service could either shorten or lengthen the lead time of the on-demand service, depending on whether or not the on-demand service has a shorter lead time in the hybrid system. 
Note that System O may not serve all customers, 
while the entire market will be covered after adding the standard service in the hybrid system. Namely, adding the standard service has two effects: the {\it competition effect} and {\it market-expansion effect}. The competition effect generally reduces the price of the on-demand service, which in turn attracts more customers thus negatively reducing its attractiveness due to the enlarged lead time. On the other hand, the SP jointly adjusts the prices of these two services so that all customers, including those impatient customers that would balk when the standard service is not available, are served. 
Generally speaking, when the on-demand service is less (or more) cost-effective in Regime $\mathcal{H}_1$ (or $\mathcal{H}_2$), introducing the standard service reduces (or increases) the price thus enlarging (or shortening) the lead time of the on-demand service to target serving the patient (or impatient) customers, with the remaining customers served by the newly introduced standard service.

\subsection{Value of Service Proliferation} \label{subsec-value-OS}
In this subsection, we examine the economic value of service proliferation, which supplements one service with the other service, in the context of the SP profit, labor welfare, consumer surplus and total social welfare. We first investigate the value of service proliferation on the profit of the SP. We define the value of the on-demand service and the value of the standard service to the SP as: 
$
\Delta_o  = \frac{\pi^* - \pi^S}{\pi^*}  = 1 - \frac{\pi^S}{\pi^*} \mbox{ and } \Delta_s  = \frac{\pi^* - \pi^O}{\pi^*} = 1 -\frac{\pi^O}{\pi^*},
$
which measures the relative improvement of profit due to adding the on-demand (or standard) service to System S (or System O) that only provides the standard (or on-demand) service.

\begin{thm} {\sc (Value of Service Proliferation)}\label{thm-pi-K-ws}  
The value of the on-demand service increases in the supply pool $K$ and the value of the standard service decreases in the hourly wage $w_s$, with the following bounds:
\beq
0 \leq \Delta_o 
\leq \frac{\overline{C}_s}{V}, \quad \mbox{ and } \quad 
0 \leq \Delta_s 
\leq \frac{\overline{C}_o}{V}
\eeq
where $\overline{C}_s, \overline{C}_o$ are given in \eqref{eq-barK-barLo}. 
\end{thm}

\begin{figure}[!hbt]%
	\vspace{-0.1in}
	\FIGURE
         {
            \includegraphics[width=0.7\textwidth]{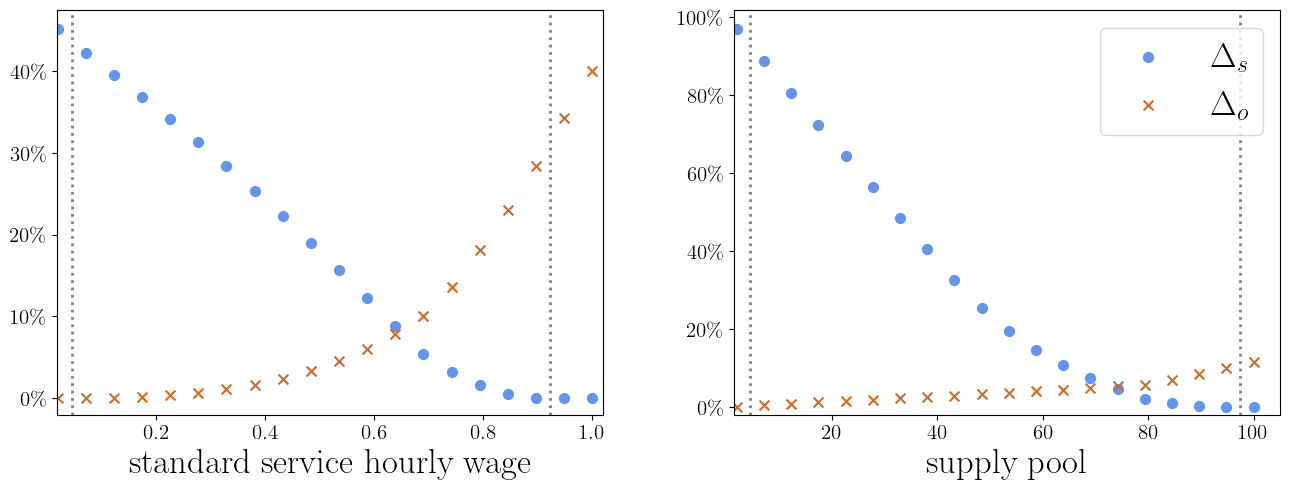} 
        }
	{\bf Value of service proliferation \label{20240820-1}}
	{This figure illustrates the impacts of employees' hourly wage $(w_s)$ and contractors' supply pool $(K)$ on the value of service proliferation to the SP profit. The two dotted vertical lines in each figure indicate the switching thresholds between two deployment policies (see Figure \ref{fig-two-tier}), and the parameters are $K = 55, \Lambda = 30, V = 2, \mu_s = 1, \mu_o = 1, w_s = 0.5$.}
\vspace{-0.2in}
\end{figure}

Theorem \ref{thm-pi-K-ws} derives the upper bounds on the value of service proliferation that adds the second service on the SP profit. We show that these two upper bounds are tight in the same sense as Theorem \ref{thm-compare-SO}, i.e., for any $\epsilon>0$, there exist scenarios $(w_s^\epsilon,K^\epsilon)$ such that $\Delta_o (w_s^\epsilon,K^\epsilon)> \frac{\overline{C}_s}{V} - \epsilon$ and $\Delta_s (w_s^\epsilon,K^\epsilon) > \frac{\overline{C}_o}{V} - \epsilon$. We remark that these upper bounds reflect the inefficiency associated with a single service mechanism, since $\overline{C}_s$ or $\overline{C}_o$ measures the per-service cost in the system with employees or contractors, respectively.  Therefore, the above results establish the intuition that proliferating with another service is more valuable only when the current system is less efficient. The monotonic properties of the value of service proliferation follow from the corresponding monotonicity of the inefficiency of a given service mechanism. For example, as we illustrate in Figure \ref{20240820-1}, with a smaller supply pool (or a higher hourly wage), the mechanism with contractors (or employees) is less efficient; thus, adding the standard (or on-demand) service could be more valuable. 
 


Next, we explore the value of service proliferation from various other perspectives, including the customers, the agents, and the entire society. Analogous to the SP profit, we define the following values of service proliferation:  
$
\Delta_o M 
= 1 - \frac{M^S}{M^*} \mbox{ and } \Delta_s M 
= 1 - \frac{M^O}{M^*}, \  M \in \{CS, LW, SW\}, 
$
where $CS^*,CS^O,CS^S$ are the consumer surplus under the hybrid system with both employees and contractors, the System O with contractors, and the System S with employees, respectively, $LW^*,LW^O,LW^S$ and $SW^*,SW^O,SW^S$ are the corresponding labor welfare and social welfare. The value of adding the on-demand (or standard) service is measured by $\Delta_o M$ (or $\Delta_s M$) from the perspective of performance measure $M$. For example, $\Delta_o CS$ quantifies the value of adding the on-demand service into the existing standard service from the customers' perspective.

\begin{figure}[htbp]
	\vspace{-0.1in}
	\FIGURE
	{
            \begin{minipage}{0.96\textwidth}
                \includegraphics[width=0.24\textwidth]{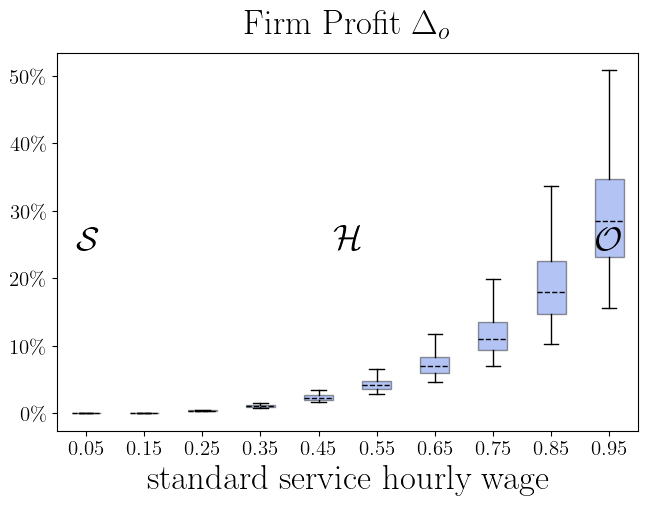} 
                \includegraphics[width=0.24\textwidth]{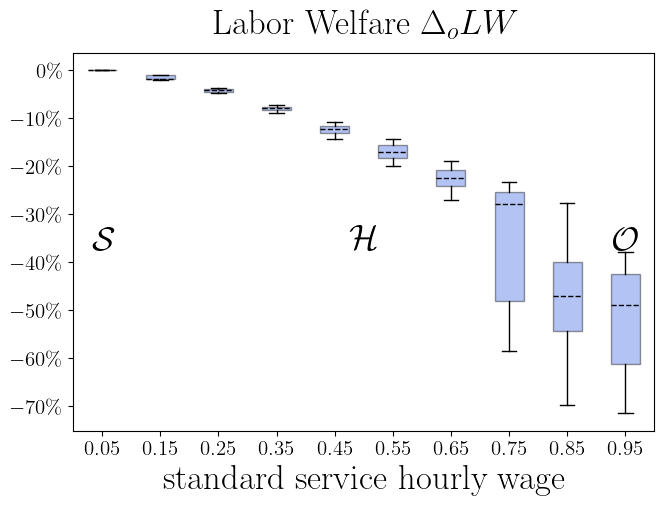} 
                \includegraphics[width=0.24\textwidth]{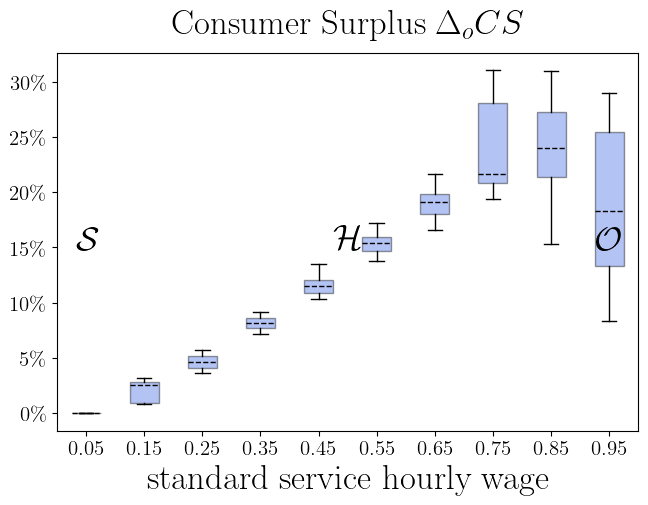} 
                \includegraphics[width=0.24\textwidth]{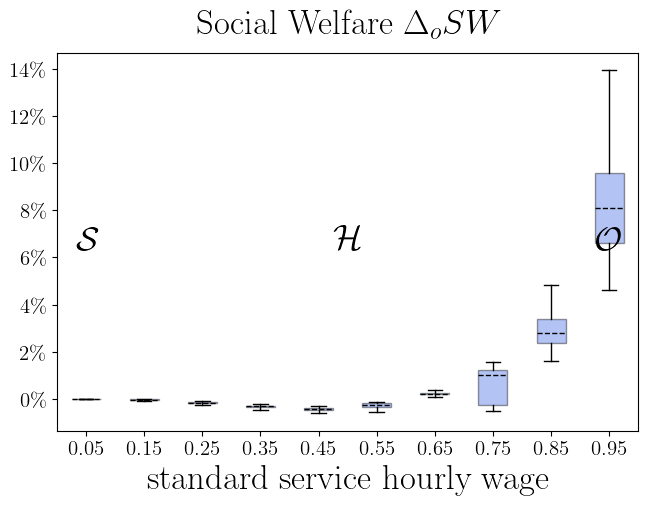} 
                \vfill
                \includegraphics[width=0.24\textwidth]{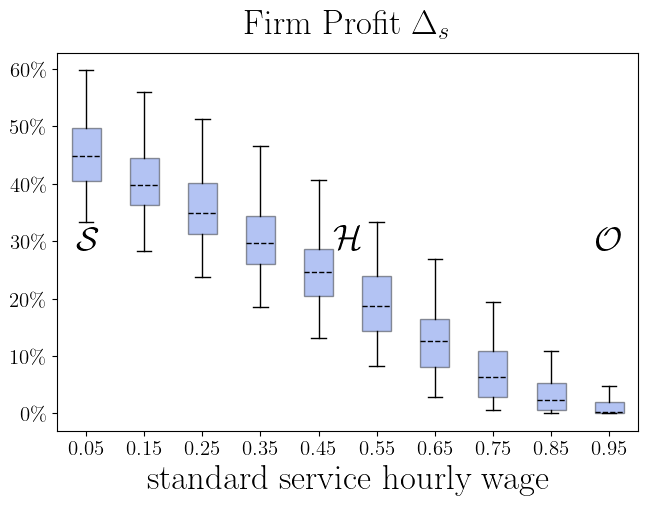} 
                \includegraphics[width=0.24\textwidth]{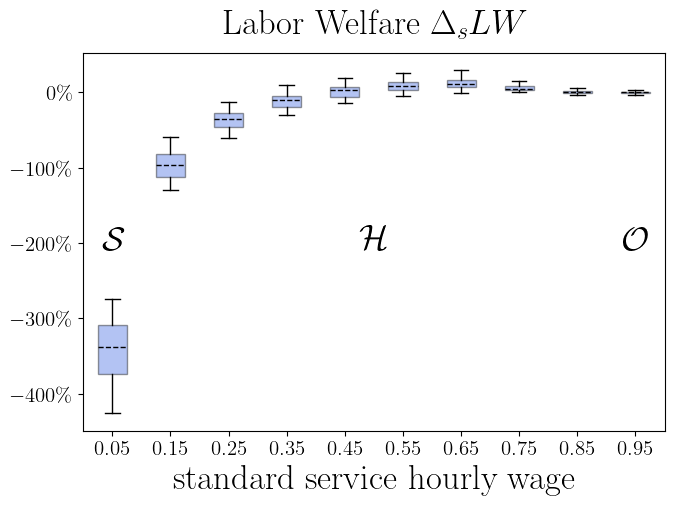} 
                \includegraphics[width=0.24\textwidth]{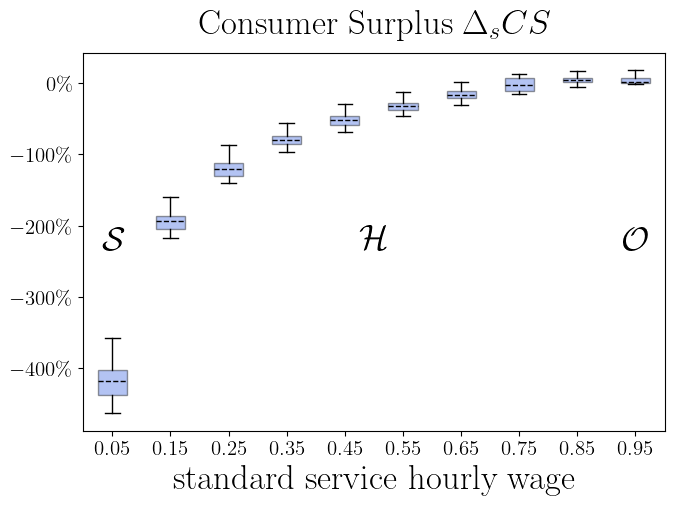} 
                \includegraphics[width=0.24\textwidth]{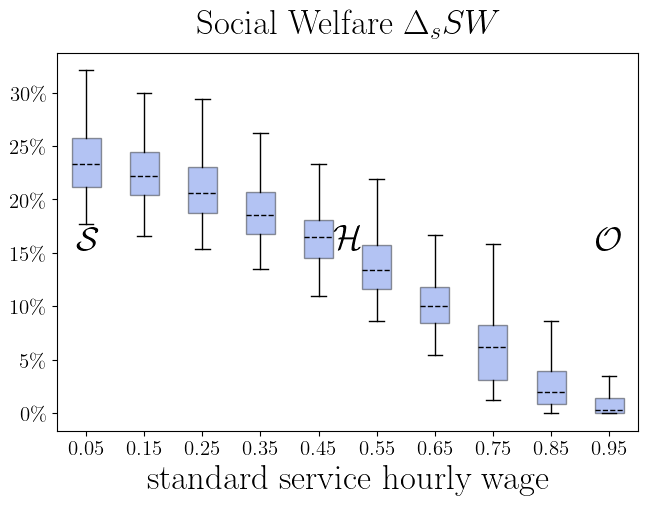}
            \end{minipage}
	}
	{\bf Value of service proliferation: Impacts of the hourly wage of employees. 
 \label{fig-value-base}}
	{
 This figure illustrates the impacts of employees' hourly wage $(w_s)$ on the value of service proliferation from the perspective of the firm, consumers, labors, and total social welfare. Specifically, for each $w_s$, 
 we randomly sample $100$ instances of service value $V$ and market size $\Lambda$, 
 and use a boxplot to report the summary statistics that include the minimum, first quantile, median, third quantile, and maximum.
 }
\vspace{-0.2in}
\end{figure}

\begin{figure}[htbp]
	\vspace{-0.1in}
	\FIGURE
	{
            \begin{minipage}{0.96\textwidth}
                \includegraphics[width=0.24\textwidth]{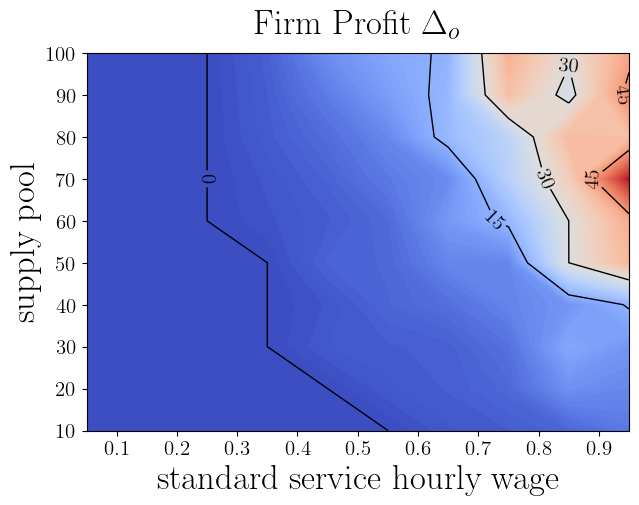} 
                \includegraphics[width=0.24\textwidth]{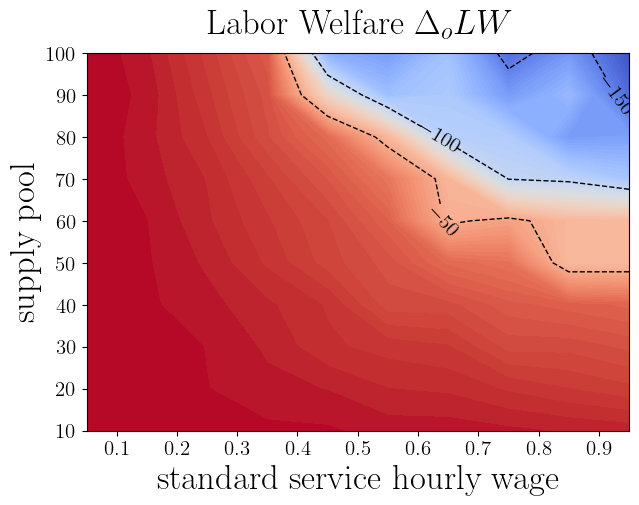} 
                \includegraphics[width=0.24\textwidth]{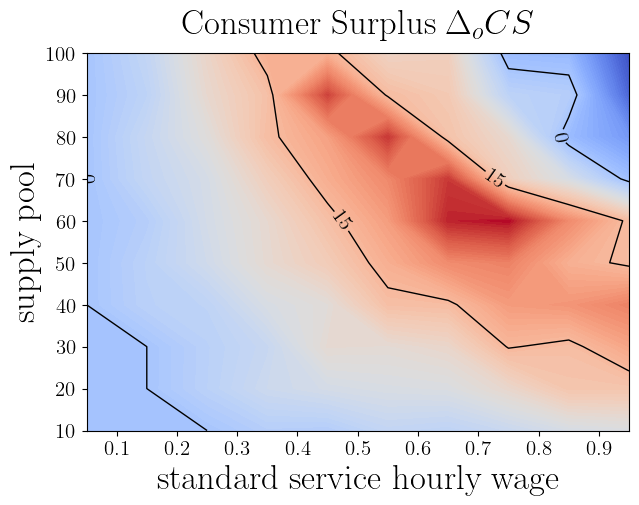} 
                \includegraphics[width=0.24\textwidth]{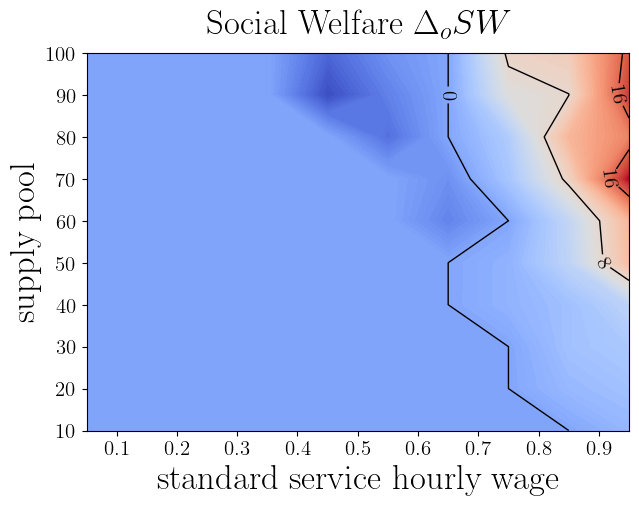} 
                \vfill
                \includegraphics[width=0.24\textwidth]{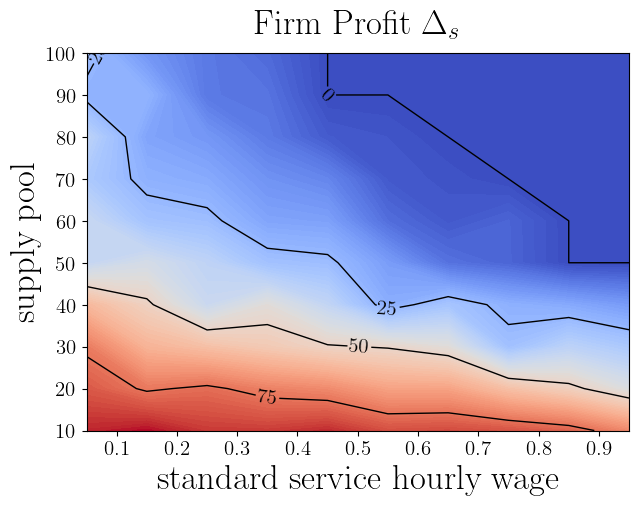} 
                \includegraphics[width=0.24\textwidth]{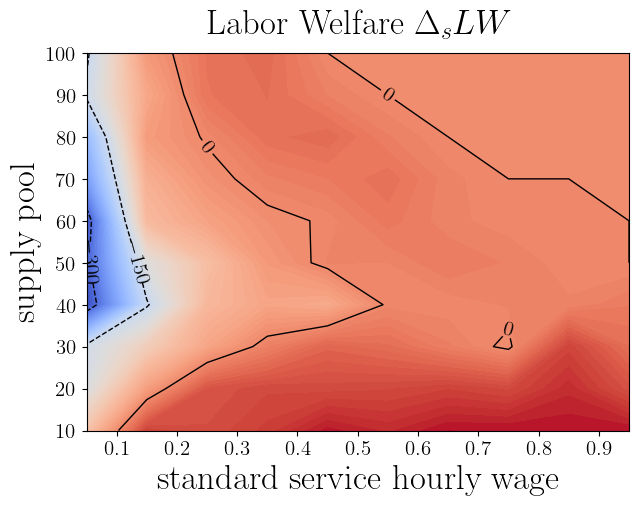} 
                \includegraphics[width=0.24\textwidth]{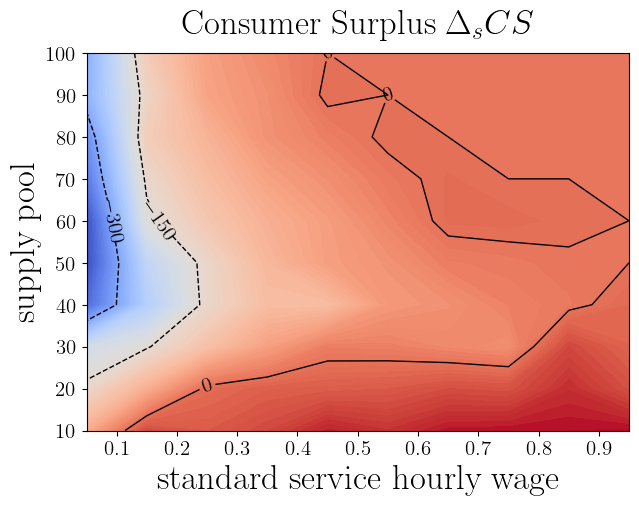} 
                \includegraphics[width=0.24\textwidth]{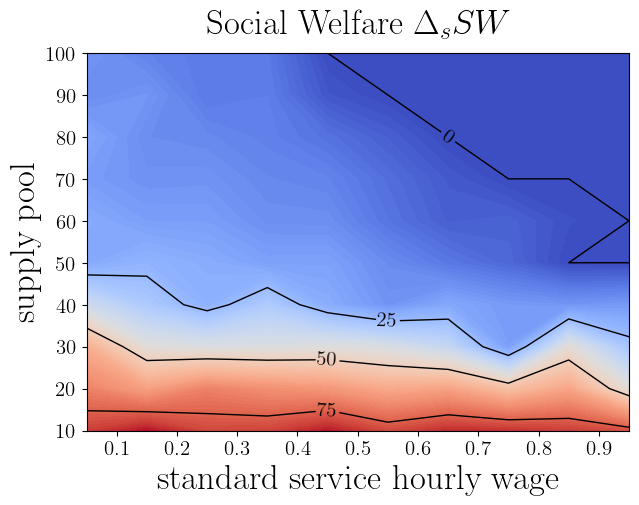}
            \end{minipage}
	}
	{Value of service proliferation: Interplay between employees and contractors \label{fig-value-heatmap}}
	{This figure demonstrates the interplay between employees' hourly wage $(w_s)$ and contractors' supply pool $(K)$ on the value of service proliferation from various perspectives. Specifically, for each pair of $(w_s,K)$, 
 we randomly sample $100$ instances of service value $V$ and market size $\Lambda$, and report the average value across these $100$ samples in the above heatmap.
 }
\vspace{-0.2in}
\end{figure}

In what follows, we will conduct extensive numerical simulations to thoroughly examine the value of service proliferation from various perspectives: consumer surplus, labor welfare, and social welfare. We remark that, as illustrated in Theorem \ref{prop-relation}, a similar analysis could be applied to characterize these values analytically. Recall that the reservation rate of contractors is normalized to the unit interval $[0,1]$. Unless explicitly stated, we set the parameters $
w_s \in (0,1],  K = 50, V=2, \mu_s=1, \mu_o=1.$ 
For the results to be general and representative, for each scenario, we randomly sample $100$ instances of service value $V \in [1.7,2.5]$ and customer arrival rate $\Lambda \in [25,35]$, and report key summary statistics among these $100$ instances in Figures \ref{fig-value-base} and \ref{fig-value-heatmap}.
We also demonstrate the robustness of our research findings when the sample size is increased from $100$ to $1000$; details are relegated to Figure \ref{fig-value-size} in the Appendix \ref{ec-sec-numerical}.

We next discuss the in-depth insights into the implications of service proliferation. We remark that all three different service deployments are possible within the range of parameters; see Figures \ref{fig-value-base} and \ref{fig-value-heatmap}, which justifies the rationality of the parameters selected. We observe that service proliferation, in general, benefits the SP but hurts labor welfare, and may either improve or reduce consumer surplus and social welfare. This is because the SP usually over-exploits the agents by optimizing his labor deployment between the employees and the contractors, thus hurting total labor welfare. On the other hand, customers generally benefit from adding the on-demand service while being hurt from adding the standard service, due to a lower or higher price induced by introducing the cost-efficient per-service payment, or the cost-inefficient per-time payment mechanism; see also Theorem \ref{pro:compareTS} or Theorem \ref{pro:compareTO}. Consequently, the net effects combined may improve or reduce social welfare. Specifically, service proliferation improves the social welfare only if the existing service is cost-inefficient and the new service is cost-efficient. We further show the robustness of these findings by changing the market size and supply pool, see details in Figures \ref{fig-Lambda-K-2} and \ref{fig-Lambda-K-1} in Online Appendix Section \ref{ec-sec-numerical}.   

\section{Extensions} \label{sec-extension}
In this section, we relax key assumptions in the base model to demonstrate the robustness of our research findings via extensive numerical simulations. We will first consider a multi-server queue for the lead time, then study general heterogeneity in customers and agents, and finally allow for heterogeneous quality between the standard and on-demand services. 

\subsection{An $M/M/k$ Queue} \label{subsec-multi}
For analytical tractability, the base model approximates the lead time using an $M/M/1$ queue. In this subsection, we relax this assumption and use an $M/M/k$ queue to derive the lead time of each service, i.e., the lead time is given by $W = \frac{1}{\mu } + \frac{1}{1+ \frac{c!(1-\rho)}{k^k \rho^k} \sum_{i=0}^{k-1}\frac{k^i \rho^i}{i!}} \left[ \frac{\rho}{\lambda (1-\rho)} \right]$, 
where $\lambda, \mu, c, \rho = \frac{\lambda }{c\mu}$ are the arrival rate, service rate per server, number of servers, and utilization, respectively.  

\begin{figure}[htbp]
	\vspace{-0.1in}
	\FIGURE
	{
            \begin{minipage}{0.9\textwidth}
                \includegraphics[width=0.32\textwidth]{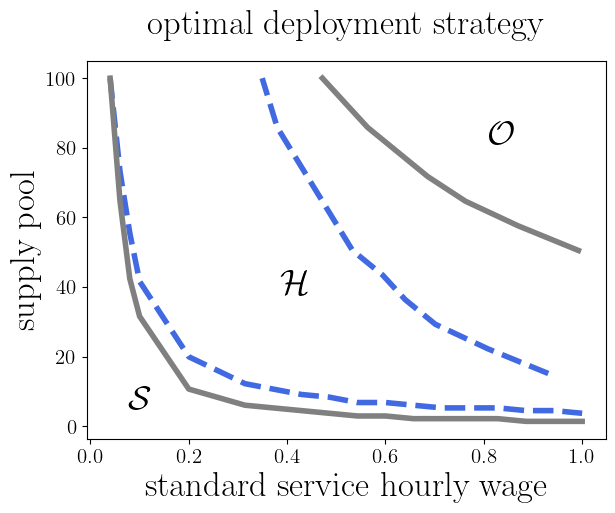} \hfill
                \includegraphics[width=0.32\textwidth]{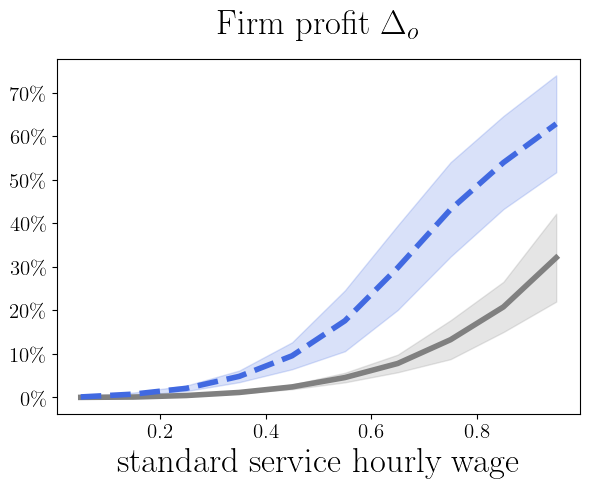} \hfill 
                \includegraphics[width=0.32\textwidth]{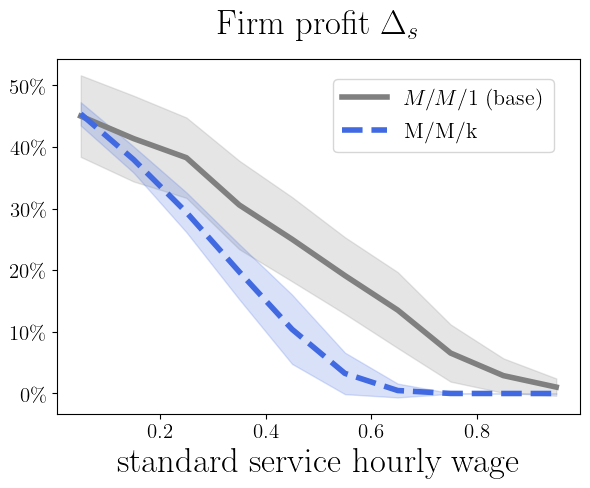}
            \end{minipage}
	}
	{\bf Optimal deployment and value of service proliferation: $M/M/1$ vs $M/M/k$ queue \label{fig-mmk}}
	{This figure illustrates the robustness of the research findings derived under the $M/M/1$ approximation. The optimal deployment strategy is shown in the left panel, in which the switching thresholds are represented by the gray solid lines for the base model with an $M/M/1$ queue, and the blue dashed lines for an $M/M/k$ queue. The value of adding the on-demand service or standard service to the firm is reported in the mid- or right panel, respectively, where the shaded area represents the interval that is within one standard deviation from the average value. }
	\vspace{-0.2in}
\end{figure}

Our analysis reveals that the findings derived under an $M/M/1$ approximation continue to hold qualitatively in an $M/M/k$ multi-server queue. Specifically, as we illustrate in Figure \ref{fig-mmk}, the $M/M/k$ queue does not change the structure of the optimal deployment policy given by Theorem \ref{opt-two} and the key insights regarding the value of service proliferation characterized in Theorem \ref{thm-pi-K-ws}, that is, adding another service mechanism is more valuable when the existing service is less cost-efficient and the newly added service is more efficient. In addition, we observe that the effectiveness of the standard mechanism relative to the on-demand mechanism, and the impacts of service proliferation on consumers, labors, and social welfare under the $M/M/k$ queue are also qualitatively similar to those under the base model, see details in Figure \ref{fig-mmk-lw} in Appendix \ref{ec-sec-numerical}. These observations indicate that approximating the lead time using an $M/M/k$ queue does not qualitatively change the key insights, while providing a tractable model to generate in-depth and comprehensive understanding of the impacts of service proliferation with different supply mechanisms.
 
\subsection{Heterogeneity of Customers and Agents}  \label{subsec-distribution}
Our base model adopts uniform distributions to capture the heterogeneity in the customers' waiting sensitivities and the agents' reservation rates. In this section, we relax this assumption and consider beta distributions, denoted by $\beta(a,b)$, with different parameters, that have the same support $[0,1]$, to investigate the effects of variability and skewness of the heterogeneity. Note that the uniform distribution $\mathcal{U}[0,1]$ is a special case of $\beta(1,1)$.  


\begin{figure}[htbp]
	\vspace{-0.1in}
    \FIGURE
	{
            \begin{minipage}{0.9\textwidth}
                \includegraphics[width=0.32\textwidth]{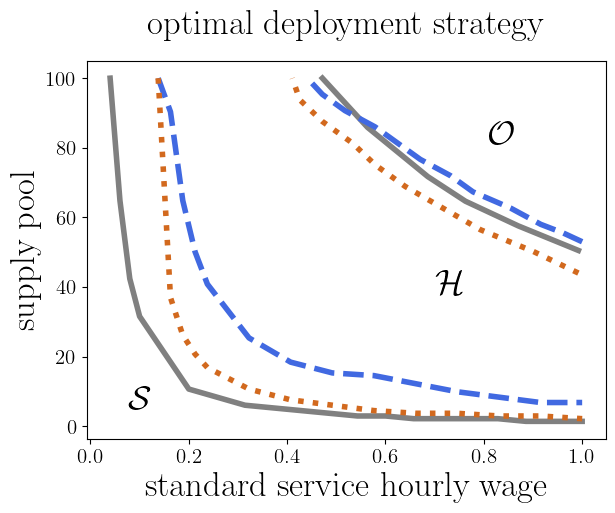}
                \includegraphics[width=0.32\textwidth]{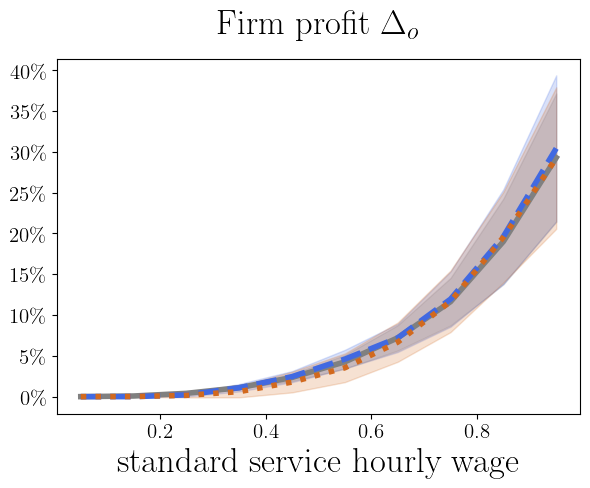}
                \includegraphics[width=0.32\textwidth]{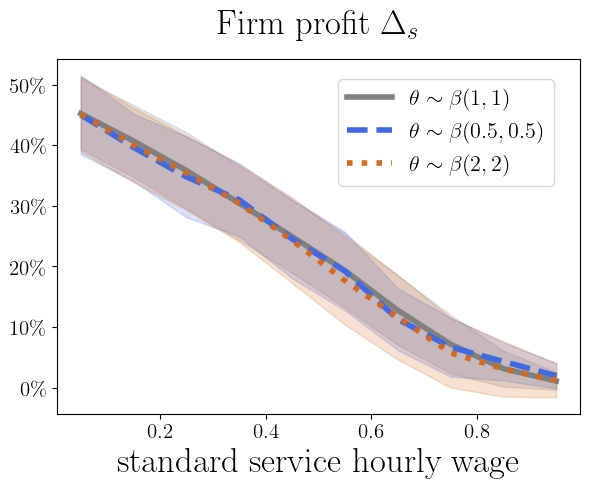}
                \vfill
                \includegraphics[width=0.32\textwidth]{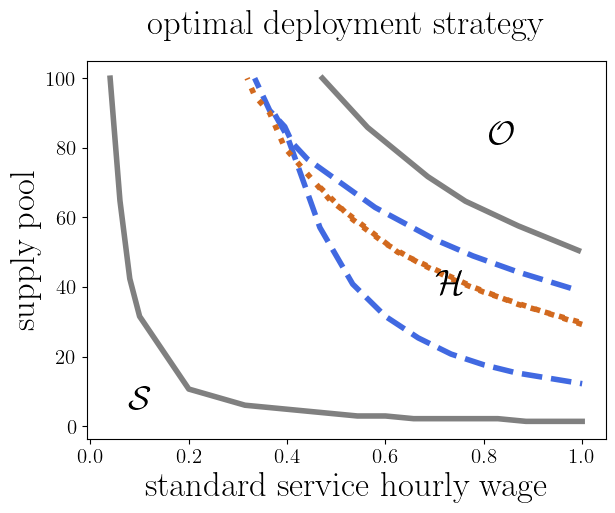} 
                \includegraphics[width=0.32\textwidth]{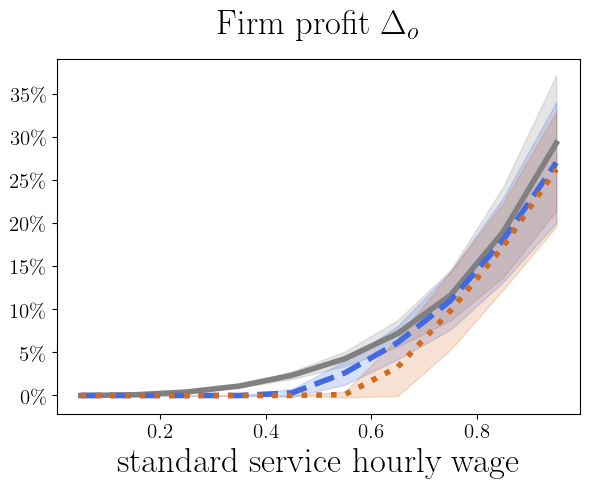}
                \includegraphics[width=0.32\textwidth]{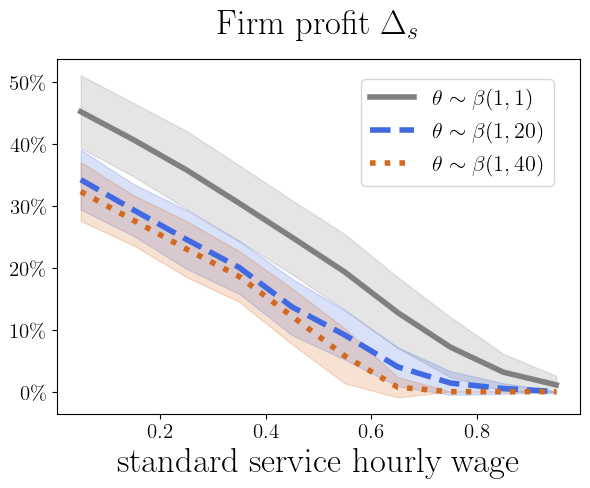}
                \vfill
                \includegraphics[width=0.32\textwidth]{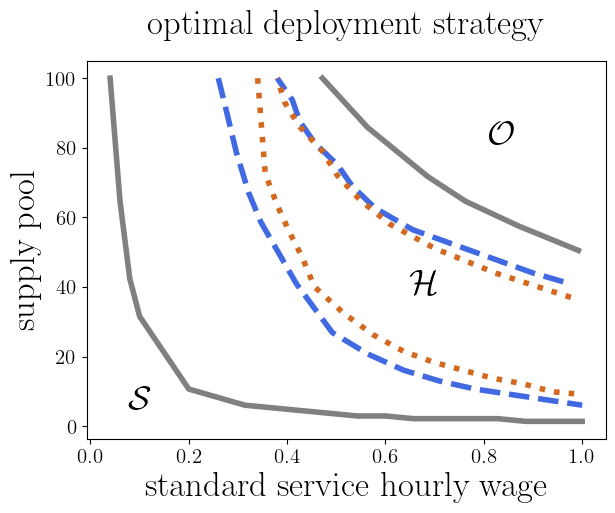} 
                \includegraphics[width=0.32\textwidth]{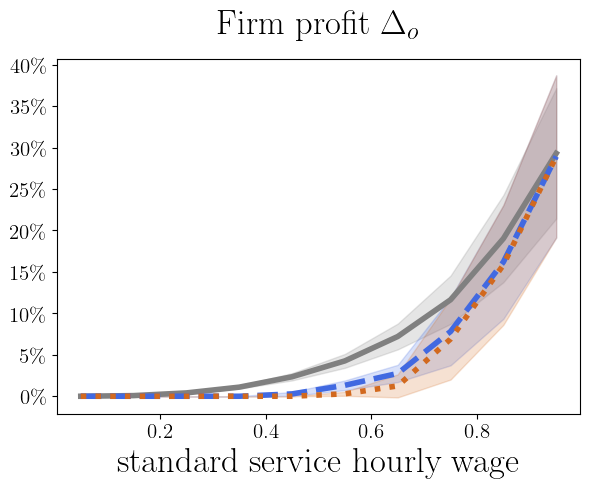}
                \includegraphics[width=0.32\textwidth]{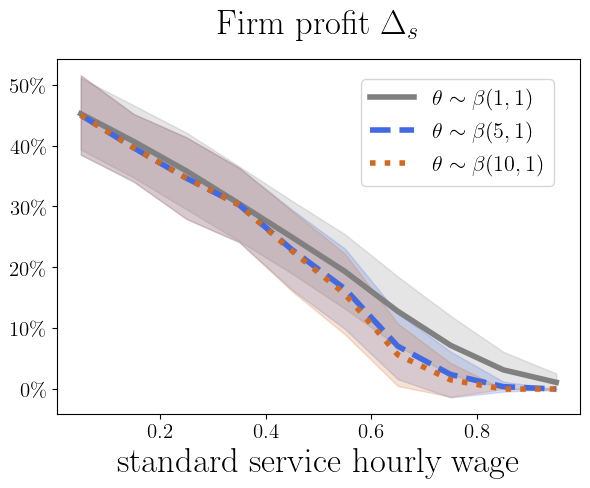}
            \end{minipage}
	}
        {Optimal deployment and value of service proliferation: Heterogeneous customers' waiting sensitivity  \label{fig-robust}}
	{This figure illustrates the impacts of heterogeneous customers' waiting sensitivities on the optimal service deployment and the value of service proliferation. 
 The skewness of the Beta distributions in the top panel, middle panel, and bottom panel is zero, positive, and negative, respectively.
 }
	\vspace{-0.2in}
\end{figure}

Our first set of numerical analysis examines heterogeneous customers' sensitivities of waiting, as seen in the results reported in Figures \ref{fig-robust} and \ref{fig-robust-zero-lw}. Our analysis reveals that the main research findings and insights are very robust across different models of customer heterogeneity. Specifically, though the variability and skewness of customer heterogeneity may change the values, the structure of optimal service deployment and the key insights on the value of service proliferation from various perspectives remain the same. We find that a positive or negative skewness generally reduces the value of service proliferation. This is because customers are less heterogeneous with a large positive or negative skewness, and thus, a single service mechanism is able to capture this highly concentrated waiting sensitivity without introducing internal competition between the two services.  For example, under the $\beta(1,40)$ with a high positive skewness, most customers are not sensitive to waiting (see also Figure \ref{ec-fig-beta-pdf}), in this case, the SP has almost no incentive to operate two services (the Regime $\mathcal{H}$ shrinks or even disappears). We remark that our results, with a large positive skewness of customers' waiting sensitivities, are consistent with \cite{lobel2024frontiers}, who show that the value of hybrid staffing is marginal in the context of a single service with waiting-insensitive customers that could be fulfilled by either full-time employees or part-time contractors or both.

Our second set of numerical analyses tests the impacts of the variability of independent agents' reservation rates, as shown in Figure  \ref{fig-robust-r}. Again, we find that our main results are robust across different distributions of the reservation rates, and that service proliferation is generally more valuable to the SP when the existing service is less cost-efficient and the new service is more cost-efficient. For example, adding on-demand (or standard) service is more (or less) valuable when the hourly wage of employees is high. A higher variability of contractors' reservation rates generally increases the value of adding the standard service due to the increased uncertainty in the supply of contractors, which reduces the efficiency of the on-demand service mechanism.

\subsection{Heterogeneous Service Quality}\label{subsec:quality}
The base model assumes that customers have the same valuation for the standard and the on-demand service, reflecting a homogeneous quality between these two services. In this subsection, we extend the base model to consider heterogeneous qualities between these two services by allowing different valuations of service. For example, an Uber XL service is usually provided by a premium car that has a higher quality than an economy car fulfilling a Uber X service. 
Specifically, customers derive a different value from the on-demand service, with the following utilities of a type-$\theta$ consumer from the standard and the on-demand service:  $U_s(p_s,W_s|\theta) = V - p_s - \theta W_s$ and $U_o(p_o,W_o|\theta) =\alpha V - p_o - \theta W_o$,
where $\alpha>0$ measures the value of an on-demand service relative to the value of a standard service, and $\alpha<1$ ($>1$)  means that the on-demand service is less (or more) valuable to customers than the standard service. 

\begin{figure}[htbp]
	\vspace{-0.1in}
	\FIGURE
	{
        \begin{minipage}{0.9\textwidth}
                \includegraphics[width=0.32\textwidth]{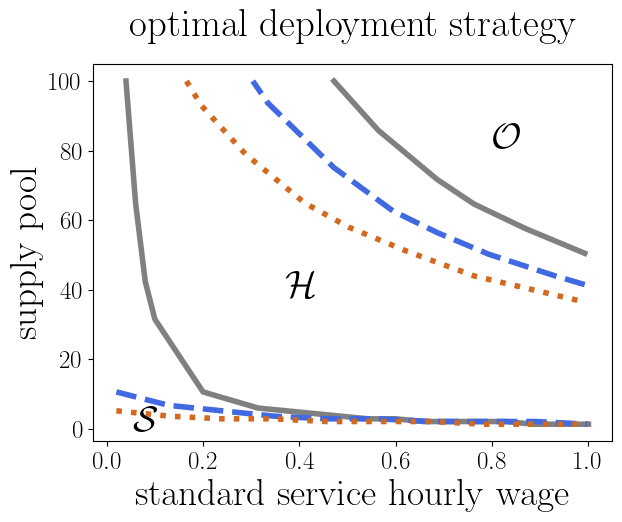} \hfill
                \includegraphics[width=0.32\textwidth]{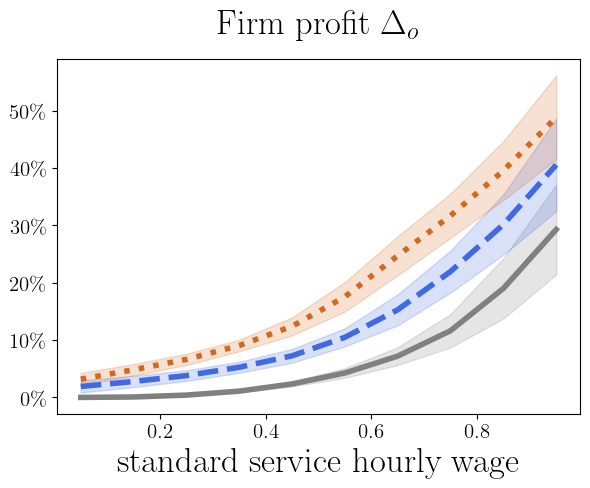} \hfill 
                \includegraphics[width=0.32\textwidth]{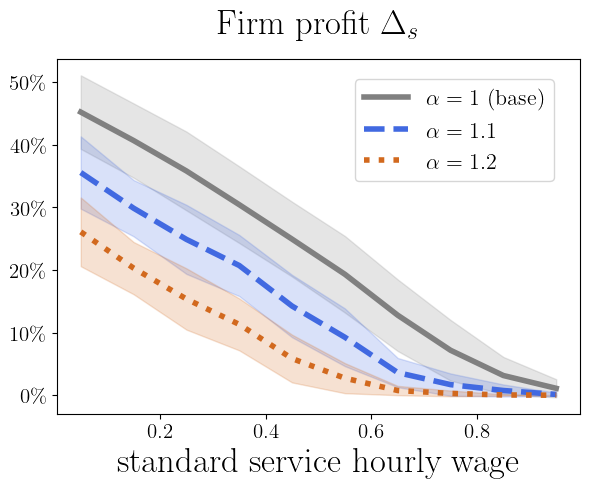}
            \end{minipage}
	}
	{\bf Optimal deployment and value of service proliferation: Heterogeneous service quality  \label{fig-diffV}}
	{}
	\vspace{-0.2in}
\end{figure}

Our analysis shows that our key research findings derived from the base model continue to hold when services have heterogeneous qualities. Specifically, though the quality differentiation may change the value of service proliferation and alter the specific partition of regimes in the optimal service deployment, it does not qualitatively change the structure of optimal service deployment and the results on the value of service proliferation from various perspectives. As we illustrate in  Figure \ref{fig-diffV} and Figure \ref{fig-diffV-lw}, when the on-demand service has a higher quality, it is more valuable to proliferate with the on-demand service while being less valuable to add the standard service from the perspective of the SP, consumers, and society as a whole.

\section{Concluding Remarks} \label{sec-conclusion} 
Our paper studies the service deployment for a profit-maximizing service provider (SP) that decides whether to offer a single service or two differentiated services, along with pricing and staffing of a workforce with employees, or contractors, or both, to price- and waiting-sensitive customers. Using a queueing-theoretic framework, we show that the SP is optimal to offer either a single service with one resource or two differentiated services with a blended workforce depending on the interaction between three factors: (i) customer valuation; (ii) the hourly wage of employees; (iii) the supply pool of contractors. Our analysis suggests that the service proliferation with blended workforce could improve the profit of the SP significantly and establish tight upper bounds on the value of service differentiation. 
 
There are fruitful research opportunities in the design of gig platforms to improve efficiency and sustainability. One future direction is to incorporate the strategic interaction between customers and the platform, e.g., how the information about waiting time should be disclosed and what the impacts of information disclosure are on the value of the blended workforce. It will also be interesting to incorporate fairness regulation into the SP to limit the price differentials between differentiated services, and examine how it changes the insights regarding the impacts of the blended workforce derived in this work.


\theendnotes


\normalem
\bibliographystyle{informs2020}
\bibliography{reference}


%
%
%

 \begin{APPENDICES}



\newpage
\setcounter{page}{1}
\setcounter{section}{0}
\renewcommand{\theequation}{A-\arabic{equation}}
\renewcommand{\thethm}{A-\arabic{thm}}
\renewcommand{\thelem}{A-\arabic{lem}}
\renewcommand{\theprop}{A-\arabic{prop}}
\renewcommand{\thecor}{A-\arabic{cor}}
\renewcommand{\thedefi}{A-\arabic{defi}}
\renewcommand{\thefigure}{A-\arabic{figure}}
\setcounter{equation}{0}
\setcounter{thm}{0}
\setcounter{lem}{0}
\setcounter{prop}{0}
\setcounter{cor}{0}
\setcounter{defi}{0}
\setcounter{figure}{0}

\begin{center}
	{\large Online Appendix \\
		\vspace{.1in}
		\Large \PaperTitle
		\vspace{.1in}
	}
\end{center}


\section{Proofs} 


\noindent
{\bf Proof of Proposition \ref{S-only}.}
It follows from $\lambda_s = \Lambda {\sf Pr}(U_s \geq 0 ) = \Lambda {\sf Pr}(V - p_s - \theta W_s \geq 0 ) $ that 
\bea
\lambda_s =\Lambda   \frac{ V-p_s}{W_s} . \label{eq-s-lambda-1}
\eea  
Following \eqref{eq-s-lambda-1} and \eqref{W} we have $k_s = \frac{\Lambda (V-p_s)+1}{\mu_s W_s}$.  
The optimization problem \eqref{model-S} can be equivalently written as 
 \[
\pi^S = 
\max\limits_{0\leq \lambda_s \leq \Lambda,W_s\geq 0} 
 \pi(\lambda_s,W_s) =  \lambda_s \left(V- \frac{\lambda_sW_s}{\Lambda}\right) - \frac{w_s}{\mu_s W_s}-\frac{w_s\lambda_s}{\mu_s} \equiv \pi_1(\lambda_s,W_s) . 
\]
 
Since $\pi(\lambda_s,W_s)$  is concave in $W_s$. By the first order condition, i.e., $\frac{\partial \pi(\lambda_s,W_s)}{\partial W_s}|_{W_s=\bar W_s(\lambda_s)}=0$, we derive $\bar W_s(\lambda_s)= \frac{1}{\lambda_s}\sqrt{\frac{\Lambda w_s}{\mu_s}}$, and consequently
\begin{align}\label{pi1s}
 \pi(\lambda_s,\bar W_s(\lambda_s)) =\lambda_s \left(V-2\sqrt{\frac{w_s}{\Lambda \mu_s}}-\frac{w_s}{\mu_s}\right)\equiv \pi_1(\lambda_s).
\end{align}
Denote the optimal solution  as $\lambda_s^S=\arg\max_{0\leq \lambda_s \leq \Lambda} \pi(\lambda_s, \bar W_s(\lambda_s))$.
\qed

\ 

\noindent
{\bf Proof of Proposition \ref{O-only}.} 
It follows from $\lambda_o = \Lambda {\sf Pr}(U_o \geq 0 ) = \Lambda {\sf Pr}(V - p_o - \theta W_o \geq 0 ) $ that 
\beq
\lambda_o =\Lambda   \frac{ V-p_o}{W_o}  
\eeq  
Since $r$ follows a uniform distribution, \eqref{ko} becomes 
\begin{align}
    k_o^2=\lambda_o K w_o.\label{2024-May-3}
\end{align}
It follows from \eqref{W} and \eqref{2024-May-3} that $C_o(w_o)= \frac{1}{K\mu_o^2}(\frac{1}{W_o} + \lambda_o )^2$, then the optimization problem \eqref{model-O} can be equivalently written as 
\bea 
\pi^O = \max \limits_{0 \leq \lambda_o \leq \Lambda,   W_o\geq 0}   \lambda_o \left(V - \frac{\lambda_o W_o}{\Lambda}\right) - \frac{1}{K \mu_o^2 }\left( \frac{1}{W_o} + \lambda_o\right)^2 \equiv \pi_2(\lambda_o,W_o) . \label{model-case2}
\eea 
Let $L_o=\lambda_o W_o$, we can rewrite \eqref{model-case2} as
 \[ 
\pi^O =  
\max\limits_{0\leq \lambda_o \leq \Lambda ,L_o\geq 0}  \pi(\lambda_o,L_o) = \lambda_o\left(V-\frac{L_o}{\Lambda}\right) - \frac{\lambda_o^2}{K\mu_o^2}\left(\frac{1}{L_o} + 1 \right)^2.
\] 
Since  $\pi(\lambda_o,L_o) $ is concave in 
$\lambda_o$, by the first order condition, i.e., $\frac{\partial \pi(\lambda_o,L_o)}{ \partial \lambda_o}|_{\lambda_o=\bar \lambda_o(L_o)}=0$, we derive  $\bar \lambda_o(L_o)=\frac{K\mu_o^2L_o^2(\Lambda V-L_o)}{2\Lambda (1+L_o)^2}$. Then, for fixed $L_o$, the best $\lambda_o(L_o) = \max\{ 0, \min\{\bar \lambda_o(L_o) ,\Lambda\} \}$.  
The constraints $ \bar \lambda_o(L_o) \geq 0$ becomes
\begin{align*}
\Lambda V \geq L_o,
\end{align*}
and $\bar \lambda_o(L_o)  \leq \Lambda$ becomes
\begin{align*}
 K\mu_o^2 L_o^2 (V \Lambda - L_o)  \leq    2\Lambda^2 (1 + L_o)^2. 
\end{align*}
Note that $\frac{L_o^2 (V \Lambda - L_o)}{(1 + L_o)^2}$ first increases and then decreases in $L_o \in [0,V\Lambda]$. 
Therefore, if $K > \frac{2\Lambda^2}{\mu_o^2} \left( \max_{  0\leq L_o\leq V\Lambda } \frac{L_o^2 (V \Lambda - L_o)}{(1 + L_o)^2} \right)^{-1} \equiv \Tilde{K}$, then $\bar \lambda_o(L_o)  \leq \Lambda$ can be equivalently written as $L_o \notin (L_1, L_2)$, where $0 < L_1 \leq L_2 < V\Lambda$ are two roots of $K\mu_o^2 L_o^2 (V \Lambda - L_o)  =  2\Lambda^2 (1 + L_o)^2$.  Then, we can rewrite the optimization as 
\bea 
\pi^O =   \max\left\{  \begin{array}{rl} 
& \max\limits_{ L_o \in [0,L_1] \cup [L_2,V\Lambda] }  \pi( \bar \lambda_o(L_o),L_o)=\frac{K\mu_o^2 L_o^2 (V \Lambda - L_o)^2}{4\Lambda^2 (1 + L_o)^2}  \equiv\pi_1 (L_o), \\
& \max\limits_{L_o \in (L_1,L_2)} \pi( \Lambda,L_o)= V\Lambda - L_o - \frac{\Lambda^2 (1+L_o)^2}{K\mu_o^2 L_o^2}\equiv \pi_2 (L_o) , \\
& \max\limits_{L_o \in  (V\Lambda, \infty) }  \pi(0,L_o) = 0.
\end{array} \right\} \label{2024-May-7}
\eea 
If $K < \Tilde{K}$, then $\bar \lambda_o(L_o)  \leq \Lambda$ holds for any $  L_o \in [0,V\Lambda]$, and the optimization is $\pi^O = \max\left\{ \max\limits_{ L_o \in [0,V\Lambda] }  \pi_1 (L_o), \ 0 \right\} $. 

The first order derivative of $\pi_1(L_o)$ is
\[ \frac{ {\sf d} \pi_1(L_o)}{{\sf d} L_o}  = \frac{K\mu_o^2 L_o(V\Lambda -L_o)(-L_o^2 -2L_o +V\Lambda)}{2\Lambda^2(1+L_o)^3}, \]
and the first order derivative of $\pi_2(L_o)$ is
\[  \frac{ {\sf d} \pi_2(L_o)}{{\sf d} L_o} = -1 + \frac{2\Lambda^2 (1+L_o)}{K\mu_o^2 L_o^3} . \]
We can show $\pi_1( L_o)$ first increases and then decreases in $L_o \in [0,V\Lambda]$ with changeover point $\sqrt{V \Lambda + 1} - 1$, and $\pi_2( L_o)$ first increases and then decreases in $L_o$ with changeover point $\bar L_o$ that is defined in \eqref{eq-barK-barLo}. In addition, $\frac{\partial \pi_1(L_o)}{\partial L_o} = \frac{\partial \pi_2(L_o)}{\partial L_o} 
$ when $L_o=  L_1  $ and $L_o=  L_2 $. Therefore, if $\sqrt{V \Lambda + 1} - 1 \in [0,L_1]$, then $\pi_2(L_o)$ must decrease on $(L_1,L_2)$; if $\sqrt{V \Lambda + 1} - 1 \in  [L_2,V\Lambda]$,  then $\pi_2(L_o)$ must increase on $(L_1,L_2)$; if $\sqrt{V \Lambda + 1} - 1 \in (L_1,L_2)$, then $\pi_2(L_o)$ must first increase then decrease on $(L_1,L_2)$. 
Note that $\sqrt{V \Lambda + 1} - 1 \leq \Lambda V$ always holds. Therefore, $\sqrt{V \Lambda + 1} - 1  \in [0,L_1] \cup [L_2,V\Lambda]$ is equivalent to $K\le \bar K$ according to definition of $\bar K$ in \eqref{eq-barK-barLo}.
In addition, we have $\Tilde{K} \leq \bar K$ since  $\max_{  0\leq L_o\leq V\Lambda } \frac{L_o^2 (V \Lambda - L_o)}{(1 + L_o)^2} \geq  \frac{L_o^2 (V \Lambda - L_o)}{(1 + L_o)^2}|_{L_o = \sqrt{V\Lambda+1}-1} = \frac{(\sqrt{V\Lambda + 1}-1)^3}{\sqrt{V \Lambda +1}} $.

Denote the optimal solution of \eqref{2024-May-7} as $L^O_o$ and $\lambda_o^O=\lambda_o(L_o^O)$. 
Therefore,  if $ K \leq \bar K$, then
$\lambda_o^O = \bar\lambda_o(\sqrt{V \Lambda + 1}-1 )= \frac{K\mu_o^2 (\sqrt{V\Lambda + 1}-1)^3}{2\Lambda \sqrt{V\Lambda+1}}$, $W_o^O = \frac{2\Lambda \sqrt{V\Lambda+1}}{K\mu_o^2 (\sqrt{V\Lambda + 1}-1)^2}$, $p_o^O = V - \frac{\lambda_o^O W_o^O}{\Lambda} = V - \frac{\sqrt{V\Lambda+1}-1}{\Lambda}$, $k_o^O= \frac{1}{\mu_o} (\frac{1}{W_o^O}+\lambda_o^O) $, $w_o^O = \frac{k_o^{O2}}{\lambda_o^O K}=\frac{1}{2}(V - \frac{\sqrt{V\Lambda+1}-1}{\Lambda})$, 
and $\pi^O = \pi_1(\sqrt{V \Lambda + 1}-1)
= \frac{K\mu_o^2 (\sqrt{V\Lambda +1}-1)^4}{4\Lambda^2}$.   
If $  K > \bar K$, then $\lambda_o^O=\Lambda$, $W_o^O = \frac{\bar L_o}{\Lambda}$, $p_o^O = V - \frac{\bar L_o}{\Lambda}$, $k_o^O= \frac{\Lambda}{\mu_o} (\frac{1}{\bar L_o}+1)$, $w_o^O = \frac{ (1+\bar L_o) \bar L_o }{2 \Lambda} = \frac{\bar L_o (1+\bar L_o)}{2\Lambda}$, 
and  $\pi^O = \pi_2(\bar L_o) = V\Lambda - \frac{\bar L_o ( \bar L_o +3)}{2}$ where  
$\bar L_o $ is defined in \eqref{eq-barK-barLo}, which obviously decreases in $K$, implying that the profit $ V\Lambda - \frac{\bar L_o ( \bar L_o +3)}{2}$ increases with $K$.

In addition, if $K \leq \bar K$, it follows from $p_o^O = V - \frac{\sqrt{V\Lambda+1}-1}{\Lambda}$ and $w_o^O = \frac{1}{2}(V - \frac{\sqrt{V\Lambda+1}-1}{\Lambda})$ that $\frac{w_o^O}{p_o^O} = \frac{1}{2}$. 
\qed 

\

\noindent{\bf Proof of Corollary \ref{cor-staff}.}
Following Proposition \ref{S-only} and Table \ref{ec-table-opt-solutions}, the optimal staffing solution in System S satisfies $k_s = \frac{\lambda_s}{\mu_s} + \sqrt{\frac{\lambda_s}{\mu_s w_s}}$. 

Following \eqref{W}, we have $k_o = \frac{\lambda_o}{\mu_o} + \frac{1}{W_o \mu_o}$, therefore, to prove $k_o = \frac{\lambda_o}{\mu_o} + \sqrt{\frac{\lambda_o}{\mu_o \frac{w_o\lambda_o}{k_o} \max\{ \frac{2\bar K}{K},2\} }}$, it is equivalent to prove $W_o = \sqrt{ \frac{w_o}{k_o \mu_o} \max\{ \frac{2\bar K}{K},2\}}$.  
By Proposition \ref{O-only} and Table \ref{ec-table-opt-solutions}, when $K > \bar K$, it follows from $W_o^O = \frac{\bar L_o}{\Lambda}$, $k_o^O = \frac{\Lambda (1 +\bar L_o)}{\mu_o \bar L_o}$ and $w_o^O = \frac{\bar L_o(1 + \bar L_o)}{2\Lambda}$ that $ W_o^O = \sqrt{ \frac{ w_o^O}{k_o^O \mu_o}\cdot 2}$; when $K \leq \bar K$, it follows from $W_o^O = \frac{\bar K (\sqrt{V\Lambda+1}-1)}{K \Lambda}$, $k_o^O =  \frac{K\Lambda \sqrt{V\Lambda+1}}{\bar K \mu_o (\sqrt{V\Lambda+1}-1)}$ and $w_o^O = \frac{\sqrt{V\Lambda+1}(\sqrt{V\Lambda+1}-1)}{2\Lambda}$ that $ W_o^O = \sqrt{ \frac{w_o^O}{k_o^O \mu_o} \cdot\frac{2\bar K}{K}}$. 
\qed 

\

\noindent
{\bf Proof of Theorem \ref{thm-compare-SO}. } 
Denote $R(w_s , K)$ as the performance ratio of the SP profit under the System S relative to System O, i.e.,
\bea \label{eq-ratio}
R(w_s , K) = \frac{\pi^S}{\pi^O} = \left\{  \begin{array}{ll} 
\frac{ 4 \Lambda^3 (V - \overline{C}_s )}{K\mu_o^2 (\sqrt{V\Lambda +1}-1)^4}, & \text{ if }  K  <  \bar K,  \\
 \frac{V - \overline{C}_s}{V - \overline{C}_o}, &  \text{ if } K \ge \bar K,
\end{array} \right.  
\eea 
where $\overline{C}_s, \overline{C}_o$ and $\bar K$ are defined in \eqref{eq-barK-barLo}. 
Note that $\overline{C}_o$ decreases in $K$. Hence, for a fixed $w_s$, $R(w_s,K) \leq 1$ if and only if $K \geq K_F(w_s)$, where $K_F(w_s)$ is given by $\frac{K\mu_o^2 (\sqrt{V\Lambda +1}-1)^4}{4\Lambda^2} = \pi^S$ or $\Lambda ( V - \overline{C}_o)= \pi^S$. Since $\overline{C}_s$ increases in $w_s$, we have $K_F(w_s)$ decreases in $w_s$. 

Following from Proposition \ref{S-only},  $\pi^S=V\Lambda$ when $w_s = 0$. Since $\pi^S$ is a decreasing function of $w_s$, it follows that $\pi^S \leq V\Lambda$ for any $w_s\geq 0$. Following from Proposition \ref{O-only}, we have $\Lambda (V - \overline{C}_o) \leq  \frac{K\mu_o^2 (\sqrt{V\Lambda +1}-1)^4}{4\Lambda^2}$ for any $K$ because $\pi_1(L_o)$ and $\pi_2(L_o)$ are continuous at $L_1$ and $L_2$. Therefore, we have 
\[ R(w_s,K) \leq \frac{V}{V- \overline{C}_o}.\]
Following from Proposition \ref{O-only}, as $K \rightarrow \infty$, we have  $\bar L_o \rightarrow 0$, and then $\pi^O \rightarrow V \Lambda$.  Additionally, since $\pi^O$ is an increasing function of $K$, it follows that $\pi^O$ is upper bounded by $V\Lambda$. Then, 
\[ R(w_s,K) \geq \frac{\pi^S}{V\Lambda} = 1 - \frac{\overline{C}_s}{V}.\]
Last, according to the continuity of $R(w_s,K)$  and above analysis, we can derive the conditions under which the limits are tight. 
\qed 

\

\noindent
{\bf Proof of Theorem \ref{thm-compare-LW-CS}.}
We first show the labor welfare. 
From  Proposition \ref{S-only}, we have $k_s^S = \frac{1}{\mu_s}(\frac{1}{W_s^S} + \lambda_s^S) = \frac{1}{\mu_s} (\sqrt{\frac{\mu_s \Lambda}{w_s}} + \Lambda)$. Then the labor welfare in System S is 
\[ LW^S = k_s^S w_s -  k_s^S \int_{0}^{1} r \mathbb{I}_{\{r \leq w_s \}}dr =k_s^S w_s - \frac{ k_s^Sw_s^2}{2} = (1 - \frac{w_s}{2}) (\sqrt{\frac{w_s \Lambda}{\mu_s}} + \frac{w_s \Lambda}{\mu_s}). \]
As a function of $w_s$, $LW^S =0$ has two real roots, 0 and 2. Therefore $LW^S$ first increases and then decreases in $w_s$. 
Following from Proposition \ref{O-only}, the  labor welfare in System O is 
\beq
LW^O = K \int_0^1 ( \frac{\lambda_o^O w_o^O}{k_o^O} -r)^+ dr = K \int_0^1 ( \frac{k_o^O}{K} -r)^+ dr  = \frac{(k_o^{O})^2}{2K} 
=\left\{  \begin{array}{ll} 
K\mu_o^2(\sqrt{V\Lambda+1}-1)^4 / (8\Lambda^2) & \text{ if } K \leq \bar K, \\
\Lambda^2 (1 + \bar L_o)^2 / (2K \mu_o^2 \bar L_o^2) & \text{ if } K > \bar K.
\end{array} \right. 
\eeq
Note that $\bar L_o$ is given by $ 2\Lambda^2 (1+\bar L_o) = K\mu_o^2 \bar L_o^3 $, so $\Lambda^2 (1 + \bar L_o)^2 / (2K \mu_o^2 \bar L_o^2)= \bar L_o ( 1+ \bar L_o) / 4 $. Since $\bar L_o$ decreases in $K$, we have $ \bar L_o ( 1+ \bar L_o) /4$ decreases in $K$. Hence, $LW^O$ first increases and then decreases in $K$ with changeover point $\bar K$, and the largest labor welfare is $ \sqrt{V \Lambda +1}(\sqrt{V\Lambda+1}-1) /4$. 

Therefore, 
$LW^O > LW^S$ if and only if $K \in (\underline{K}_L(w_s), \overline{K}_L(w_s))$, where $\underline{K}_L(w_s)=  \frac{8\Lambda^2 LW^S}{\mu_o^2(\sqrt{V\Lambda+1}-1)^4}$ if $LW^S \leq  \frac{\sqrt{V \Lambda +1}(\sqrt{V\Lambda+1}-1)}{4}$ and $\underline{K}_L(w_s)=\infty$ if $LW^S >  \frac{\sqrt{V \Lambda +1}(\sqrt{V\Lambda+1}-1)}{4}$,  and $\overline{K}_L(w_s)$ is given by $\bar L_o ( 1+ \bar L_o) = 4 LW^S$ if $LW^S \leq  \frac{\sqrt{V \Lambda +1}(\sqrt{V\Lambda+1}-1)}{4}$ and  $\overline{K}_L(w_s)= \infty$ if $LW^S> \frac{\sqrt{V \Lambda +1}(\sqrt{V\Lambda+1}-1)}{4}$. 
Since $LW^S$ first increases and then decreases in $w_s$, we have $\underline{K}_L(w_s)$ first increases and then decreases, or increases in $w_s$, and $\overline{K}_L(w_s)$ first decreases and then increases, or decreases in $w_s$. 

Next, we consider the consumer surplus. Note that, the consumer surplus in System S is 
\[ CS^S = \Lambda \int_0^1 (V - p_s^S - \theta W_s^S)^+ d \theta  = \Lambda \int_0^1 (W_s^S - \theta W_s^S) d \theta = \frac{\Lambda W_s^S}{2} = \frac{1}{2}\sqrt{\frac{w_s\Lambda}{\mu_s}}, \]
and  the consumer surplus in System O is 
\[
CS^O = \Lambda \int_0^1 (V - p_o^O - \theta W_o^O)^+ d \theta  = \Lambda \int_0^{\frac{V-p_o^O}{W_o^O}} (V - p_o^O - \theta W_o^O) d \theta = \frac{\Lambda (V-p_o^O)^2}{2W_o^O} 
=\left\{  \begin{array}{ll} 
\frac{K\mu_o^2 (\sqrt{V\Lambda + 1}-1)^4}{4\Lambda^2 \sqrt{V\Lambda+1}} & \text{ if } K \leq \bar K,  \nonumber \\
\bar L_o  /2 & \text{ if } K > \bar K.
\end{array} \right. 
\]
Since $\bar L_o$ decreases in $K$, we have $CS^O$ first increases and then decreases in $K$ with changeover point $\bar K$, and the highest consumer surplus is $\frac{\sqrt{V\Lambda+1}-1}{2}$. 

Note that $\frac{\sqrt{V\Lambda+1}-1}{2} \geq \frac{1}{2}\sqrt{\frac{w_s\Lambda}{\mu_s}}$ is equivalent to $\pi_1(\Lambda) = V\Lambda - \frac{w_s\Lambda}{\mu_s} - 2\sqrt{\frac{w_s\Lambda}{\mu_s}} \geq 0$. Therefore, for a fixed $w_s$, $CS^O > CS^S$ if and only if $K \in (\underline{K}_C(w_s), \overline{K}_C(w_s))$, where $\underline{K}_C(w_s)$ is given by $\frac{K\mu_o^2 (\sqrt{V\Lambda + 1}-1)^4}{4\Lambda^2 \sqrt{V\Lambda+1}} = CS^S$, and $\overline{K}_C(w_s)$ is given by $\frac{\bar L_o}{2} = CS^S$. 
Since $CS^S$ increases in $w_s$, we have $\underline{K}_C(w_s)$ increases in $w_s$ and $\overline{K}_C(w_s)$ decreases in $w_s$. 

Finally, we prove the results on the social welfare. Note that, the social welfare in System S is 
\beq
SW^S &=& \pi^S+ LW^S + CS^S 
=  V\Lambda -   \frac{(w_s)^2\Lambda}{2\mu_s} - \frac{1+w_s}{2} \sqrt{\frac{w_s\Lambda}{\mu_s}},
\eeq
which decreases on $w_s$. If $K \leq \bar K$, then the social welfare in System O is 
\beq
SW^O &=& \pi^O + LW^O + CS^O 
= \frac{K\mu_o^2 (\sqrt{V\Lambda + 1}-1)^4 ( 3\sqrt{V\Lambda+1} + 2)}{8\Lambda^2 \sqrt{V\Lambda+1}},
\eeq
which increases in $K$. 
If $K > \bar K$, then  \[SW^O = \pi^O + LW^O + CS^O 
= V\Lambda - \frac{\bar L_o ( 3+\bar L_o)}{4},\]which increases in $K$ since $\bar L_o$ decreases in $K$. Hence, for a fixed $w_s$, $SW^O > SW^S$ if and only if $K > K_S(w_s)$, where $K_S(w_s)$ is given by $\frac{K\mu_o^2 (\sqrt{V\Lambda + 1}-1)^4 ( 3\sqrt{V\Lambda+1} + 2)}{8\Lambda^2 \sqrt{V\Lambda+1}} = SW^S$ or $V\Lambda - \frac{\bar L_o ( 3+\bar L_o)}{4} = SW^S$. Since $SW^S$ decreases in $w_s$, we have $K_S(w_s)$ decreases in $w_s$. 
\qed 

\

\noindent
{\bf Proof of Corollary \ref{cor-joint-compare}.}
Following from Theorems \ref{thm-compare-SO}-\ref{thm-compare-LW-CS}, System S induces an all-win situation if and only if
$K \leq \min\{K_F(w_s),K_S(w_s), \underline{K}_L(w_s),\underline{K}_C(w_s)\}$, or $    \max\{\overline{K}_L(w_s),\overline{K}_C(w_s)\} \leq K \leq \min\{ K_F(w_s),K_S(w_s)\} $.  We will prove the second condition is empty by showing that $K_F(w_s)< \overline{K}_C(w_s)$. 

First, $\overline{K}_C(w_s)$ is given by $CS^O = CS^S$ when $K \geq \bar K_o$, i.e., $ \bar L_{o,C} = \sqrt{\frac{w_s\Lambda}{\mu_s}}$ where  
$\bar L_{o,C} $ satisfies $2\Lambda^2 (1+\bar L_{o,C}) = \overline{K}_C \mu_o^2 \bar L_{o,C}^3$. Second, when $K \geq \bar K_o$, $K_F(w_s)$ is given by $\pi^O  = \pi^S$, i.e., $ \bar L_{o,F} (\bar L_{o,F} + 3) = 2 \frac{w_s\Lambda}{\mu_s} + 4 \sqrt{\frac{w_s\Lambda}{\mu_s}}$, where $2\Lambda^2 (1+\bar L_{o,F}) = K_F \mu_o^2 \bar L_{o,F}^3$. Therefore, $K_F(w_s)< \overline{K}_C(w_s)$ if and only if  $\bar L_{o,F} > \bar L_{o,C}$, or equivalently,  $\bar L_{o,F} (\bar L_{o,F}+3) > \bar L_{o,C}(\bar L_{o,C}+3)$, which is equivalent to $ \frac{w_s \Lambda}{\mu_s} + \sqrt{\frac{w_s \Lambda}{\mu_s}} > 0$. Then we have  $K_F(w_s)< \overline{K}_C(w_s)$ for any $w_s >  0$. 
\qed

\ 

\noindent{\bf Proof of Theorem \ref{opt-two}.} 
The proof relies on five auxiliary Lemmas \ref{lem-subproblem}-\ref{lem-pi-mono} that are presented and proved in Appendix \ref{ec-sec-lemma}. We establish this result by following five steps:

First, we simplify the optimization problem \eqref{opt} by variable transformation, derive a one-to-one mapping between $(p_s,p_o,k_s,w_o)$ and $(\lambda_s,W_s,\lambda_o,W_o)$, and then rewrite $\pi^*$ in a simpler way. The equivalent form of $\pi^*$ is given by   Lemma  \ref{lem-subproblem}. The subscript in $\pi^H_{\lambda_s, W_s, \lambda_o,W_o}$ indicates that the decision variables are $(\lambda_s,W_s,\lambda_o,W_o)$.  

Second, we solve the best $\lambda_s$, $W_s$ and $W_o$ sequentially. Specifically, we solve $\lambda_s$ in Lemma \ref{lem-solve-ls}, then solve $W_s$ in Lemma \ref{lem-solve-Ws}, and finally solve $W_o$ in Lemma \ref{lem-solve-loWo}. 

Third, we show the monotonicity of $\pi^H$, $\pi^S$ and $\pi^O$, see  Lemma \ref{lem-pi-mono}, and then characterize the switching curves to define the regions for each possible optimal solution.  Since $\pi^S$ is independent of $K$ and $\pi^H$ increases in $K$, for each given $w_s$, there exists a unique switching curve $K_{SH}(w_s)$ satisfying $\pi^S = \pi^H$, and 
\[ \pi^S \geq \pi^H   \text{ if and only if } K \leq K_{SH}(w_s).\]
Since $\pi^O$ is independent of $w_s$ and $\pi^H$ decreases in $w_s$, for each given $K$, there exists a unique $w_{s,OH}(K)$ satisfying $\pi^O  = \pi^H$, and 
\[ \pi^O \geq \pi^H  \text{ if and only if } w_s \geq w_{s,OH}(K).\]  
Following from Theorem \ref{thm-compare-SO}, we have  
\[\pi^O \geq \pi^S \text{  if and only if }  K \geq K_F(w_s).\] 

We next show that $ K_{SH}(w_s)$ decreases in $w_s$ by showing that  $\frac{\partial \pi^S}{\partial w_s} \leq \frac{\partial \pi^H}{\partial w_s}$, i.e., as $w_s$ increases, $\pi^S$ decreases faster than $\pi^H$. Following from Proposition \ref{S-only}, we have $ \frac{\partial \pi^S}{\partial w_s} = - \frac{\Lambda}{\mu_s} - \sqrt{\frac{\Lambda}{w_s \mu_s}}$.  Following from the proof of Lemma \ref{lem-solve-loWo}, for a fixed $W_o$, we have
\[\frac{\partial \pi_{ol}(\lambda_o,w_s)}{\partial w_s} = - \frac{\Lambda - \lambda_o}{\mu_s } - \sqrt{\frac{(\Lambda + \lambda_o)(\Lambda - \lambda_o)}{ \Lambda w_s \mu_s}}) \geq \frac{\partial \pi_{ol}(\lambda_o,w_s)}{\partial w_s} \Big|_{\{\lambda_o =0 \}} =   - \frac{\Lambda}{\mu_s} - \sqrt{\frac{\Lambda}{w_s \mu_s}}=\frac{\partial \pi^S}{\partial w_s},  \]
where $\pi_{ol}(\lambda_o,w_s)$ is defined by \eqref{2024-0819-1} and the inequality comes from $\frac{\partial \pi_{ol}(\lambda_o,w_s)}{\partial w_s}$ increases in $\lambda_o$. 
Similarly, for a fixed $W_o$, we have
\[\frac{\partial \pi_{sl}(\lambda_o,w_s) }{\partial w_s} = -(\Lambda - \lambda_o) ( \frac{1}{\mu_s} + \sqrt{\frac{1}{\Lambda w_s \mu_s}}) \geq \frac{\partial \pi_{sl}(\lambda_o,w_s)}{\partial w_s} \Big|_{\{\lambda_o =0\}} = - \frac{\Lambda}{\mu_s} - \sqrt{\frac{\Lambda}{w_s \mu_s}} = \frac{\partial \pi^S}{\partial w_s}, \]
where the inequality comes from $\frac{\partial \pi_{sl}(\lambda_o,w_s)}{\partial w_s}$ increases in $\lambda_o$. Therefore, we can conclude that  $\pi^S$ has a faster decreasing rate than  $\pi^H$ given any $w_s > 0$. 

Last, we can define regions as follows, 
\beq
\mathcal{S} &=& \{(w_s,K):  K \leq \min\{K_{SH}(w_s), K_F(w_s) \} \} ,\\ 
\mathcal{O} &=& \{(w_s,K):  K \geq  K_F(w_s), \ w_s \geq w_{s,OH}(K) \}, \\
\mathcal{H} &=& \{(w_s, K): K \geq  K_{SH} (w_s),   w_s \leq w_{s,OH}(K)  \}.
\eeq
\qed

\ 

\noindent
{\bf Proof of Theorem  \ref{pro:compareTS}.} 
Following from Proposition \ref{S-only}, in System S, the best price is $p_s^S = V - \sqrt{\frac{w_s}{\Lambda \mu_s}}$, and the best lead time is $W_s^S = \sqrt{\frac{w_s}{\Lambda \mu_s}}$. 
Following from Lemma \ref{lem-solve-Ws}, we consider the following two cases.

If $(\lambda_o,W_o) \in \mathcal{R}_1$, then $\lambda_s = \Lambda - \lambda_o$ and $W_s = \sqrt{\frac{\Lambda w_s}{\mu_s (\Lambda + \lambda_o)(\Lambda - \lambda_o)}}$.  Note that $W_o \geq W_s$ in this region, so that  $p_s = V - \frac{\lambda_o W_s + \lambda_s W_s}{\Lambda} = V - W_s $ following \eqref{eq-price-ol}. Since $\sqrt{\frac{\Lambda w_s}{\mu_s (\Lambda + \lambda_o)(\Lambda - \lambda_o)}} > \sqrt{\frac{w_s}{\Lambda \mu_s}}$ for any $ 0 < \lambda < \Lambda$, we have $W_s > W_s^S$ and $p_s < p_s^S$ in this case.

If  $(\lambda_o,W_o) \in \mathcal{R}_2$, then $\lambda_s = \Lambda - \lambda_o$ and $W_s = \frac{1}{\Lambda - \lambda_o} \sqrt{\frac{\Lambda w_s}{\mu_s}}$. Note that $W_o \leq  W_s$ in this region, so that $p_s = V - \frac{\lambda_s W_s +\lambda_o W_o}{\Lambda} = V - \sqrt{\frac{w_s}{\Lambda \mu_s}} - \frac{\lambda_o W_o}{\Lambda}$ following \eqref{eq-price-sl}. Since $\frac{1}{\Lambda - \lambda_o} \sqrt{\frac{\Lambda w_s}{\mu_s}} > \sqrt{\frac{w_s}{\Lambda \mu_s}}$ for any $ 0 < \lambda < \Lambda$, we have $W_s > W_s^S$, and $p_s = V - W_s^S - \frac{\lambda_o W_o}{\Lambda} < V - W_s^S = p_s^S$. 

To summarize, for any fixed $(\lambda_o,W_o)$, we have $p_s < p_s^S$ and $W_s > W_s^S$.
\qed 

\

\noindent
{\bf Proof of Theorem \ref{pro:compareTO}.} 
Following from Proposition \ref{O-only}, in System O,  when $K\leq \bar K$, the best price is $p_o^O =  \frac{\Lambda V+1-\sqrt{\Lambda V+1}}{\Lambda}$ and the best lead time is $W_o^O =  \frac{2\Lambda \sqrt{\Lambda V+1}}{K\mu_o^2(\sqrt{\Lambda V+1}-1)^2} = \frac{\bar K ( \sqrt{V\Lambda+1}-1)}{K \Lambda}$. When $K > \bar K$, the best price is $p_o^O = V - \frac{\bar L_o}{\Lambda}$ and the best lead time is $W_o^O = \frac{\bar L_o}{\Lambda}$. Following from  Lemma \ref{lem-solve-Ws}, we consider the following two cases.

Case 1: $W_o \geq W_s$. The objective function in this case is
$\pi_{ol}(\lambda_o,W_o)= \pi_2(\lambda_o,W_o) + (\Lambda - \lambda_o) (V - \frac{w_s }{\mu_s} ) - 2 \sqrt{\frac{w_s (\Lambda + \lambda_o)(\Lambda - \lambda_o)}{\Lambda \mu_s }}$, where $\pi_2(\lambda_o,W_o)$ is defined in  \eqref{model-case2}. 
Following from Lemma \ref{lem-solve-Ws} and Lemma \ref{lem-solve-loWo}, we know that if the optimal $(\lambda_o,W_o)$ is on the boundary, the optimal profit equals $\pi^O$. Hence, we only focus on interior optimum. That is, $(\lambda_o^*, W_o^*)$ is given by $\frac{\partial \pi_{ol}(\lambda_o, W_o)}{\partial W_o} = - \frac{\lambda_o^2}{\Lambda} + \frac{2(1+ \lambda_o W_o)}{K \mu_o^2 W_o^3} = 0$ and $\frac{\partial \pi_{ol}(\lambda_o, W_o)}{\partial \lambda_o} = - \frac{2\lambda_o W_o}{\Lambda} - \frac{2(1+\lambda_o W_o)}{K \mu_o^2 W_o} + \frac{w_s}{\mu_s} + 2 \sqrt{\frac{w_s}{\Lambda \mu_s}} \frac{\lambda_o}{\sqrt{\Lambda^2 - \lambda_o^2}} = 0$. 
Let $L_o =  \lambda_o W_o$, we can rewrite $\frac{\partial \pi_{ol}(\lambda_o, W_o)}{\partial W_o} = 0$ and $\frac{\partial \pi_{ol}(\lambda_o, W_o)}{\partial \lambda_o} =0$ as 
\bea 
&& -\frac{L_o^2}{\Lambda} + \frac{2 \lambda_o (1+L_o)}{K\mu_o^2 L_o} = 0, \  -\frac{2L_o}{\Lambda} - \frac{2\lambda_o (1+L_o)}{K \mu_o^2 L_o} + \frac{w_s}{\mu_s} + 2 \sqrt{\frac{w_s}{\Lambda \mu_s}}\frac{\lambda_o}{\sqrt{\Lambda^2 - \lambda_o^2}} = 0, \label{eq-ol-fod-4}
\eea 
which implies $-\frac{2L_o}{\Lambda} - \frac{L_o^2}{\Lambda}+ \frac{w_s}{\mu_s} + 2 \sqrt{\frac{w_s}{\Lambda \mu_s}}\frac{\lambda_o}{\sqrt{\Lambda^2 - \lambda_o^2}} = 0$. 
Then we obtain $\lambda_o^*$ and $L_o^*$: 
\bea 
L_o^* = \sqrt{1 + \frac{w_s \Lambda}{\mu_s} + 2 \sqrt{\frac{w_s \Lambda}{\mu_s}}\frac{  \lambda_o^*}{\sqrt{\Lambda^2 - \lambda_o^{*2}}}  } - 1, \label{eq-ol-opt-Lo}
\eea 
and $\lambda_o^*$ is given by 
\bea 
\frac{(\sqrt{1 + \frac{w_s \Lambda}{\mu_s} + 2 \sqrt{\frac{w_s \Lambda}{\mu_s}}\frac{  \lambda_o^*}{\sqrt{\Lambda^2 - \lambda_o^{*2}}}  } - 1)^3}{\sqrt{1 + \frac{w_s \Lambda}{\mu_s} + 2 \sqrt{\frac{w_s  \Lambda}{\mu_s}}\frac{ \lambda_o^*}{\sqrt{\Lambda^2 - \lambda_o^{*2}}}  } } - \frac{2\Lambda \lambda_o^*}{K \mu_o^2} = 0.  \label{eq-ol-opt-lambdao}
\eea 

In this case, $\lambda_s = \Lambda - \lambda_o$ and $W_s = \frac{1}{\Lambda - \lambda_o} \sqrt{\frac{\Lambda w_s}{\mu_s}}$. Note that $W_o \geq W_s$ in this region, so  $p_o = V - \frac{\lambda_o W_o + \lambda_s W_s}{\Lambda} = V - \frac{L_o}{\Lambda} - \sqrt{\frac{ w_s (\Lambda - \lambda_o)}{\mu_s \Lambda(\Lambda + \lambda_o)}}$ following \eqref{eq-price-ol}.  
For $K\leq \bar K$, $p_o^* \leq p_o^O$  is equivalent to $ L_o^*  \geq  \sqrt{V\Lambda + 1}-1 - \sqrt{\frac{ w_s (\Lambda - \lambda_o^*)}{\mu_s (\Lambda + \lambda_o^*)}}$, and $W_o^* \geq W_o^O$  is equivalent to $L_o^*  \geq \lambda_o^* \frac{\bar K ( \sqrt{ \Lambda V+1}-1)}{K \Lambda}$. Using \eqref{eq-ol-fod-4} and \eqref{eq-ol-opt-Lo}, we can rewrite the two inequalities as
\begin{align}
\sqrt{1 + \frac{w_s \Lambda}{\mu_s} + 2 \sqrt{\frac{w_s \Lambda}{\mu_s}}\frac{  \lambda_o^*}{\sqrt{\Lambda^2 - \lambda_o^{*2}}}  }  + \sqrt{\frac{ w_s (\Lambda - \lambda_o^*)}{\mu_s (\Lambda + \lambda_o^*)}} - \sqrt{V\Lambda + 1} \geq 0,  \ \frac{2\Lambda }{K \mu_o^2} -  \frac{(\lambda_o^*)^2 (\frac{\bar K ( \sqrt{ \Lambda V+1}-1)}{K \Lambda})^3}{\lambda_o^* \frac{\bar K ( \sqrt{ \Lambda V+1}-1)}{K \Lambda} + 1} \geq 0. \label{20240923-1}
\end{align}
For $K > \bar K$, $p_o^* \leq p_o^O$ is equivalent to $L_o^* \geq \bar L_o - \sqrt{\frac{ w_s (\Lambda - \lambda_o^*)}{\mu_s (\Lambda + \lambda_o^*)}}$, and $W_o^* \geq W_o^O$ is equivalent to $L_o^* \geq \frac{ \lambda_o^* \bar L_o}{\Lambda}$. Using \eqref{eq-ol-fod-4} and \eqref{eq-ol-opt-Lo}, we can rewrite the two inequalities as 
\bea 
\sqrt{1 + \frac{w_s \Lambda}{\mu_s} + 2 \sqrt{\frac{w_s \Lambda}{\mu_s}}\frac{  \lambda_o^*}{\sqrt{\Lambda^2 - \lambda_o^{*2}}}  }  - 1  +  \sqrt{\frac{ w_s (\Lambda - \lambda_o^*)}{\mu_s (\Lambda + \lambda_o^*)}} - \bar L_o \geq 0, \ \frac{2\Lambda }{K \mu_o^2} -  \frac{(\lambda_o^*)^2 (\frac{\bar L_o}{\Lambda})^3}{\lambda_o^* \frac{\bar L_o}{\Lambda} + 1} \geq 0. \label{20240923-3}
\eea 
We define 
\bea 
\mathcal{T}_1 &=&  \{ (w_s,K) | \eqref{eq-ol-opt-lambdao},\eqref{20240923-1}, K \leq \bar K \}  \cup \{ (w_s,K) | \eqref{eq-ol-opt-lambdao},\eqref{20240923-3}, K > \bar K \}.
\eea 

Case 2: $W_o \leq W_s$. The objective function in this case is
$ \pi_{sl}(\lambda_o,W_o)= \pi_2(\lambda_o,W_o) + (\Lambda - \lambda_o)(V - \frac{w_s}{\mu_s}- 2\sqrt{\frac{w_s}{\Lambda \mu_s}} - \frac{2\lambda_o W_o}{\Lambda})$. The interior optimum $(\lambda_o^*, W_o^*)$ is given by $\frac{\partial \pi_{sl}(\lambda_o, W_o)}{\partial W_o}  = \frac{ 2 ( 1+\lambda_o W_o)}{K \mu_o^2 W_o^3} - \frac{\lambda_o (2\Lambda - \lambda_o)}{\Lambda} = 0$ and $\frac{\partial \pi_{sl}(\lambda_o, W_o)}{\partial \lambda_o}  = - \frac{2(\Lambda - \lambda_o)W_o}{\Lambda} - \frac{2(1+\lambda_o W_o)}{K \mu_o^2 W_o} + \frac{w_s}{\mu_s} + 2 \sqrt{\frac{w_s}{\Lambda \mu_s}} = 0$. 
Let $L_o = \lambda_o W_o$, we can rewrite $\frac{\partial \pi_{sl}(\lambda_o, W_o)}{\partial W_o}  = 0$ and $\frac{\partial \pi_{sl}(\lambda_o, W_o)}{\partial \lambda_o}  = 0$ as 
\bea 
&& \frac{ 2 \lambda_o ( 1+ L_o)}{K \mu_o^2 L_o} - \frac{ (2\Lambda - \lambda_o)L_o^2}{\Lambda \lambda_o} = 0,  \  - \frac{2(\Lambda - \lambda_o)L_o}{\Lambda \lambda_o} - \frac{2\lambda_o (1+L_o)}{K \mu_o^2 L_o} + \frac{w_s}{\mu_s} + 2 \sqrt{\frac{w_s}{\Lambda \mu_s}} = 0, \label{eq-sl-fod-4}
\eea 
which implies $- \frac{2(\Lambda - \lambda_o)L_o}{\Lambda \lambda_o} - \frac{(2\Lambda - \lambda_o)L_o^2 }{\Lambda \lambda_o} + \frac{w_s}{\mu_s} + 2 \sqrt{\frac{w_s}{\Lambda \mu_s}} = 0$. 
Then we can obtain $\lambda_o^*$ and $L_o^*$: 
\bea 
L_o^* = \frac{-\Lambda + \lambda_o^* + \sqrt{(\Lambda-\lambda_o^*)^2 + (2\Lambda - \lambda_o^*)  \lambda_o^*(\frac{w_s \Lambda}{\mu_s} + 2 \sqrt{\frac{w_s \Lambda}{ \mu_s}})}}{2\Lambda - \lambda_o^*}, \label{eq-sl-opt-Lo}
\eea 
and $\lambda_o^*$ is given by 
\bea 
\frac{(\frac{-\Lambda + \lambda_o^* + \sqrt{(\Lambda-\lambda_o^*)^2 + (2\Lambda - \lambda_o^*)  \lambda_o^*(\frac{w_s \Lambda}{\mu_s} + 2 \sqrt{\frac{w_s \Lambda}{ \mu_s}})}}{2\Lambda - \lambda_o^*})^3}{\frac{-\Lambda + \lambda_o^* + \sqrt{(\Lambda-\lambda_o^*)^2 + (2\Lambda - \lambda_o^*)  \lambda_o^*(\frac{w_s \Lambda}{\mu_s} + 2 \sqrt{\frac{w_s \Lambda}{ \mu_s}})}}{2\Lambda - \lambda_o^*} + 1} - \frac{ 2 \Lambda \lambda_o^2}{K \mu_o^2 (2\Lambda - \lambda_o) } = 0. \label{eq-sl-opt-lambdao}
\eea 

In this case, $\lambda_s = \Lambda - \lambda_o$ and $W_s = \frac{1}{\Lambda - \lambda_o} \sqrt{\frac{\Lambda w_s}{\mu_s}}$. Note that $W_o \leq W_s$ in this region, so   $p_o = V - \frac{\lambda_s W_o + \lambda_o W_o}{\Lambda} = V - W_o$ following \eqref{eq-price-sl}. 
For $K\leq \bar K$,   $p_o^* \geq p_o^O$  is equivalent to $L_o^* \leq \lambda_o^* \frac{\sqrt{\Lambda V+1}-1}{\Lambda}$, 
and $W_o \leq W_o^O$  is equivalent to $L_o^* \leq \lambda_o^* \frac{\bar K ( \sqrt{ \Lambda V+1}-1)}{K \Lambda}$. 
Note that $K\leq \bar K$ implies $\frac{\bar K ( \sqrt{ \Lambda V+1}-1)}{K \Lambda} \geq \frac{   \sqrt{ \Lambda V+1}-1 }{  \Lambda}$, i.e., $L_o^* \leq \lambda_o^* \frac{\sqrt{\Lambda V+1}-1}{\Lambda}$ implies $L_o^* \leq \lambda_o^* \frac{\bar K ( \sqrt{ \Lambda V+1}-1)}{K \Lambda}$. 
Using \eqref{eq-sl-fod-4}, i.e., $\frac{(L_o^*)^3}{L_o^* + 1} = \frac{2\Lambda (\lambda_o^*)^2}{K \mu_o^2 (2\Lambda - \lambda_o^*)}$, we can rewrite $L_o^* \leq \lambda_o^* \frac{\sqrt{\Lambda V+1}-1}{\Lambda}$ as
\bea 
\frac{2\Lambda }{K \mu_o^2 (2\Lambda - \lambda_o^*)} - \frac{\lambda_o^* (\frac{\sqrt{\Lambda V+1}-1}{\Lambda})^3}{\lambda_o^* \frac{\sqrt{\Lambda V+1}-1}{\Lambda} + 1 } \leq 0. \label{20240923-sl-1}
\eea 
For $K > \bar K$, both $p_o^* \geq p_o^O$ and  $W_o^* \leq W_o^O$ are equivalent to $L_o^* \leq  \lambda_o^* \frac{\bar L_o}{\Lambda}$; using \eqref{eq-sl-fod-4}, it equals to 
\bea 
\frac{2\Lambda }{K \mu_o^2 (2\Lambda - \lambda_o^*)} - \frac{\lambda_o^* (\frac{\bar L_o}{\Lambda})^3}{\lambda_o^* \frac{\bar L_o}{\Lambda} + 1 } \leq 0. \label{20240923-sl-2}
\eea 
We define 
\bea 
\mathcal{T}_2 &=&  \{ (w_s,K) | \eqref{eq-sl-opt-lambdao},\eqref{20240923-sl-1}, K \leq \bar K \}  \cup \{ (w_s,K) | \eqref{eq-sl-opt-lambdao}, \eqref{20240923-sl-2}, K > \bar K \}.
\eea 
\qed



\

\noindent
{\bf Proof of Theorem \ref{thm-pi-K-ws}.}  
Following from Lemma \ref{lem-pi-mono},   $\pi^H$ increases in $K$ and decreases in $w_s$;  $\pi^S$ decreases in $w_s$ and remains unchanged with respect to $K$; $\pi^O$ increases in $K$ and remains unchanged with respect to $w_s$. Therefore, we can summarize that,  $\Delta_o = 1 - \frac{\pi^S}{\pi^*}$ is increasing  in $K$   and $\Delta_s = 1 - \frac{\pi^O}{\pi^*}$ is decreasing in $w_s$. 

To find the upper bound of $\Delta_o$, we show $\pi^H = \pi^O$ when $K$ is large. 
Following the proof of Lemma \ref{lem-solve-loWo},  we first show that, to achieve $\pi^H$,  the best $W_o$ is either $\hat W_o(\lambda_o)$ or $\check W_o(\lambda_o)$, or the boundary on which $\pi^H = \pi^O$.  Specifically we show that the best $W_o$ will not lie on the boundary  $W_o = \frac{1}{\lambda_o}( \sqrt{\frac{\Lambda w_s (\Lambda +\lambda_o)}{\mu_s (\Lambda - \lambda_o)}} - \sqrt{\frac{\Lambda w_s}{\mu_s}}) $.
It follows from  $\frac{K \mu_o^2 \lambda_o^2}{2\Lambda} = \frac{1+ \lambda_o \hat W_o(\lambda_o)}{ \hat W_o(\lambda_o)^3}$, $\frac{K \mu_o^2 \lambda_o(2\Lambda - \lambda_o)}{2\Lambda} = \frac{1+ \lambda_o \check W_o(\lambda_o)}{\check W_o(\lambda_o)^3}$ and $\lambda_o < 2\Lambda - \lambda_o$ that $\check W_o(\lambda_o) < \hat W_o(\lambda_o) $. Therefore, if $\hat W_o(\lambda_o) <  \frac{1}{\lambda_o}( \sqrt{\frac{\Lambda w_s (\Lambda +\lambda_o)}{\mu_s (\Lambda - \lambda_o)}} - \sqrt{\frac{\Lambda w_s}{\mu_s}}) $, i.e., $\hat W_o(\lambda_o)$ is infeasible, then we have $\check W_o(\lambda_o) < \hat W_o(\lambda_o) < \frac{1}{\lambda_o}( \sqrt{\frac{\Lambda w_s (\Lambda +\lambda_o)}{\mu_s (\Lambda - \lambda_o)}} - \sqrt{\frac{\Lambda w_s}{\mu_s}})$ and the optimal solution regarding $\pi^H$ is given by $  \check W_o(\lambda_o)$. Similarly, if $\check W_o(\lambda_o) >  \frac{1}{\lambda_o}( \sqrt{\frac{\Lambda w_s (\Lambda +\lambda_o)}{\mu_s (\Lambda - \lambda_o)}} - \sqrt{\frac{\Lambda w_s}{\mu_s}}) $, i.e., $\check W_o(\lambda_o)$ is infeasible, then we have $\hat W_o(\lambda_o)  > \check W_o(\lambda_o) > \frac{1}{\lambda_o}( \sqrt{\frac{\Lambda w_s (\Lambda +\lambda_o)}{\mu_s (\Lambda - \lambda_o)}} - \sqrt{\frac{\Lambda w_s}{\mu_s}}) $ and the optimal solution regarding $\pi^H$ is given by $  \hat W_o(\lambda_o)$.  

Second, we show that given the best $W_o$ is either $\hat W_o(\lambda_o)$ or $\check W_o(\lambda_o)$, the best $\lambda_o \in \mathcal{R}_1  \cup \mathcal{R}_2$ should be as large as possible if $K$ is large. It suffices to show for sufficiently large $K$ that
\bea
\frac{\partial \pi_{ol}(\lambda_o, \hat W_o(\lambda_o))}{\partial \lambda_o} \geq 0 \text{ and } \frac{\partial \pi_{sl}(\lambda_o,\check W_o(\lambda_o))}{\partial \lambda_o}  \geq 0.
\eea
Since 
$\hat W_o(\lambda_o)$ is given by $\frac{K \mu_o^2 \lambda_o^2}{2\Lambda} = \frac{1+ \lambda_o W_o}{ W_o^3}$ and $\frac{1+ \lambda_o W_o}{ W_o^3}$ decreases in $W_o$, we have $\hat W_o(\lambda_o)$ decreases in $K$. Similarly,  
$\check W_o(\lambda_o)$ is given by $\frac{K \mu_o^2 \lambda_o(2\Lambda - \lambda_o)}{2\Lambda} = \frac{1+ \lambda_o W_o}{W_o^3}$ and   $\check W_o(\lambda_o)$ decreases in $K$.  
According to the proof of Lemma \ref{lem-solve-loWo}, we have that for sufficiently large $K$
\beq 
\frac{\partial \pi_{ol}(\lambda_o, \hat W_o(\lambda_o))}{\partial \lambda_o} &=& - \frac{2(1+\lambda_o \hat W_o(\lambda_o))}{K \mu_o^2 \hat W_o(\lambda_o)} + \frac{w_s}{\mu_s} - \frac{2\lambda_o}{\Lambda} ( \hat W_o(\lambda_o) - \sqrt{\frac{\Lambda w_s}{\mu_s (\Lambda^2 -\lambda_o^2)}})\geq 0,\\
\frac{\partial \pi_{sl}(\lambda_o,\check W_o(\lambda_o))}{\partial \lambda_o} &=& - \frac{2(1+\lambda_o \check W_o(\lambda_o))}{K \mu_o^2 \check W_o(\lambda_o)} + \frac{w_s}{\mu_s} - \frac{2}{\Lambda} (\check W_o(\lambda_o) (\Lambda - \lambda_o) - \sqrt{\frac{\Lambda w_s}{\mu_s}})\geq0 .
\eeq 


Third, we define $K_{OH}(w_s)$ as follows
\[
K_{OH}(w_s)=\min \{K:  
  \frac{\partial \pi_{ol}(\lambda_o,\hat{W_o}(\lambda_o))}{\partial \lambda_o} \geq 0 \ \forall (\lambda_o,\hat{W_o}(\lambda_o)) \in \mathcal{R}_1,    \text{ and } \ \frac{\partial \pi_{sl}(\lambda_o,\check W_o \lambda_o)}{\partial \lambda_o} > 0 \ \forall (\lambda_o,\check{W_o}(\lambda_o)) \in \mathcal{R}_2\} \]
Therefore $\pi^H=\pi^O$ when $K \geq K_{OH}(w_s)$. 
Since $\pi^H$ increases in $K$, we can conclude that $\pi^H$ is upper bounded by $\pi^O  $ when $K = K_{OH}(w_s)$, which further implies for any $K\geq K_{OH}(w_s) $
\[ \Delta_o \leq 
1 - R(w_s, K_{OH}(w_s))\leq  \frac{\overline{C}_s}{V}, \]
where $R(w_s,K) = \pi^S/\pi^O$  and tha last inequality is due to Theorem \ref{thm-compare-SO}.

Next, we  find the upper bound of $\Delta_s$. In the proof of Theorem \ref{opt-two} we showed that  $\frac{\partial \pi^S}{\partial w_s} \leq \frac{\partial \pi^H}{\partial w_s}$, i.e., as $w_s$ increases, $\pi^S$ decreases faster than $\pi^H$. In addition, it is intuitive that $\pi^S \geq \pi^H$ when $w_s = 0$, since firm-scheduled employees are free. Together with the previous conclusion that $\pi^S$ has a faster decreasing rate than $\pi^H$ when $w_s > 0$, there must exist a unique intersection point between the two profit functions. Let $w_{s,SH}(K)$ denote the unique intersection point between $\pi^S$ and $\pi^H$, i.e., 
$w_{s,SH}(K)=\min \{w_s:  \pi^S\leq \pi^H \}.$
Then  we have $\pi^S \geq \pi^H$ for any standard labor wage $w_s \leq w_{s,SH}(K)$, which further implies 
$\Delta_s \leq  
1 - \frac{1}{R(w_{s,SH}(K),K)}\leq \frac{\overline{C}_o}{ V}, $ 
where the last inequality is due to Theorem \ref{thm-compare-SO}.
\qed 

\ 

\section{Supplemental Results}\label{ec-sec-supplements}
In this section, we first present a few supporting analytical results, then consider the flexible setting in which the SP can assign the standard service (or on-demand service) to contractors (or employees), and finally provide supplemental numerical reports to illustrate the robustness of our research findings.  All proofs of supplemental results in this section are relegated to Section \ref{ec-sec-supplemental-proofs} in the Online E-companion.
\subsection{Supplemental Lemmas and Theorems}\label{ec-sec-lemma}
In this subsection, we first present five auxiliary lemmas that are used to establish Theorem \ref{opt-two}, and provide a sufficient condition under which adding on-demand service is valuable to consumers. 

\begin{lem} \label{lem-subproblem}
$\pi^* = \max \{ \pi^S ,\pi^O ,\pi^H_{\lambda_s, W_s, \lambda_o,W_o}\} $ where 
\bea \label{model-H}
\pi^H_{\lambda_s, W_s, \lambda_o,W_o} = \left\{  \begin{array}{rl}
\max_{} & \pi(\lambda_s, W_s, \lambda_o,W_o) \\
\mbox{s.t.} & \lambda_s, \lambda_o, W_s,W_o > 0, \ \lambda_s + \lambda_o \leq \Lambda ,
\end{array}\right.  
\eea 
and $\pi(\lambda_s, W_s, \lambda_o,W_o) = V(\lambda_s + \lambda_o) - \frac{ \lambda_s^2 W_s + 2 \lambda_s \lambda_o \min\{W_o,W_s\} + \lambda_o^2 W_o  }{\Lambda} - \frac{w_s}{\mu_s} (\frac{1}{W_s} + \lambda_s) - \frac{ 1}{K\mu_o^2}(\frac{1}{W_o} + \lambda_o )^2$.
\end{lem} 

\begin{lem} \label{lem-solve-ls}
$\pi^* = \max\{\pi^S, \pi^O, \pi^H_{W_s, \lambda_o,W_o}\}$, where 
\bea 
\pi^H_{W_s, \lambda_o,W_o} = 
\max_{(W_s,\lambda_o, W_o) \in \mathcal{D}_3 \cup \mathcal{D}_6} \  \pi(\lambda_s, W_s, \lambda_o,W_o)\big|_{ \{\lambda_s= \Lambda - \lambda_o\}}  ,
\eea 
and $\mathcal{D}_3$ and  $\mathcal{D}_6$ are formally defined in the proof.
\end{lem}

\begin{lem} \label{lem-solve-Ws} 
$\pi^* = \{\pi^S, \pi^O, \pi^H_{\lambda_o,W_o} \} $, where 
\bea 
\pi^H_{\lambda_o,W_o} = \max\left\{  \begin{array}{rl} 
& \max_{(\lambda_o,W_o)\in \mathcal{R}_1} \pi(\lambda_s, W_s, \lambda_o,W_o) \big|_{ \{ \lambda_s = \Lambda - \lambda_o, W_s = \sqrt{\frac{\Lambda w_s}{\mu_s (\Lambda + \lambda_o)(\Lambda - \lambda_o)}}\}},   \\
&  \max_{(\lambda_o,W_o)\in \mathcal{R}_2}  \pi(\lambda_s, W_s, \lambda_o,W_o)\big|_{ \{\lambda_s = \Lambda - \lambda_o, W_s = \frac{1}{\Lambda - \lambda_o} \sqrt{\frac{\Lambda w_s}{\mu_s}}\} }   
\end{array} \right\}  \label{pilaWo}
\eea 
and $\mathcal{R}_1$ and $\mathcal{R}_2$ are formally defined in the proof.
\end{lem}

\begin{lem} \label{lem-solve-loWo}
$\pi^* = \max\{\pi^S, \pi^O,  \pi^{H}    \}$, where $\pi^{H} = \max\{ \pi^{H}_{ol} , \pi^{H}_{sl} \}$ and 
\beq 
\pi^{H}_{ol} &=& \max_{0\leq  \lambda_o  \leq \sqrt{\Lambda^2 - \frac{4\Lambda w_s \mu_s}{(V\mu_s -w_s)^2}} } \pi(\lambda_s, W_s, \lambda_o,W_o) \big|_{ \{ \lambda_s = \Lambda - \lambda_o, W_s = \sqrt{\frac{\Lambda w_s}{\mu_s (\Lambda + \lambda_o)(\Lambda - \lambda_o)}}, W_o = W_{o,ol}(\lambda_o) \}} , \nonumber \\
\pi^{H}_{sl} &=& \max_{0\leq \lambda_o\leq \Lambda} \pi(\lambda_s, W_s, \lambda_o,W_o)\big|_{ \{\lambda_s = \Lambda - \lambda_o, W_s = \frac{1}{\Lambda - \lambda_o} \sqrt{\frac{\Lambda w_s}{\mu_s}},  W_o = W_{o,sl}(\lambda_o)\} } ,
\eeq 
and $W_{o,ol}(\lambda_o) $ and $ W_{o,sl}(\lambda_o)$ are formally defined in the proof.
\end{lem}

\begin{lem} \label{lem-pi-mono}
    $\pi^H$ and $\pi^O$ increase in $K$,  $\pi^H$ and $\pi^S$ decrease in $w_s$. 
\end{lem}

\begin{thm} \label{prop-relation}
 {\sc (Value of the On-demand Service to the consumers)} 
     If $(w_s,K) \in \mathcal{H}$ and on-demand service has a longer lead time, then $CS^* \geq CS^S$. 
\end{thm}

Theorem \ref{prop-relation} indicates that introducing an on-demand service can benefit customers when the on-demand service aims to serve impatient customers. Note that we already showed in Theorem \ref{pro:compareTS} that introducing the on-demand service makes the standard service cheaper with a long lead time. As such, under the condition that the on-demand service has a longer lead time than the standard service in the hybrid system, adding the on-demand service lowers the prices of both services, thus leading to a higher consumer surplus. 

\subsection{Flexible Servers}\label{ec-sec-flexible}
Our base model assumes that each type of service could only be served by one specific type of server. That is, the standard-only service is only provided by standard employees, who are paid at $w_s$ per hour, and the on-demand service is only provided by on-demand contractors, who are paid at $w_o$ per service. In this section, we relax this assumption of dedicated servers and consider flexible servers whereby the SP could utilize on-demand contractors for standard services and/or employees for on-demand services. For example, it is not uncommon to observe that Uber occasionally provides a premium car, which is outsourced and usually serves a premium order of an Uber XL service, to serve an economic order of an Uber X service. Our objective in this section is to show the robustness of our research findings by extending the results to this general setting.

    \begin{figure}[!htb]
	\vspace{-0.1in}
	\FIGURE
	{
		\includegraphics[width=0.9\textwidth]{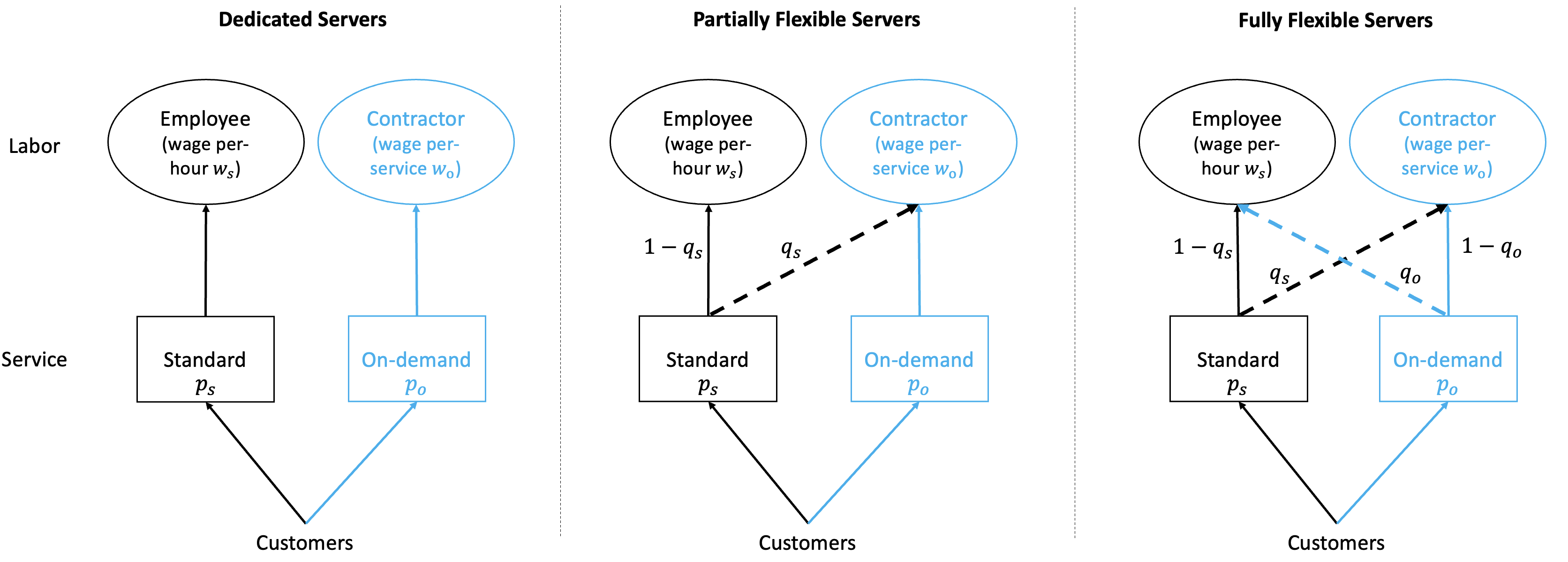}
	}
	{Dedicated vs flexible servers \label{fig-dedicated-flexible-servers}}
	{This figure illustrates the difference between dedicated and flexible servers. The two types of labor differ in the payment schedule, i.e., the employees and contractors are paid based on per hour and per service, respectively. The two types of services are differentiated in both the prices and the average waiting time, which is endogenously determined by supply and demand. }
    \vspace{-0.2in}
    \end{figure}
    
The interplay between different types of servers and services under the dedicated and flexible servers is illustrated in Figure \ref{fig-dedicated-flexible-servers}. Specifically, under the partially flexible setting, a customer who selects the on-demand service can only be served by a contractor, while she could be served by either a contractor or an employee when selecting the standard service. 
Under the fully flexible setting, a customer selecting any service could be served by either a contractor or an employee. In what follows, we focus on the fully flexible setting, that includes the partially flexible setting as a special case. 
Unlike the base model with dedicated servers, the SP determines and pre-announces the job assignment scheduling which specifies probabilities, or the percentages, $q_s$ (for scheduling a standard service to a contractor) and $q_o$ (for scheduling an on-demand service to an employee). The other procedures regarding customers' selection of services, the SP's pricing and staffing follow exactly the same as the base model. Specifically, customers are sensitive to both prices and waiting time, with the following expected utilities associated with each service:  
\beq
U_s = V - p_s - \theta[(1-q_s) W_s + q_s W_o], \quad U_o = V - p_o - \theta [q_o W_s + (1-q_o) W_o], 
\eeq
where $\theta\in U[0,1]$ captures customers' heterogeneity in the waiting, $(1-q_s) W_s + q_s W_o$ is the expected waiting time for a standard service, and $q_o W_s + (1-q_o) W_o$ is the expected waiting time for an on-demand service. Analogous to the base model, each customer selects the service to maximize her expected utility, thus resulting in the following effective arrival rates for   each service: 
\bea 
\lambda_s &=& (1-q_s) \Lambda {\sf Pr}(U_s \ge 0, U_s \ge U_o) + q_o \Lambda {\sf Pr} (U_o \ge 0, U_o \ge U_s) , \label{eq-lambda-flexible-1}\\
\lambda_o &= & q_s  \Lambda {\sf Pr}(U_s \ge 0, U_s \ge U_o) + (1-q_o) \Lambda {\sf Pr} (U_o \ge 0, U_o \ge U_s) . \label{eq-lambda-flexible-2}
\eea 
Note that $\lambda_s$  and  $\lambda_o$  defined in \eqref{eq-lambda-flexible-1}-\eqref{eq-lambda-flexible-2} are the arrivals served by employees and  contractors respectively. Using \eqref{eq-lambda-flexible-1} and \eqref{eq-lambda-flexible-2}, we can derive the arrivals who pay for the service at price $p_s$ as $\Lambda {\sf Pr}(U_s \geq 0, U_s \geq U_o) = \frac{\lambda_s - q_o(\lambda_s+\lambda_o)}{1-q_s-q_o}$, and also the arrivals who pay for the service at price $p_o$ as $\Lambda {\sf Pr}(U_o \geq 0, U_o \geq U_s) =\frac{\lambda_o - q_s( \lambda_s + \lambda_o)}{1-q_s-q_o}  $.

Anticipating the customers' strategic responses, the SP selects the job assignment schedules, as well as the prices and staffing of these two types of services, to maximize his expected profit. Analogous to \eqref{opt} in the base model, the SP's optimization problem is formulated as follows:  
\begin{align} \label{model-flexible}
&\pi_{f}^* =   \left\{\begin{array}{rl}
 \max\limits_{p_s,p_o,k_s,w_o, q_s,q_o} & \pi_{f}(p_s,p_o,k_s,w_o,q_s,q_o)= p_s \frac{\lambda_s - q_o(\lambda_s+\lambda_o)}{1-q_s-q_o} - w_s k_s + p_o\frac{\lambda_o - q_s( \lambda_s + \lambda_o)}{1-q_s-q_o}  - w_o \lambda_o   \\
\mbox{s.t.} 
& \eqref{eq-lambda-flexible-1}, \eqref{eq-lambda-flexible-2} ,\\
& W_i =1/( k_i\mu_i - \lambda_i), \ i \in \{s,o\}, \\
& k_o = K {\sf Pr}( \frac{ \lambda_o w_o}{k_o} \geq r).
\end{array}\right. 
\end{align}
We remark that the base model \eqref{opt} is a special case of \eqref{model-flexible} by setting $ q_s=0$ and $q_o = 0$, implying the dominance between these two optimal profits, i.e., $\pi^* \le \pi_f^*$. In the following, we show that these two systems are equivalent from the perspective of the SP. 

\begin{thm} {\sc (Profit Equivalence between Dedicated and Flexible System)}\label{lem-flexible}
The service provider derives an identical optimal profit between dedicated servers and flexible servers, i.e., $\pi^* = \pi_f^*$. 
\end{thm}

We have thus shown that the flexibility of servers has no value in the service system when the staffing, scheduling of workers, and pricing can be fully endogenized. This is in sharp contrast to the literature on flexibility (see e.g., \citealt{dong2021optimal}) that show a significant value of flexibility when these features cannot be fully optimized. 
 
\subsection{Supplemental Numerical Results}\label{ec-sec-numerical}
In this section, we provide extensive numerical analysis to complement the robustness analysis summarized in the main context. 
\subsubsection{Impacts of market size and supply pool on value of service proliferation} Our analysis reveals an important interplay between the market size and supply pool on the value of service proliferation. Our first set of numerical experiments examines the impacts of market size for a given supply pool, see Figure \ref{fig-Lambda-K-1} for a heat map visualization. We find that as the market size increases: 1) The value of on-demand service for the firm and society decreases, and 2) the value of on-demand service for consumers and labor increases. This is because the increased market size provides a greater benefit to the profit and social welfare of the benchmark system, System S. 
Additionally, the value of standard service for all players increases. This indicates that, compared to System O, an increasing market size has a more substantial impact on the hybrid system.

The second experiment investigates the scenario in which the market size and supply pool expand at the same rate, see Figure \ref{fig-Lambda-K-2} for the visualization. The findings provide two interesting insights regarding the impacts of supply and demand on the value of service proliferation. Specifically, as the number of self-scheduled providers and customer demand grow, at approximately the same rate: 1) the value of on-demand service for the firm, consumers, and society becomes higher, and 2) the value of standard service for all players becomes lower. This is because the increasing supply pool has a positive impact on the on-demand service, reducing the relative advantage of standard service. Our results indicate that the positive effect of the increasing supply pool on on-demand service has dominated the effect of the increasing market size on standard service.

\subsubsection{$M/M/1$ approximation vs $M/M/k$} Figure \ref{fig-mmk-lw} compares the value of service proliferation with respect to labor, consumers, and society under the base model with an $M/M/1$ queue and the extension with an $M/M/k$ queue. Note that we have shown in Figure \ref{fig-mmk} that the threshold for switching from hybrid services to on-demand service only becomes smaller under the $M/M/k$ queueing model.
Focusing on the common region where the optimal solution is hybrid services under both settings (i.e., when $w_s$ is not high), the trends are generally consistent across the two models, with one exception, $\Delta_s CS$, exhibits a different pattern.

\subsubsection{Heterogeneity in customers' waiting sensitivity} Figure \ref{fig-robust-zero-lw} compares the value of service proliferation with respect to labor, consumers, and society under different customer waiting sensitivity assumptions. Specifically, Figure \ref{fig-robust-zero-lw}(a) represents the zero skewness case, Figure \ref{fig-robust-zero-lw}(b) represents the positive skewness case, and Figure \ref{fig-robust-zero-lw}(c) represents the negative skewness case.
The trends are generally consistent across the different scenarios, with one exception. The value of service proliferation for consumer surplus, $\Delta_o CS$, in Figure \ref{fig-robust-zero-lw}(b) exhibits a different pattern. This is because, as more customers have a low sensitivity to waiting, the value of on-demand service to consumers first increases and then decreases.

\subsubsection{Heterogeneity in agents' reservation rate} Figure \ref{fig-robust-r} shows the impacts of different distributions of agents' reservation rate on the value of hybrid services. The results indicate that, given the same mean reservation rate, the value of hybrid services is robust to different shapes of the reservation rate distribution.

\begin{figure}[!hbt]
	\vspace{0in}
	\FIGURE
	{
            \includegraphics[width=0.4\textwidth]{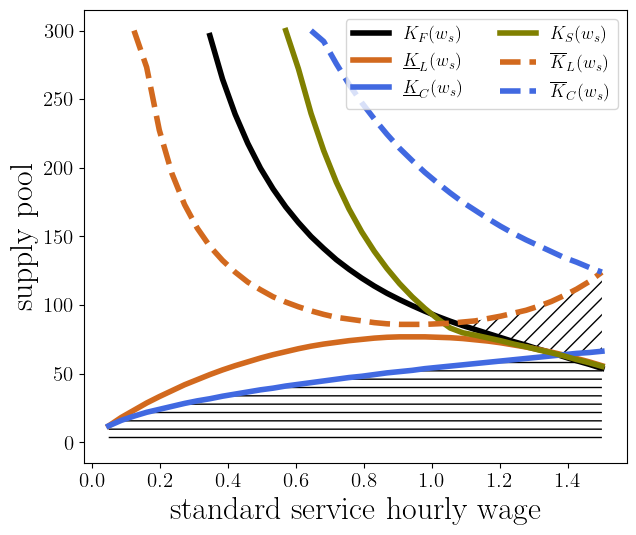} 
        }
	{\bf Incentive coordination: Standard vs on-demand \label{fig-compare-2}}
	{The region shaded with horizontal (or diagonal) lines represents the coordinated outcome where the SP, the labors, and the consumers are all better off under the standard (or on-demand) mechanism. Parameters are $\Lambda = 30, V = 2.2, \mu_s = 1.3, \mu_o = 1$. }
	\vspace{-0.2in}
\end{figure}

\subsubsection{Heterogeneous service quality} Figure \ref{fig-diffV-lw} compares the value of service proliferation with respect to labor, consumers, and society under different service quality scenarios. The key findings are:
1) As the on-demand service provides higher quality, the values of both services for labor are robust.
2) The values of the on-demand service for consumers and society become higher. This is because customers can realize a higher valuation by choosing the on-demand service when it offers superior quality.
3) The value of the standard service for society becomes lower. This is a consequence of the higher perceived value of the on-demand service option.

\begin{figure}[htbp]
	\vspace{-0.1in}
	\FIGURE
	{
            \begin{minipage}{\textwidth}
                \includegraphics[width=0.24\textwidth]{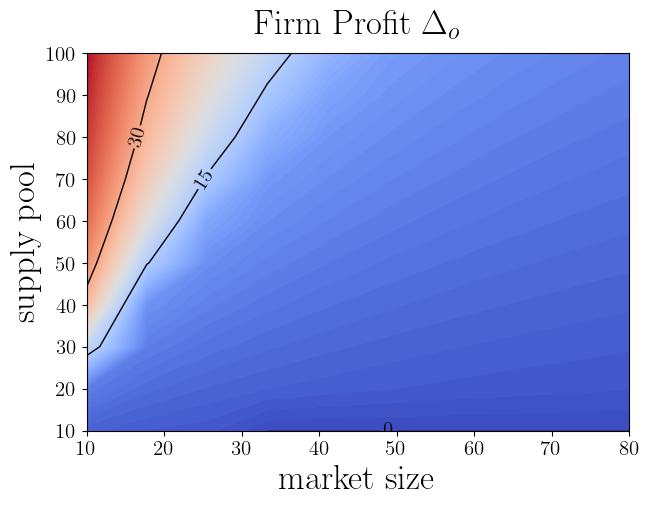} 
                \includegraphics[width=0.24\textwidth]{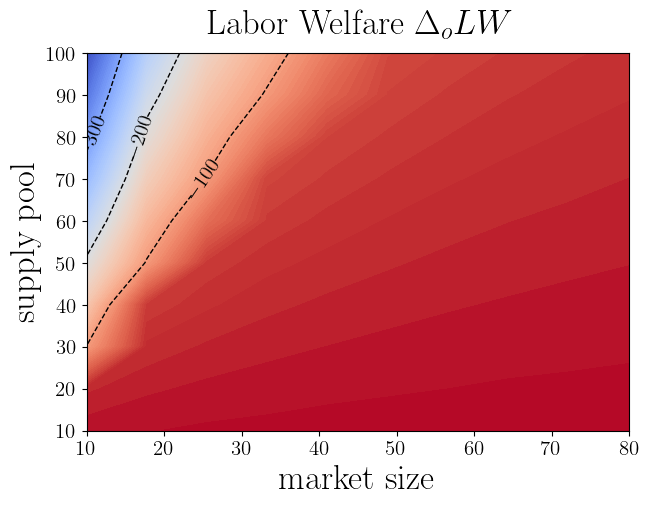}
                \includegraphics[width=0.24\textwidth]{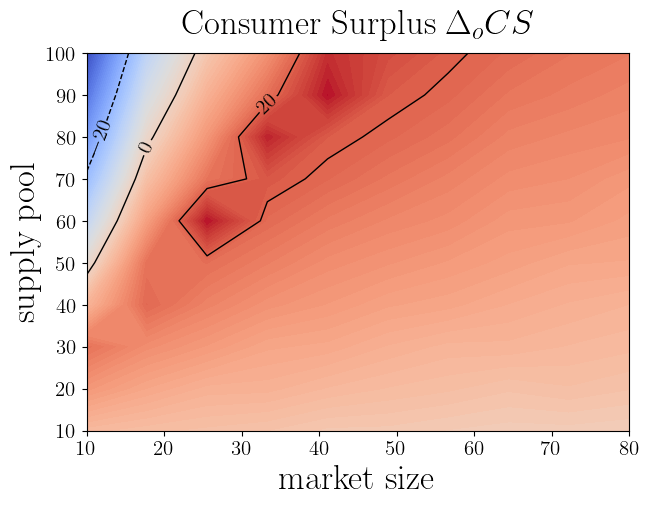} 
                \includegraphics[width=0.24\textwidth]{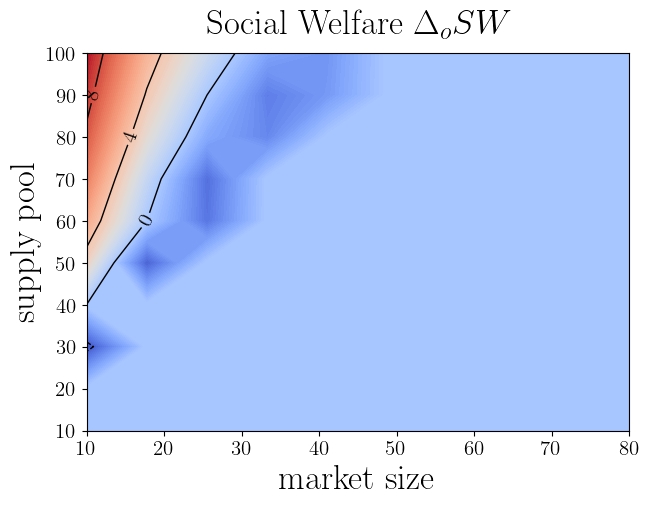} 
                \vfill
                \includegraphics[width=0.24\textwidth]{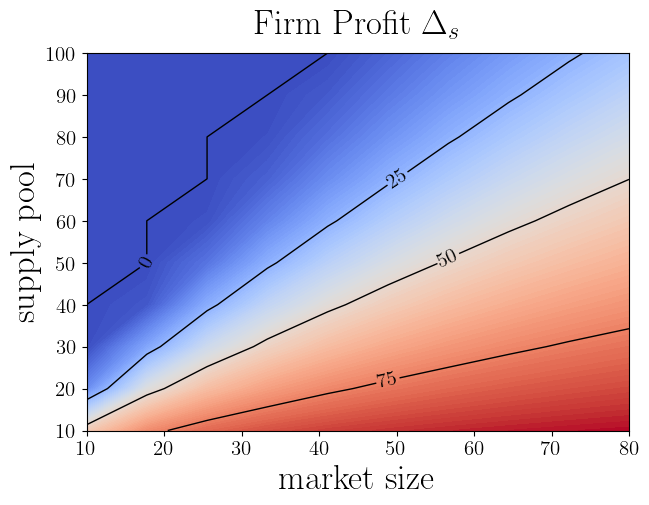} 
                \includegraphics[width=0.24\textwidth]{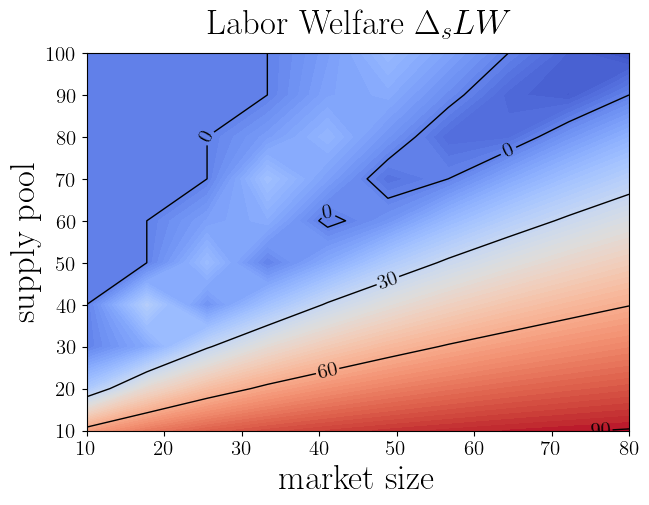} 
                \includegraphics[width=0.24\textwidth]{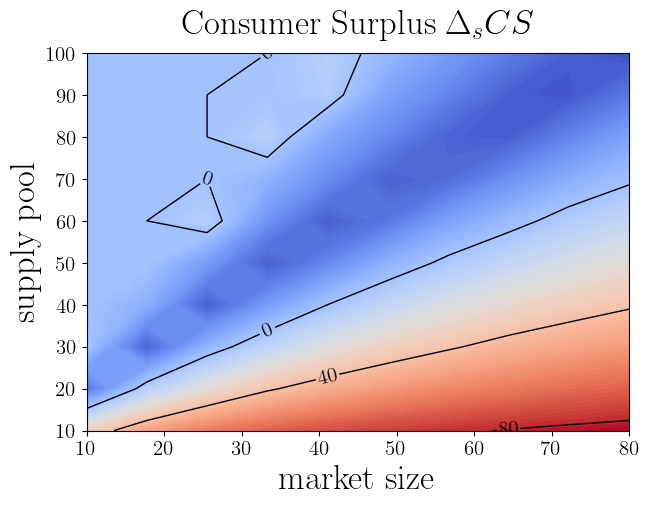} 
                \includegraphics[width=0.24\textwidth]{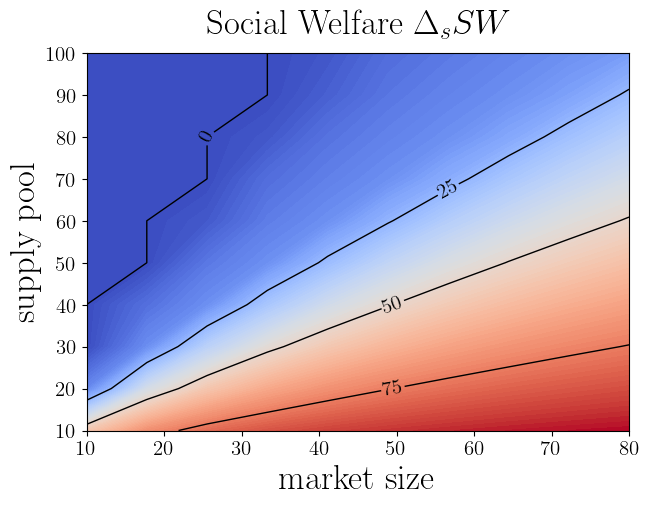}
            \end{minipage}
	}
	{\bf Value of service proliferation: The interplay between market size and supply pool. \label{fig-Lambda-K-1} }
	{This figure illustrates the interplay between market size $(\Lambda)$ and supply pool $(K)$ on the value of service proliferation. We change these two parameters independently and report the value of service proliferation from different perspectives with the hourly wage of standard employees $w_s = 0.6$.}
	\vspace{-0.1in}
\end{figure}

\begin{figure}[htbp]
	\vspace{-0.1in}
	\FIGURE
	{
            \begin{minipage}{\textwidth}
                \includegraphics[width=0.24\textwidth]{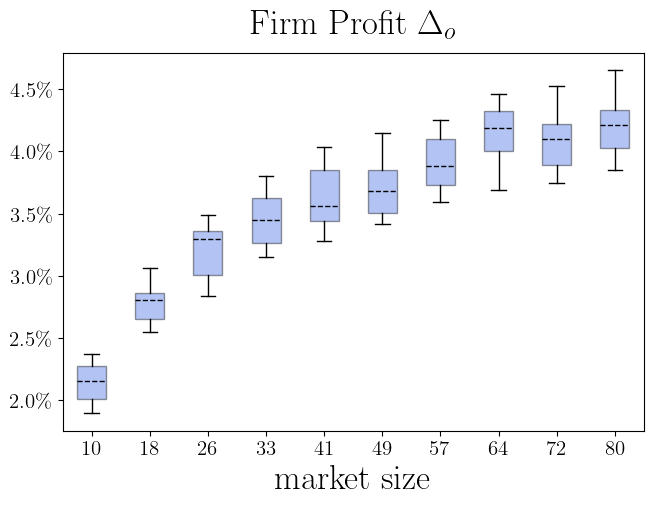} 
                \includegraphics[width=0.24\textwidth]{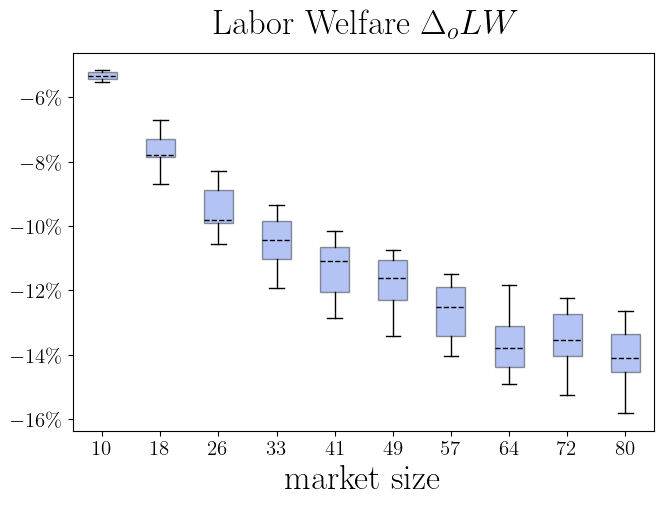} 
                \includegraphics[width=0.24\textwidth]{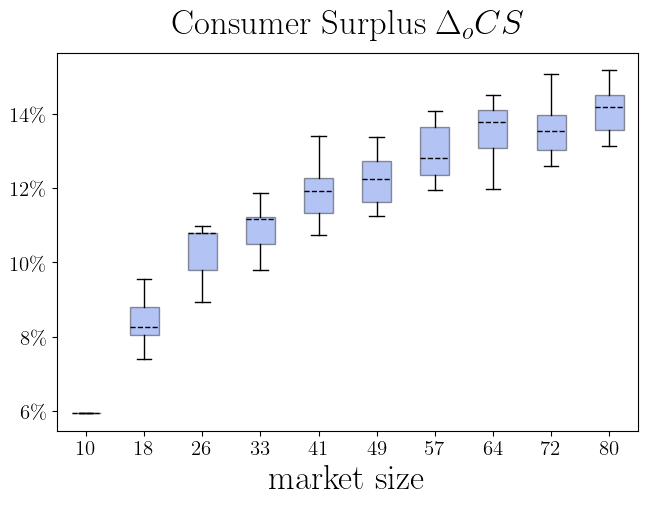} 
                \includegraphics[width=0.24\textwidth]{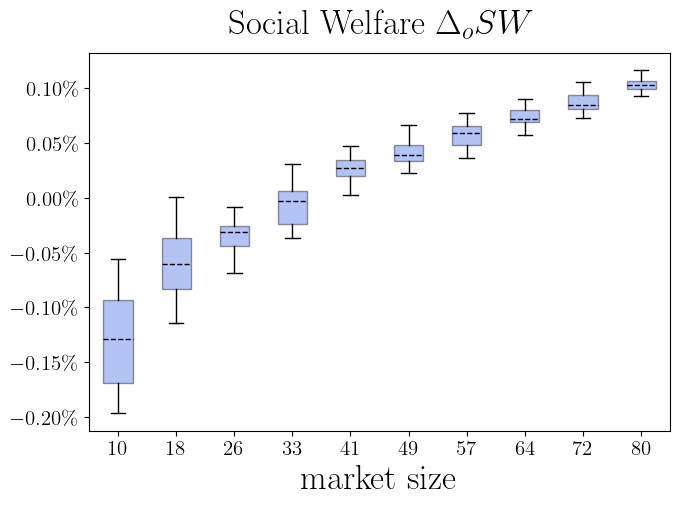} 
                \vfill
                \includegraphics[width=0.24\textwidth]{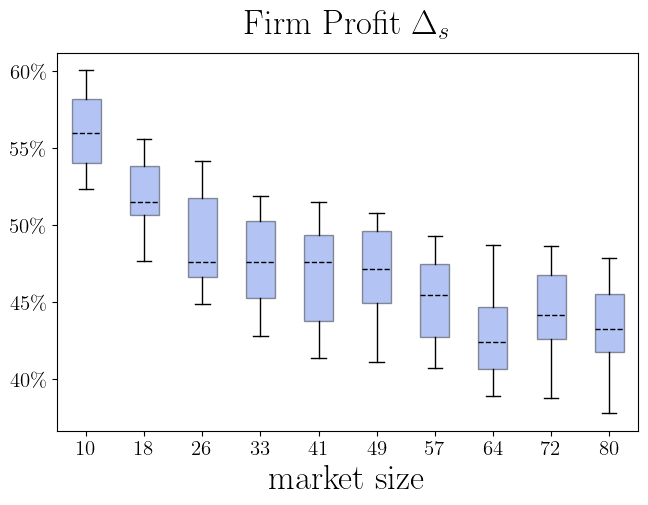}
                \includegraphics[width=0.24\textwidth]{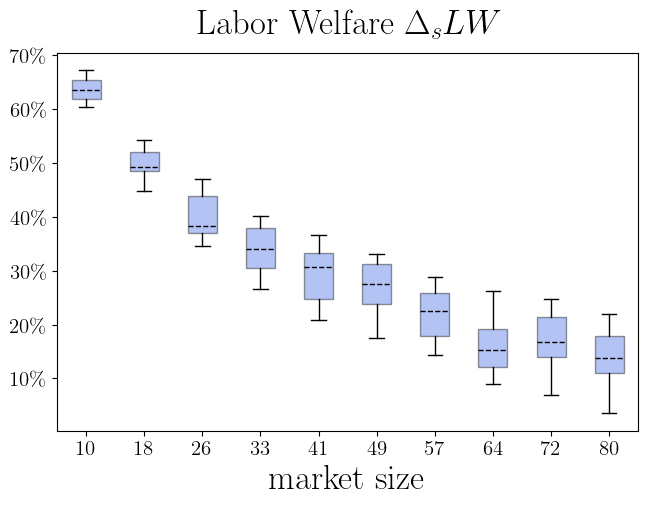} 
                \includegraphics[width=0.24\textwidth]{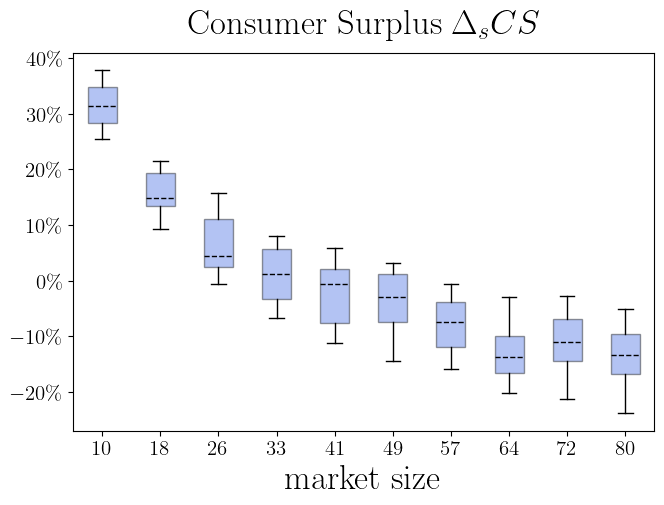} 
                \includegraphics[width=0.24\textwidth]{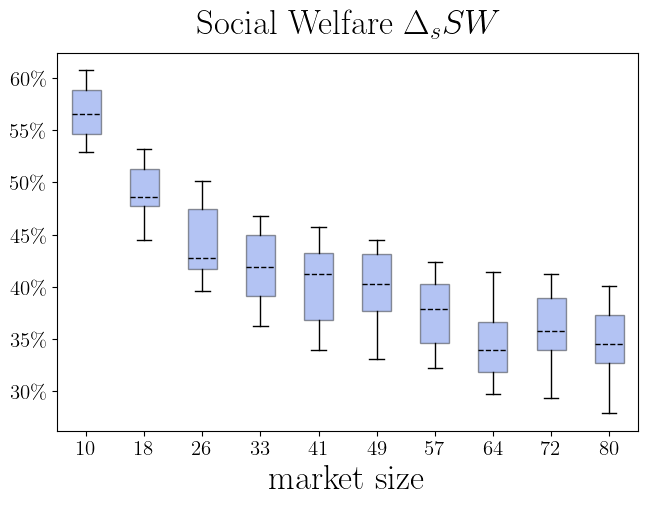}
            \end{minipage}
	}
	{\bf Value of service proliferation:  Impacts of market size  \label{fig-Lambda-K-2}}
	{This figure illustrates the interplay between market size $(\Lambda)$ and supply pool $(K)$ on the value of service proliferation. Here, we change these two parameters simultaneously with $K = \epsilon \Lambda, \epsilon \in [0.9,1.1],w_s = 0.6.$}
	\vspace{-0.1in}
\end{figure}

\begin{figure}[htbp]
	\vspace{-0.1in}
	\FIGURE
	{
            \begin{minipage}{0.95\textwidth}
                \includegraphics[width=0.24\textwidth]{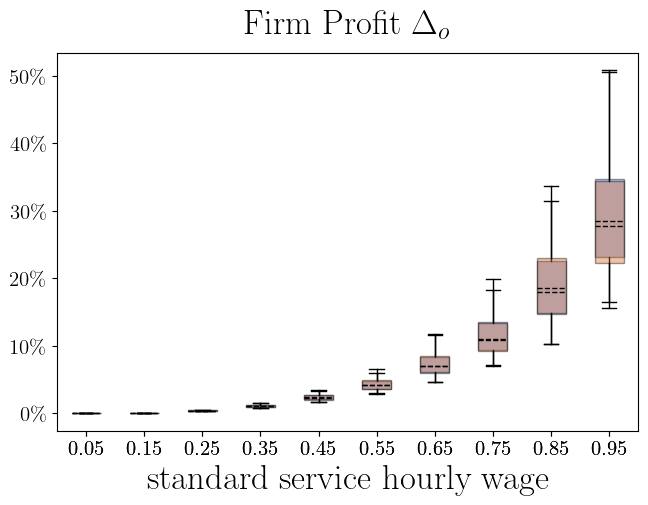} 
                \includegraphics[width=0.24\textwidth]{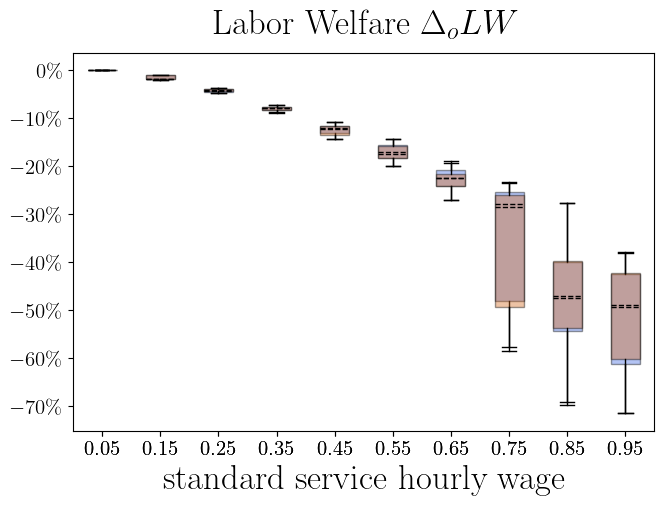} 
                \includegraphics[width=0.24\textwidth]{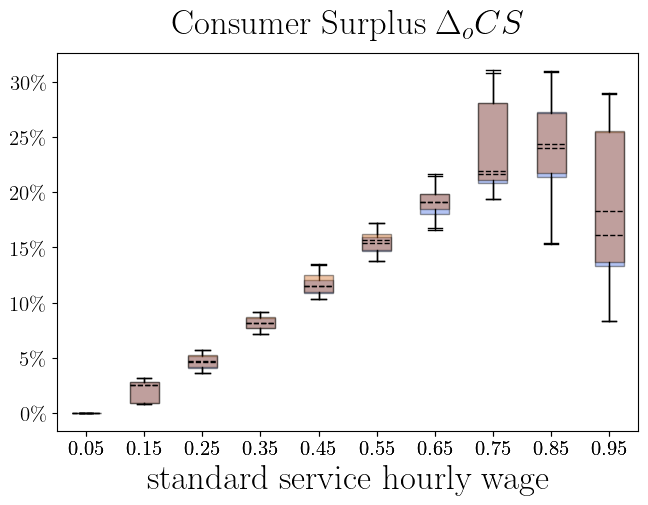} 
                \includegraphics[width=0.24\textwidth]{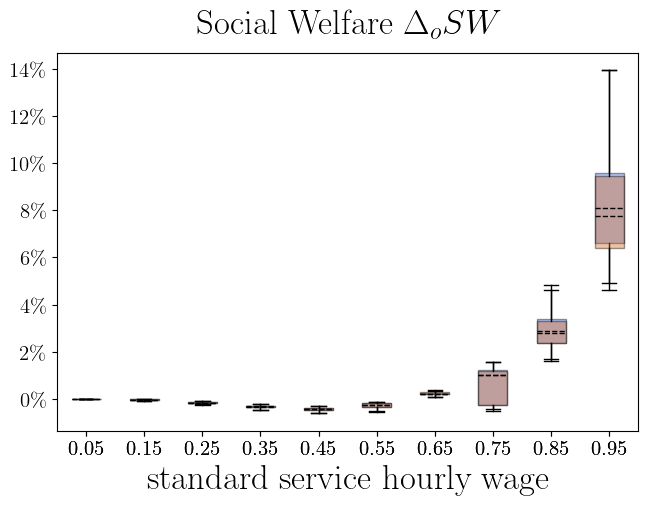} 
                \vfill
                \includegraphics[width=0.24\textwidth]{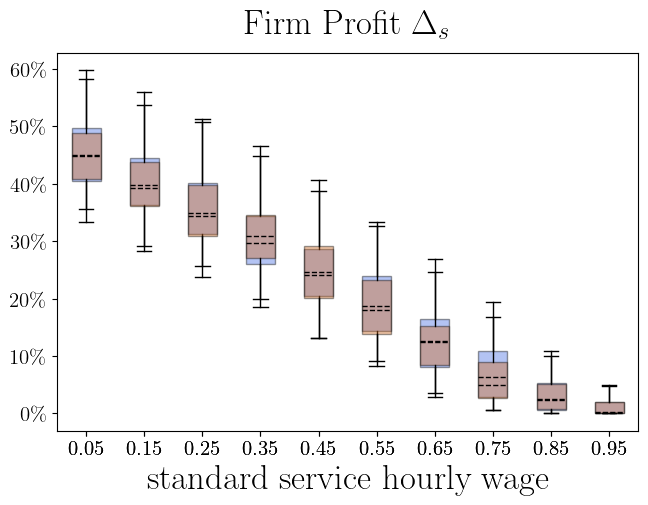} 
                \includegraphics[width=0.24\textwidth]{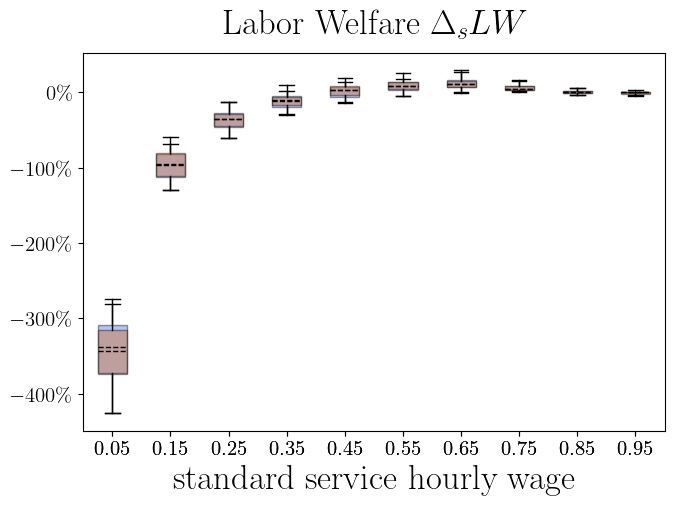} 
                \includegraphics[width=0.24\textwidth]{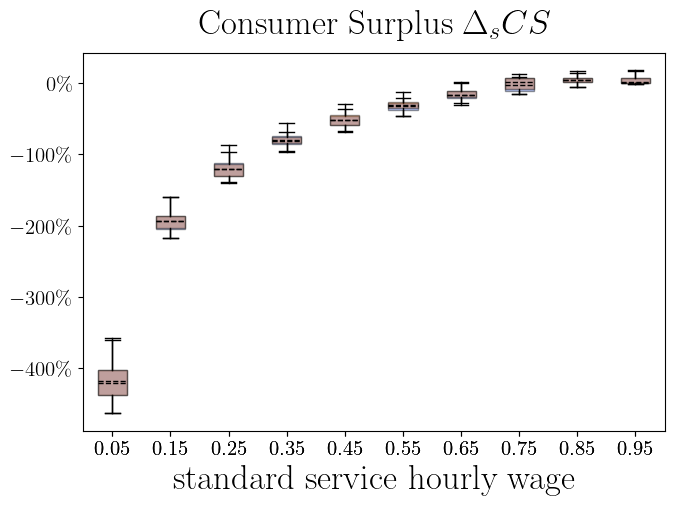}
                \includegraphics[width=0.24\textwidth]{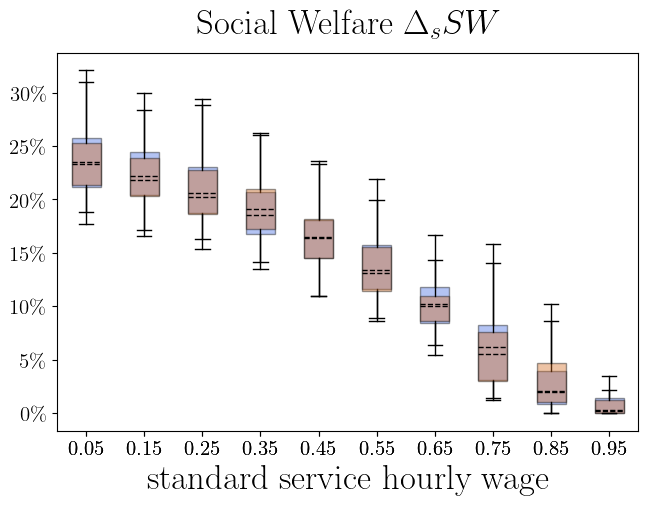}
            \end{minipage}
	}
	{\bf Value of service proliferation: Robustness of sample size
 \label{fig-value-size}}
	{This figure illustrates the robustness of the sample size on simulating the value of service proliferation. The blue boxplot is based on a sample size of 1000, and the orange boxplot is based on a smaller sample of 100 observations.}
	\vspace{-0.2in}
\end{figure}

\begin{figure}[htbp]
	\vspace{-0.1in}
	\FIGURE
	{
            \subfigure[~Standard vs on-demand]{
                \begin{minipage}{0.25\textwidth}
                    \includegraphics[width=0.96\textwidth]{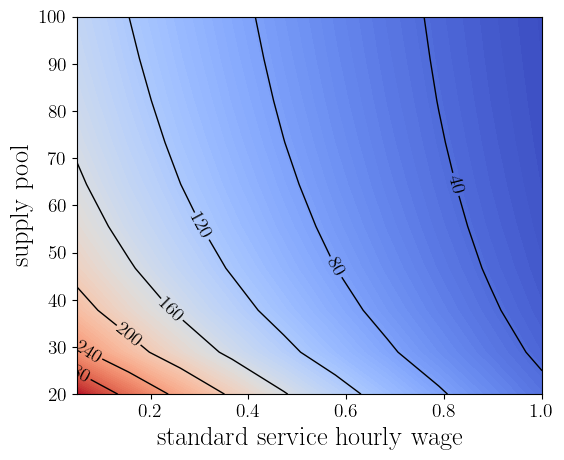} 
                    \vfill
                    \includegraphics[width=0.96\textwidth]{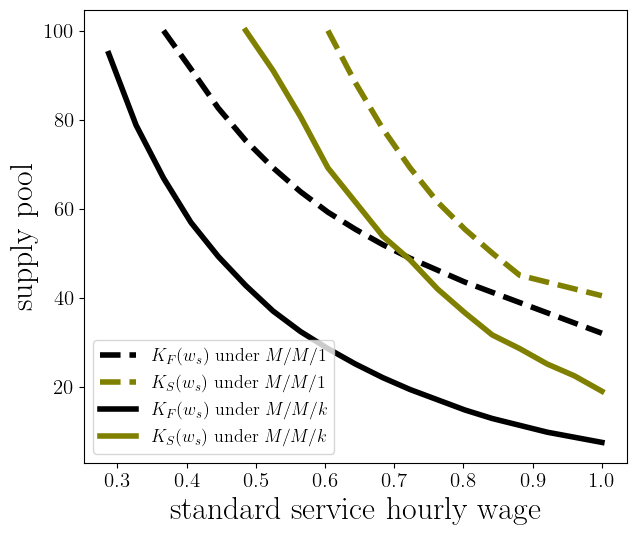} 
                \end{minipage}
            }
            \subfigure[~Labor welfare, consumer surplus and social welfare]{
                \begin{minipage}{0.73\textwidth}
                    \includegraphics[width=0.96\textwidth]{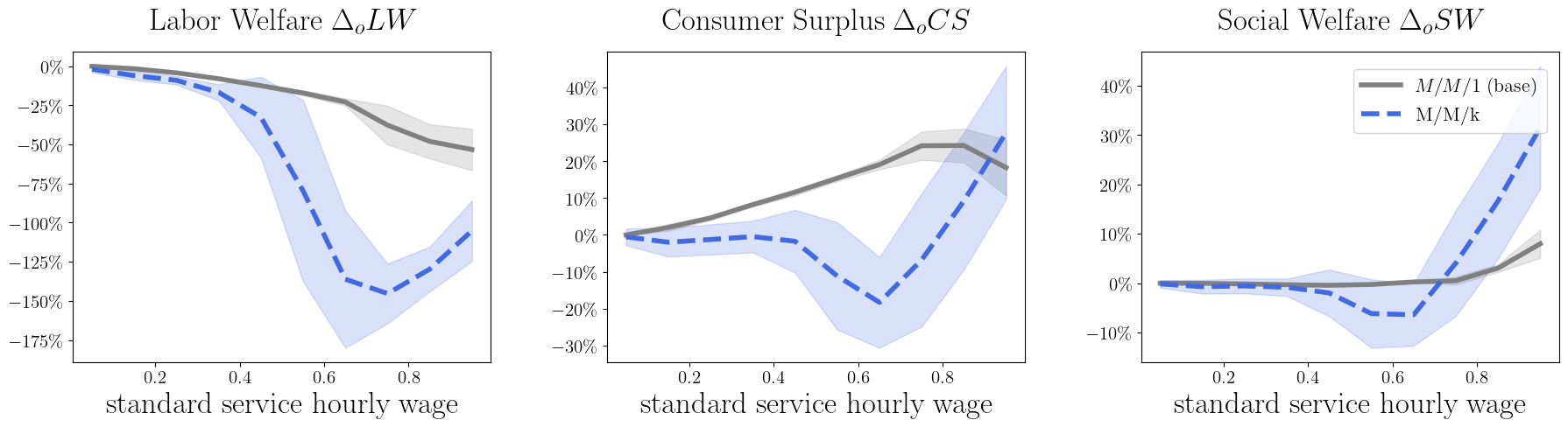} 
                    \vfill
                    \includegraphics[width=0.96\textwidth]{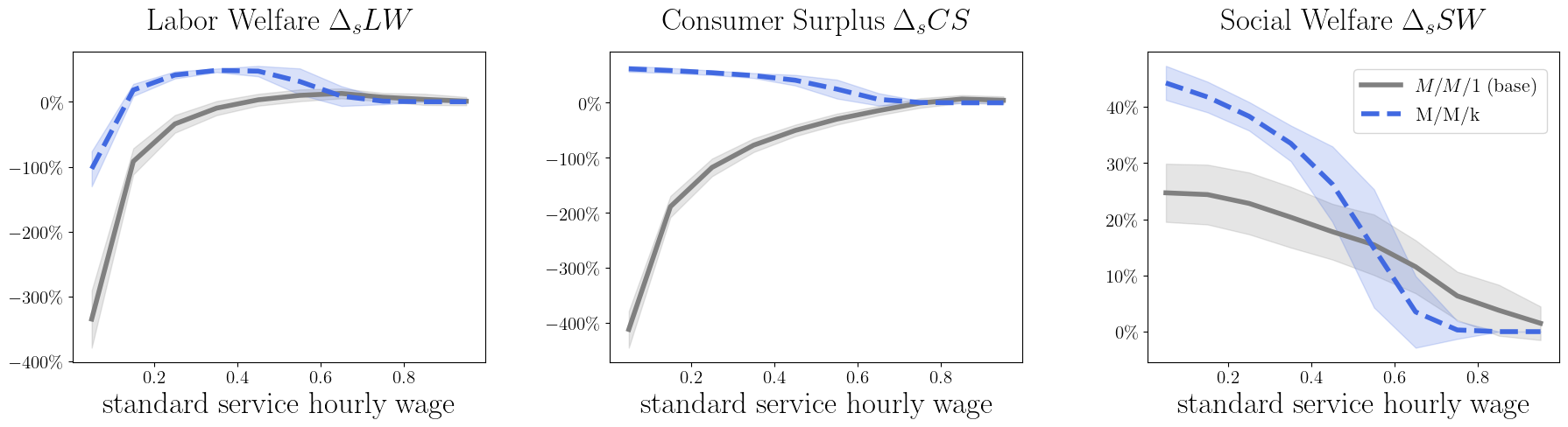} 
                \end{minipage}
            }
	}
	{Value of service proliferation: $M/M/k$ queue \label{fig-mmk-lw} }
	{This figure illustrates the impact of different service mechanisms when the lead time is evaluated by an $M/M/k$ queue. The left panel considers a single service mechanism with employees or contractors, analogous to Figure \ref{fig-ratio}, in which the heat map measures the effectiveness of employees, whereby a value of $<(>) 100$ means that the standard mechanism reduces (improves) the SP profit. The right panel depicts the value of service proliferation from the perspectives of labor welfare, consumer surplus, and total social welfare.} 
	\vspace{-0.2in}
\end{figure}

\begin{figure}[htbp]
	\vspace{-0.1in}
        \FIGURE
	{
        \begin{minipage}{0.7\textwidth}
            {
                \includegraphics[width=0.32\textwidth]{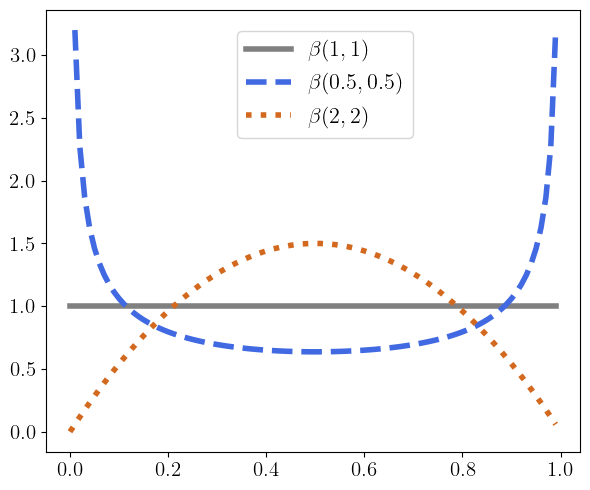}
                \includegraphics[width=0.32\textwidth]{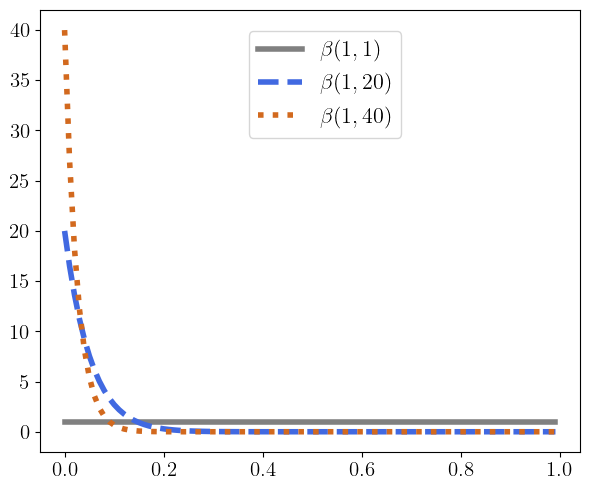}
                \includegraphics[width=0.32\textwidth]{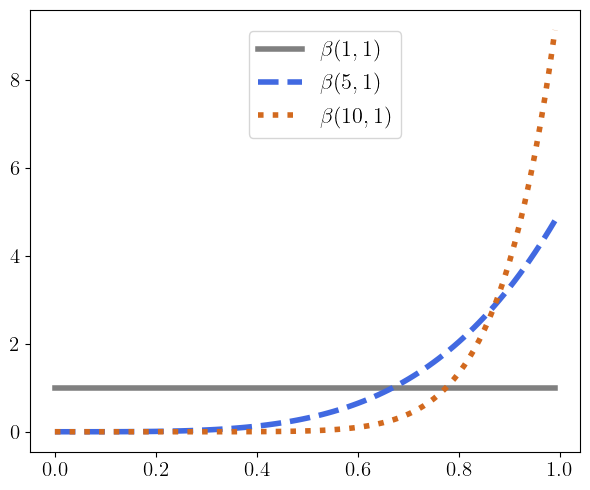}
            } 
        \end{minipage}
	}
	{Probability density functions for Beta distributions: Variability and skewness  \label{ec-fig-beta-pdf}}
	{This figure illustrates the variability and skewness of various Beta distributions used in this paper. The skewness is zero in the left panel (or is symmetric), positive in the middle panel, and negative in the right panel. The variability, measured by the coefficient of variation, for a symmetric Beta distribution $\beta(\alpha,\alpha)$ is given by $\frac{1}{\sqrt{2\alpha+1}}$.} 
	\vspace{-0.2in}
\end{figure}

\begin{figure}[htbp]
	\vspace{-0.1in}
	\FIGURE
	{
        \begin{minipage}{0.8\textwidth}
            \subfigure[~zero skewness]{
                \includegraphics[width=0.95\textwidth]{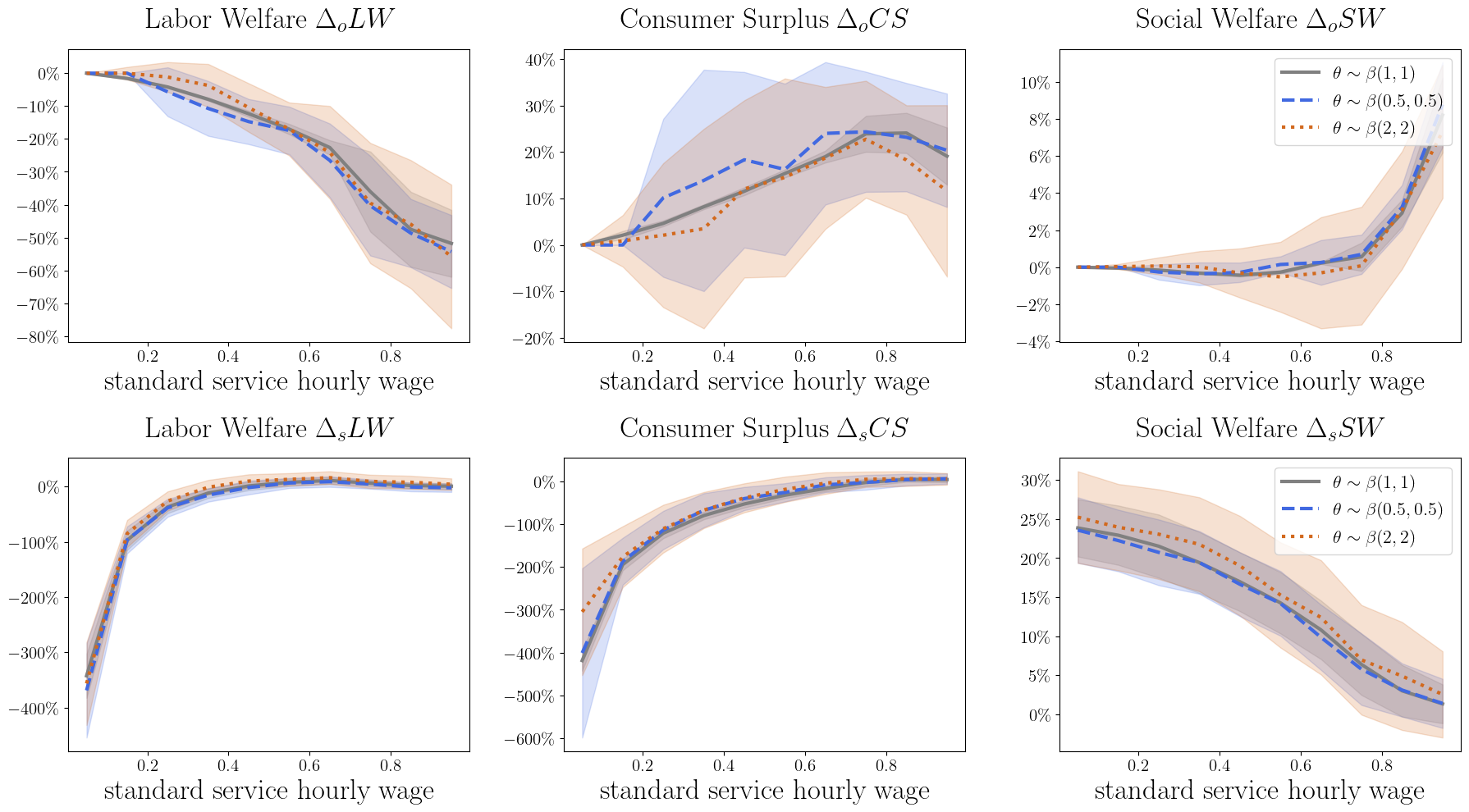} 
            }
            \vspace{-0.1in}
            \vfill
            \subfigure[~positive skewness]{
                    \includegraphics[width=0.95\textwidth]{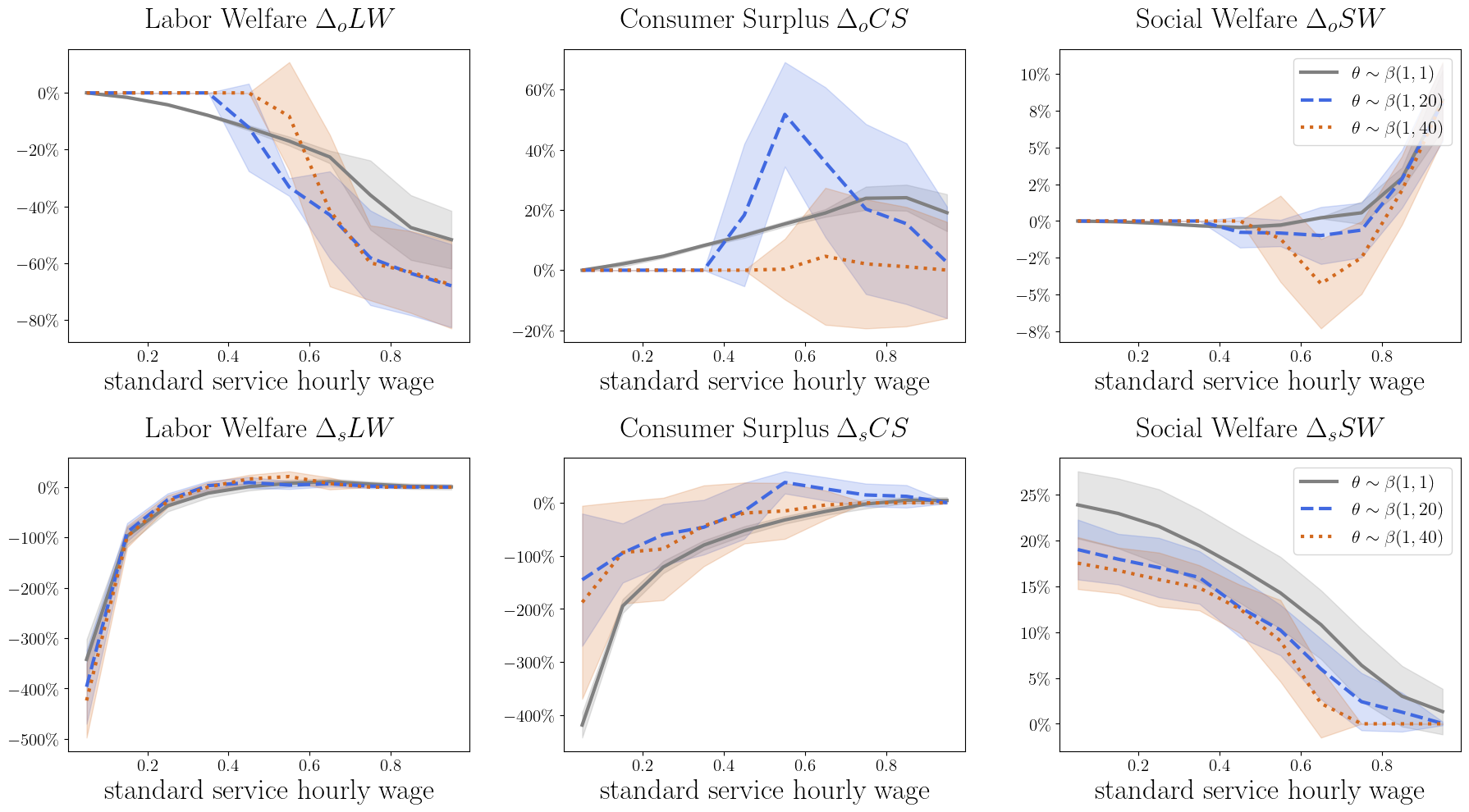} 
                }
            \vspace{-0.1in}
            \vfill
            \subfigure[~negative skewness]{
                    \includegraphics[width=0.95\textwidth]{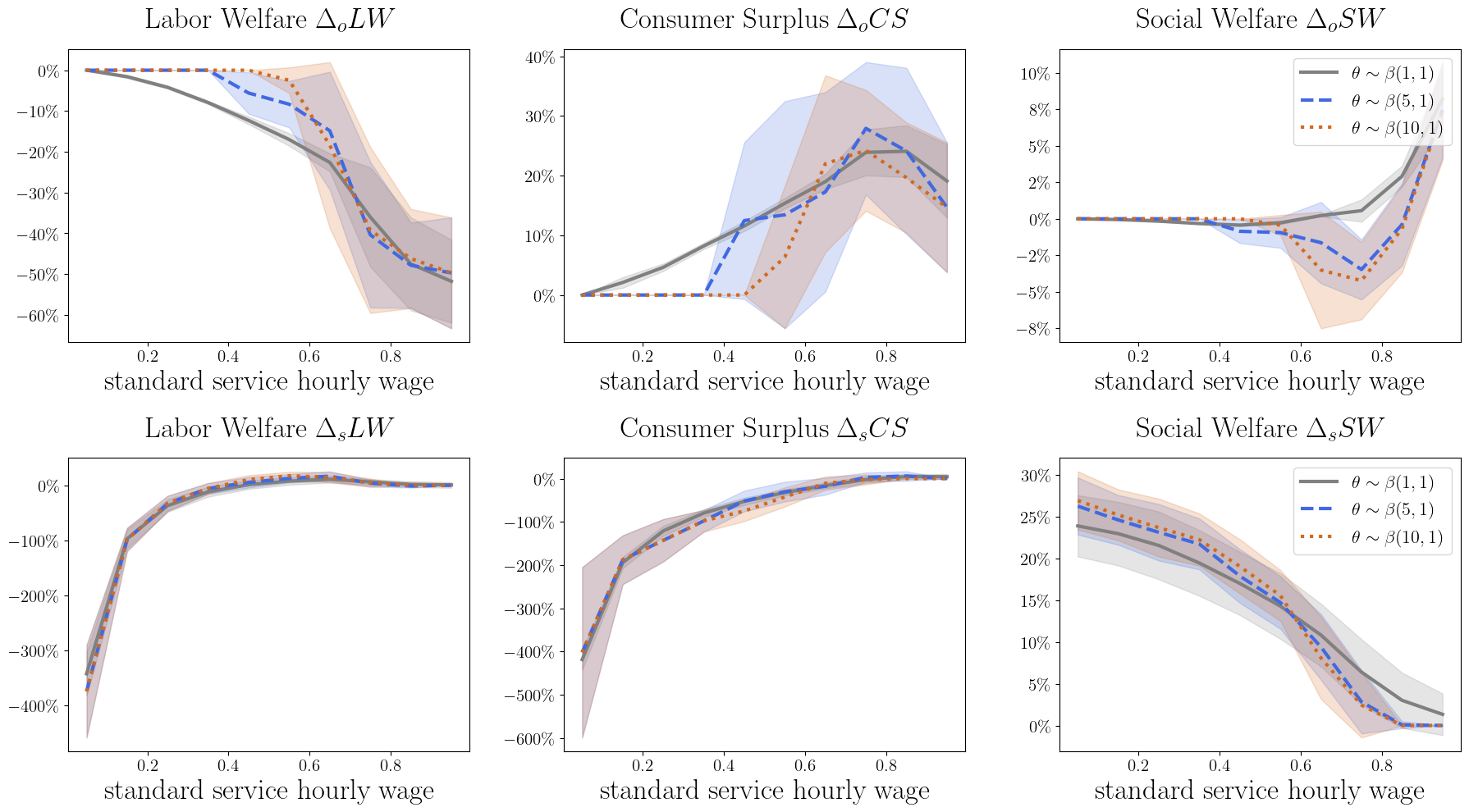} 
                }
        \end{minipage}
	}
	{Heterogeneity in customer waiting sensitivity \label{fig-robust-zero-lw}}
	{} 
	\vspace{-0.2in}
\end{figure}


\begin{figure}[htbp]
	\vspace{-0.1in}
	\FIGURE
	{
            \begin{minipage}{0.8\textwidth}
            {
                \includegraphics[width=0.32\textwidth]{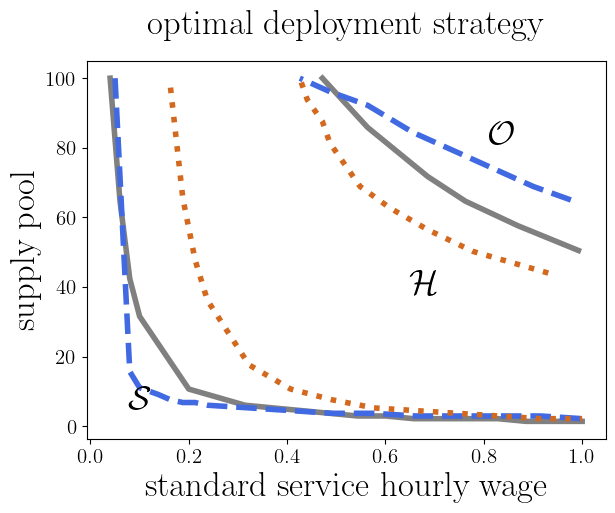} 
                \includegraphics[width=0.32\textwidth]{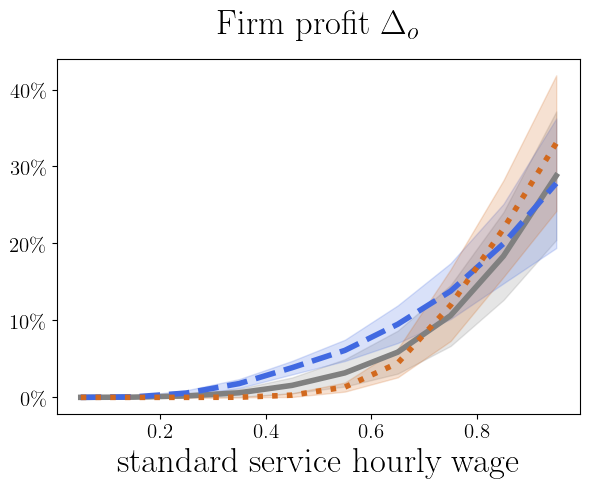} 
                \includegraphics[width=0.32\textwidth]{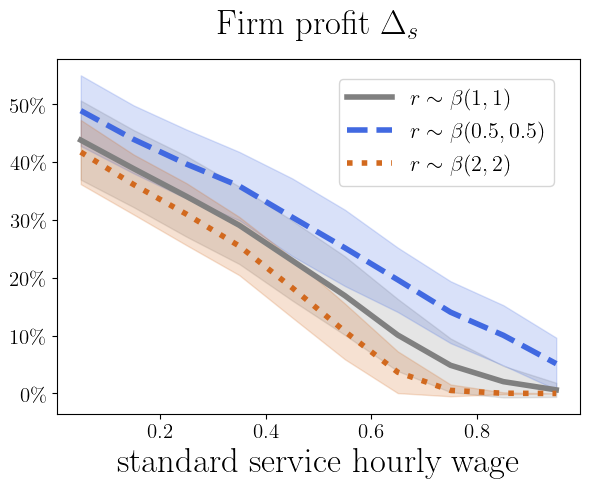}
                \vfill   
                \includegraphics[width=0.96\textwidth]{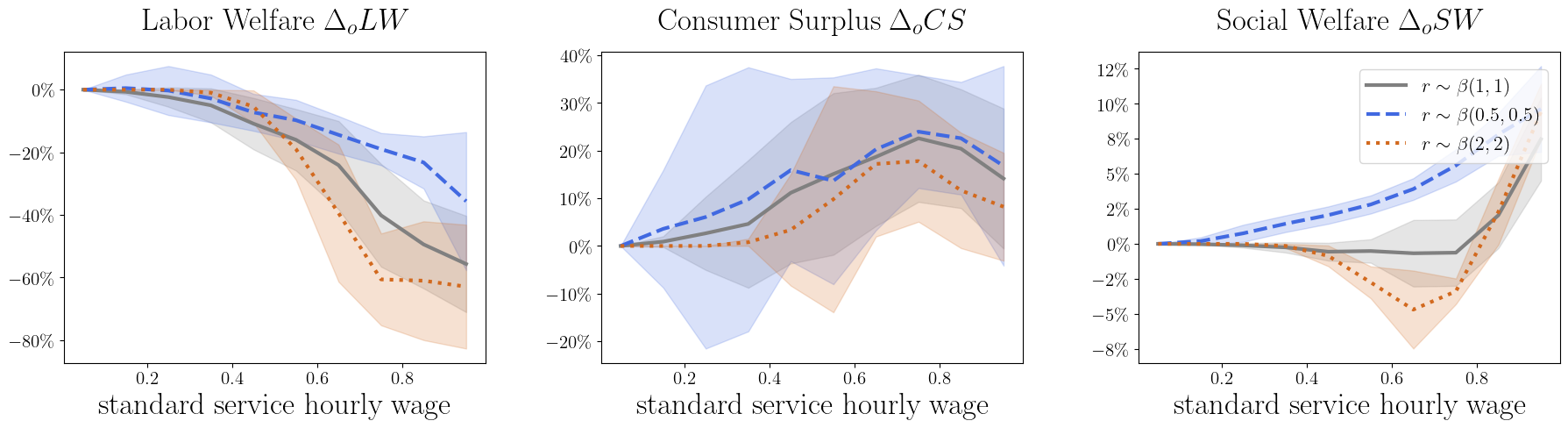} 
                \vfill        
                \includegraphics[width=0.96\textwidth]{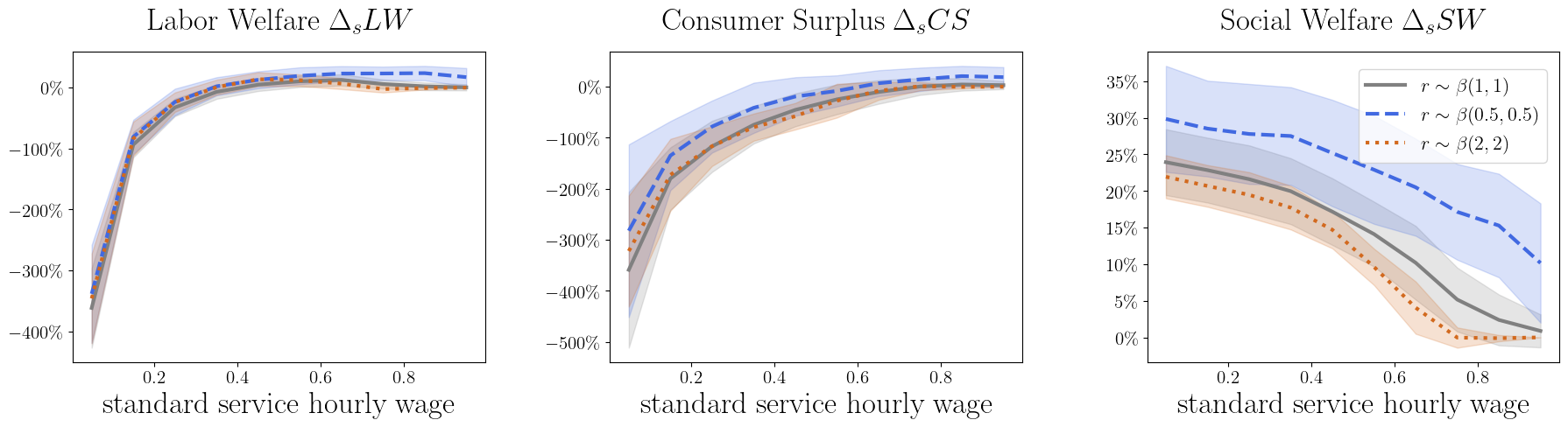}
            }           
            \end{minipage}
	}
	{\bf Heterogeneity in agent reservation rate: zero skewness  \label{fig-robust-r}}
	{}
	\vspace{-0.2in}
\end{figure}

\begin{figure}[htbp]
	\vspace{0in}
	\FIGURE
	{
        \begin{minipage}{0.8\textwidth}
            {
                \includegraphics[width=0.96\textwidth]{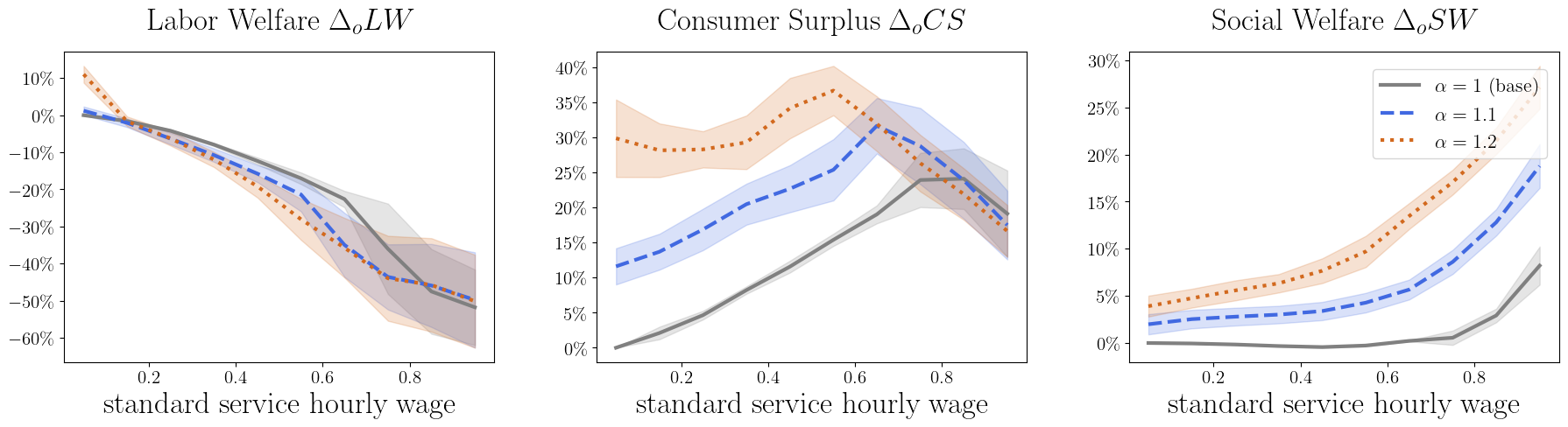} 
                \vfill
                \includegraphics[width=0.96\textwidth]{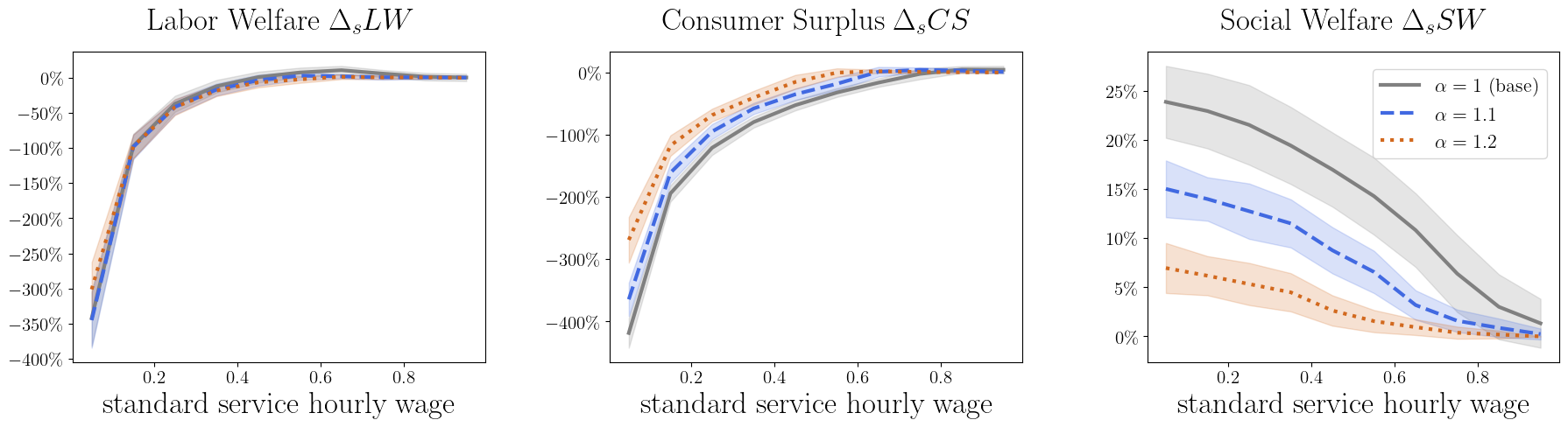}
            } 
        \end{minipage}
	}
	{\bf Heterogeneity in service quality \label{fig-diffV-lw} }
	{} 
	\vspace{-0.2in}
\end{figure}

\clearpage


\setcounter{page}{1}
\renewcommand{\thesection}{EC-\Alph{section}}
\setcounter{section}{0}

\begin{center}
	{\large Additional Online Supplement \\
		\vspace{.1in}
		\Large \PaperTitle
		\vspace{.1in}
	}
\end{center}

This note serves as the additional online supplement for the paper titled `` \PaperTitle". In this note, we provide the formal proofs for all supplemental results presented in Section \ref{ec-sec-supplements}.

\vspace{-0.1in}

\renewcommand{\theequation}{EC-\arabic{equation}}
\renewcommand{\thefigure}{EC-\arabic{figure}}
\renewcommand{\thetable}{EC-\arabic{table}}
\renewcommand{\thethm}{EC-\arabic{thm}}
\renewcommand{\thelem}{EC-\arabic{lem}}
\renewcommand{\theprop}{EC-\arabic{prop}}
\renewcommand{\thecor}{EC-\arabic{cor}}
\setcounter{equation}{0}
\setcounter{thm}{0}
\setcounter{lem}{0}
\setcounter{prop}{0}
\setcounter{cor}{0}
\setcounter{figure}{0}
\setcounter{table}{0}

\section{Supplemental Proofs}\label{ec-sec-supplemental-proofs}
\noindent
{\bf Proof of Lemma \ref{lem-subproblem}. }
We simplify the optimization problem \eqref{opt} by variable transformation. There are two cases, $W_s\leq W_o$ and $W_s\geq W_o$.  We  derive a one-to-one mapping between $(p_s,p_o,k_s,w_o)$ and $(\lambda_s,W_s,\lambda_o,W_o)$ in each case, and then rewrite $\pi^*$ in a simpler way by combing the two cases. 

Using \eqref{W} and \eqref{ko}, the cost functions become
\[C_s(k_s) = \frac{w_s}{\mu_s} (\frac{1}{W_s} + \lambda_s) ,~ 
C_o(w_o) =  \frac{1}{K\mu_o^2}(\frac{1}{W_o} + \lambda_o )^2.\] 
Given $(p_s, p_o, W_s,W_o)$, the effective arrival rate to the two channels are given by
\begin{eqnarray}
\lambda_s &= \Lambda {\sf Pr}(U_s \ge 0, U_s \ge U_o)  = \Lambda {\sf Pr} \left(\theta \leq \frac{V-p_s}{W_s}, p_s -p_o + (W_s - W_o) \theta \leq 0 \right), \nonumber \\
\lambda_o &= \Lambda {\sf Pr} (U_o \ge 0, U_o \ge U_s) =  \Lambda {\sf Pr} \left(\theta \leq \frac{V-p_o}{W_o}, p_s -p_o + (W_s - W_o) \theta \geq 0  \right). \label{eq-lambda}
\end{eqnarray}
We consider two cases for \eqref{opt}:  $W_o \geq W_s$  and  $W_o \leq W_s$. 

Case 1: $W_o \geq W_s$.  
It follows from $W_o \geq W_s$  that  \eqref{eq-lambda}  can be written as
\bea 
\lambda_o = \Lambda {\sf Pr} ( \theta \leq \frac{V-p_o}{W_o} , \theta \leq \frac{p_s - p_o}{W_o - W_s } ), \ \lambda_s = \Lambda {\sf Pr} ( \frac{p_s - p_o}{W_o - W_s } < \theta \leq  \frac{V-p_s}{W_s}).  \label{eq-lambda-ol}
\eea
Note that, if $W_s \leq W_o$, $p_o\leq p_s$, and $V\geq  \max\{p_o,p_s\} $, then  $\frac{V-p_s}{W_s}\geq  \frac{p_s- p_o }{W_o - W_s} $ is equivalent to $\frac{W_o}{W_s} \geq  \frac{V-p_o }{V-p_s} $, $\frac{p_s- p_o }{W_o - W_s} \geq \frac{V-p_o}{W_o}  $ is equivalent to $\frac{W_o}{W_s} \leq  \frac{V-p_o }{V-p_s} $. Therefore, there are three subcases:
\begin{itemize}
    \item If $p_s<p_o$, then \eqref{eq-lambda-ol} becomes $\lambda_o =0$ and \eqref{opt} reduces to \eqref{model-S}, which is System S. Hence the optimal profit in this subcase is $\pi^S$.
    \item If $p_o\leq p_s$ and $\frac{W_o}{W_s} \leq  \frac{V-p_o }{V-p_s} $, then \eqref{eq-lambda-ol} becomes  $\lambda_s=0$  and \eqref{opt} reduces to \eqref{model-O},   which is System O. Hence the optimal profit in this subcase is $\pi^O$.
    \item  If $p_o\leq p_s$ and $\frac{W_o}{W_s}  \geq \frac{V-p_o }{V-p_s} $, then \eqref{eq-lambda-ol} becomes
\[\lambda_o = \Lambda \frac{p_s- p_o }{W_o - W_s},~
\lambda_s =  \Lambda   \left( \frac{V-p_s}{W_s}-  \frac{p_s- p_o }{W_o - W_s}  \right),\]
which can be equivalently written as
\bea 
p_o = V - \frac{\lambda_o W_o + \lambda_s W_s}{\Lambda}, \ p_s = V - \frac{\lambda_o W_s + \lambda_s W_s}{\Lambda}. \label{eq-price-ol}
\eea
The optimization problem \eqref{opt} in this case becomes 
\begin{align} \label{model-H-ol} 
\pi_{ol}=
\left\{  \begin{array}{rl} \max & \pi(\lambda_s,W_s, \lambda_o,W_o) 
= V(\lambda_s + \lambda_o) - \frac{ \lambda_s^2 W_s + 2 \lambda_s \lambda_o W_s + \lambda_o^2 W_o  }{\Lambda} - \frac{w_s}{\mu_s} (\frac{1}{W_s} + \lambda_s) - \frac{ 1}{K\mu_o^2}(\frac{1}{W_o} + \lambda_o )^2   \\
\mbox{s.t.} &    \lambda_s,\lambda_o,W_s > 0, W_s \leq W_o,  \lambda_s+\lambda_o \leq \Lambda.
\end{array}
\right. ,
\end{align}
where the subscript in $\pi_{ol}$ indicates that on-demand service has a longer lead time. 
\end{itemize}

 Case 2: $W_o \leq W_s$. 
It follows from $W_o \leq W_s$ that \eqref{eq-lambda}  can be written as
\begin{align}
\lambda_s = \Lambda {\sf Pr} \left(\theta \leq \frac{V-p_s}{W_s},    \theta \leq  \frac{p_o- p_s }{W_s - W_o}\right),  ~
\lambda_o =  \Lambda {\sf Pr} \left(\frac{p_o- p_s }{W_s - W_o} \leq \theta \leq \frac{V-p_o}{W_o}   \right). \label{eq-lambda-sl}
\end{align} 
Note that, if $W_s\geq W_o$, $p_o\geq p_s$, and $V\geq  \max\{p_o,p_s\} $, then  $\frac{V-p_s}{W_s}\leq  \frac{p_o- p_s }{W_s - W_o} $ is equivalent to $\frac{W_o}{W_s} \geq  \frac{V-p_o }{V-p_s} $, $\frac{p_o- p_s }{W_s - W_o}    \leq \frac{V-p_o}{W_o}  $ is equivalent to $\frac{W_o}{W_s} \leq  \frac{V-p_o }{V-p_s} $. Therefore, we consider three subcases: 
\begin{itemize}
    \item If $p_o<p_s$, then \eqref{eq-lambda-sl} becomes $\lambda_s =0$  and  \eqref{opt} reduces to \eqref{model-O}, which is System O. Hence the optimal profit in this subcase is $\pi^O$.
    \item If $p_o\geq p_s$ and $\frac{W_o}{W_s} \geq  \frac{V-p_o }{V-p_s} $, then \eqref{eq-lambda-sl} becomes 
  $\lambda_o=0$ and \eqref{opt} reduces to \eqref{model-S}, which is System S. Hence the optimal profit in this subcase is $\pi^S$.
    \item If $p_o\geq p_s$ and $\frac{W_o}{W_s}  \leq \frac{V-p_o }{V-p_s} $, then \eqref{eq-lambda-sl} becomes
\[\lambda_s= \Lambda \frac{p_o- p_s }{W_s - W_o},~
\lambda_o =  \Lambda   \left( \frac{V-p_o}{W_o}-  \frac{p_o- p_s }{W_s - W_o}       \right),\]
which can be equivalently written as
\bea 
p_s = V - \frac{\lambda_s W_s + \lambda_o W_o}{\Lambda}, \ p_o = V - \frac{\lambda_s W_o + \lambda_o W_o}{\Lambda}. \label{eq-price-sl}
\eea
The optimization problem \eqref{opt} in this case is 
\begin{align}\label{model-H-sl} 
\pi_{sl}=
\left\{  \begin{array}{rl} \max & \pi(\lambda_s,W_s, \lambda_o,W_o) 
= V(\lambda_s + \lambda_o) - \frac{ \lambda_s^2 W_s + 2 \lambda_s \lambda_o W_o + \lambda_o^2 W_o  }{\Lambda} - \frac{w_s}{\mu_s} (\frac{1}{W_s} + \lambda_s) - \frac{ 1}{K\mu_o^2}(\frac{1}{W_o} + \lambda_o )^2   \\
\mbox{s.t.} &    \lambda_s,\lambda_o,W_o > 0,W_s\geq W_o, \lambda_s+\lambda_o \leq \Lambda.
\end{array}
\right. 
\end{align} 
where the subscript in $\pi_{sl}$ indicates that standard service has a longer lead time. 
\end{itemize}

We combine the optimization problems \eqref{model-H-ol} and \eqref{model-H-sl} by defining the following problem:
\beq
\pi^H_{\lambda_s, W_s, \lambda_o,W_o} = \left\{  \begin{array}{rl}
\max_{} & \pi(\lambda_s, W_s, \lambda_o,W_o) = V(\lambda_s + \lambda_o) - \frac{ \lambda_s^2 W_s + 2 \lambda_s \lambda_o \min\{W_o,W_s\} + \lambda_o^2 W_o  }{\Lambda} - \frac{w_s}{\mu_s} (\frac{1}{W_s} + \lambda_s) - \frac{ 1}{K\mu_o^2}(\frac{1}{W_o} + \lambda_o )^2. \\
\mbox{s.t.} & \lambda_s, \lambda_o, W_s,W_o > 0, \ \lambda_s + \lambda_o \leq \Lambda. 
\end{array}\right. 
\eeq
That is, $\max\{\pi_{ol}, \pi_{sl}\} = \pi^H_{\lambda_s, W_s, \lambda_o,W_o}$.  Based on above analysis, we can summarize that the optimal profit of optimization problem \eqref{opt} is given by  $\pi^* = \max\{ \pi^S, \pi^O, \pi^H_{\lambda_s, W_s, \lambda_o,W_o} \}$. 
\qed

\noindent
{\bf Proof of Lemma \ref{lem-solve-ls}.}
In this Lemma, we aim to solve for the optimal value of $\lambda_s$   in \eqref{model-H} for any given $(W_s, \lambda_o, W_o)$. Let $\mathcal{D} = \{( W_s, \lambda_o, W_o )| W_s, W_o > 0, 0\leq \lambda_o \leq \Lambda \}$, and 
\beq
\mathcal{D}_1 &=& \{(W_s, \lambda_o, W_o )|  0\leq \lambda_o \leq \Lambda,  W_o \geq W_s > 0, W_s \geq   \frac{\Lambda (V \mu_s - w_s)}{2\mu_s \lambda_o}\}, 
\\
\mathcal{D}_2 &=& \{(W_s, \lambda_o, W_o )|  0\leq \lambda_o \leq \Lambda,  W_o \geq W_s > 0, \frac{V \mu_s - w_s}{2\mu_s } \leq W_s \leq \frac{\Lambda ( V \mu_s -w_s)}{2\mu_s \lambda_o}\}, 
\\
\mathcal{D}_3 &=& \{(W_s, \lambda_o, W_o )|   0\leq \lambda_o \leq \Lambda,   W_o \geq W_s>0 , W_s < \frac{V \mu_s - w_s}{2\mu_s }\}, 
\\
\mathcal{D}_4 &=& \{(W_s, \lambda_o, W_o )|  0\leq \lambda_o \leq \Lambda,   W_s \geq  W_o > 0, V\Lambda - \frac{\Lambda w_s}{\mu_s } - 2 \lambda_o W_o \leq 0\}, 
\\
\mathcal{D}_5 &=& \{(W_s, \lambda_o, W_o )|   0\leq \lambda_o \leq \Lambda, W_s \geq  W_o > 0,  V\Lambda - \frac{\Lambda w_s}{\mu_s } - 2 \lambda_o W_o \geq 0, W_s \geq \frac{V\Lambda - \frac{\Lambda w_s}{\mu_s } - 2 \lambda_o W_o }{2(\Lambda - \lambda_o) } \}, 
\\
\mathcal{D}_6 &=& \{(W_s, \lambda_o, W_o )|  0\leq \lambda_o \leq \Lambda,  W_s \geq  W_o > 0,  V\Lambda - \frac{\Lambda w_s}{\mu_s } - 2 \lambda_o W_o \geq 0, W_s \leq \frac{V\Lambda - \frac{\Lambda w_s}{\mu_s } - 2 \lambda_o W_o }{2(\Lambda - \lambda_o) } \}.
\eeq
Note that $\cup_{i=1}^3 \mathcal{D}_i = \{( W_s, \lambda_o, W_o )| 0\leq  \lambda_o \leq \Lambda , W_o \geq W_s > 0 \}$, $\cup_{i=4}^6 \mathcal{D}_i = \{( W_s, \lambda_o, W_o )| 0\leq  \lambda_o \leq \Lambda , W_s \geq W_o >0 \}$,  and
$\cup_{i=1}^6 \mathcal{D}_i = \mathcal{D}$, see Figure \ref{fig-solve-ls} for the illustration in $(W_s,W_o)$ space. 
\begin{figure}[!hbt]
	\vspace{-0.1in}
	\FIGURE
	{
            \includegraphics[width=0.45\textwidth]{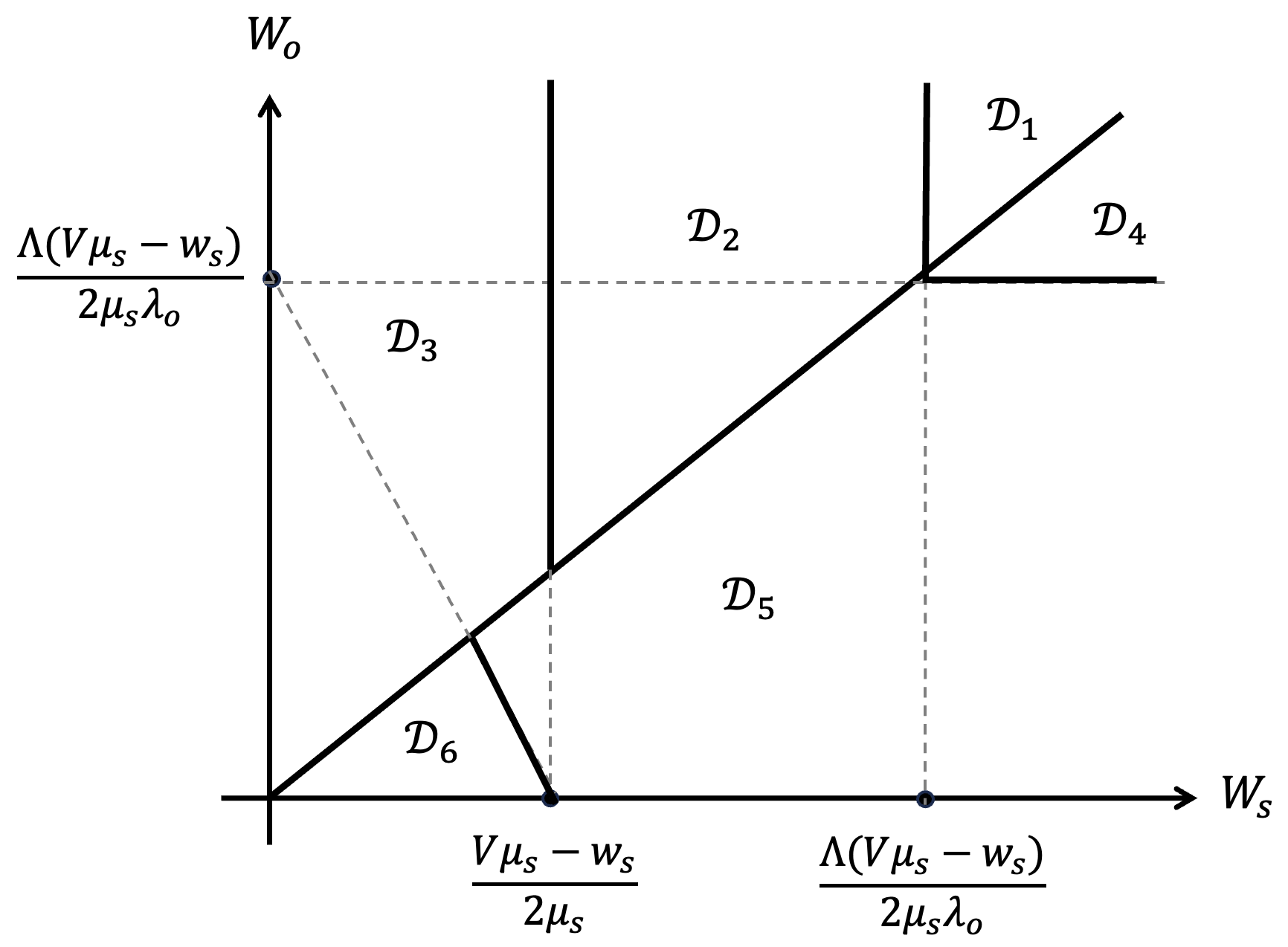} 
	}
	{\bf Illustration of $\mathcal{D}_i(i=1,2,...,6)$ in $(W_s,W_o)$ space \label{fig-solve-ls}}
	{}
	\vspace{-0.2in}
\end{figure}



Consistent with Lemma \ref{lem-subproblem}, we consider two cases: $W_o \geq W_s$ and $W_o \leq W_s$. 

Case 1: $W_o \geq W_s$. Then the objective function in \eqref{model-H}
can be rewritten as $\pi(\lambda_s,W_s, \lambda_o,W_o )  = \pi_2(\lambda_o,W_o) + V \lambda_s - \frac{ \lambda_s^2 W_s + 2 \lambda_s \lambda_o W_s  }{\Lambda} - \frac{w_s}{\mu_s} (\frac{1}{W_s} + \lambda_s)$. For fixed $(W_s,\lambda_o,W_o)$, the objective function is concave in $\lambda_s$. 
Solving the first order condition, i.e., $\frac{\partial \pi(\lambda_s,W_s, \lambda_o,W_o)}{\partial \lambda_s}|_{\bar \lambda_s(W_s, \lambda_o,W_o)}=0$, we derive $\bar \lambda_s(W_s, \lambda_o,W_o)= \frac{\Lambda (V \mu_s - w_s)}{2\mu_s W_s} - \lambda_o$. Since the feasible region of $\lambda_s$ is $[0,\Lambda - \lambda_o]$, the best $\lambda_s$ has three possibilities.
\begin{itemize}
    \item $\lambda_s = 0$ if $W_s \geq \frac{\Lambda (V \mu_s - w_s)}{2\mu_s \lambda_o}$, i.e., $(W_s, \lambda_o,W_o) \in \mathcal{D}_1$. Then the optimal profit in this subcase is   $\pi^O$.
    \item $\lambda_s= \bar \lambda_s(W_s,\lambda_o,W_o)$ if $\frac{V \mu_s - w_s}{2\mu_s } \leq W_s \leq \frac{\Lambda ( V \mu_s -w_s)}{2\mu_s \lambda_o}$, i.e., $(W_s, \lambda_o,W_o) \in \mathcal{D}_2$. In this case, the market is partially covered. Let 
\bea \label{model-olp}
\pi_{ol,p} =\left\{  \begin{array}{rl}
\max_{} & \pi( \lambda_s, W_s, \lambda_o,W_o)\big|_{ \{ \lambda_s = \bar \lambda_s(W_s)\} } = \pi_2(\lambda_o ,W_o) + \frac{(V\Lambda  - \frac{\Lambda  w_s}{\mu_s}- 2\lambda_oW_o)^2 - \frac{4\Lambda w_s}{\mu_s} }{4\Lambda W_s} \\
\mbox{s.t.} &  0 \leq \lambda_o \leq \Lambda, \frac{V \mu_s -w_s}{2\mu_s } \leq W_s \leq \min\{ \frac{\Lambda ( V\mu_s-w_s)}{2\mu_s \lambda_o},W_o\}. 
\end{array}\right.
\eea
    \item $\lambda_s = \Lambda - \lambda_o$ if $W_s \leq \frac{V \mu_s - w_s}{2\mu_s }$, i.e., $(W_s, \lambda_o,W_o) \in \mathcal{D}_3$.  In this case, the market is fully covered. Let
\bea \label{model-olf}
\pi_{ol,f} =\left\{  \begin{array}{rl}
\max_{} &\pi( \lambda_s, W_s, \lambda_o,W_o)\big|_{ \{ \lambda_s =\Lambda - \lambda_o\}} \\
& = \pi_2(\lambda_o ,W_o) + V(\Lambda - \lambda_o) - \frac{(\Lambda -\lambda_o)^2 W_s}{\Lambda} - \frac{2\lambda_o (\Lambda - \lambda_o) W_s}{\Lambda} - \frac{w_s}{\mu_s}(\frac{1}{W_s} + \Lambda - \lambda_o). \\
\mbox{s.t.} & 0 \leq  \lambda_o \leq \Lambda, W_s < \min\{ \frac{ V\mu_s-w_s}{2\mu_s },W_o\}. 
\end{array}\right.
\eea    
\end{itemize}

Next we show  
\begin{align}
\pi_{ol,p} \leq\max \{\pi^O,\pi_{ol,f}\}. \label{piolp}
\end{align}
Note that in \eqref{model-olp}, the feasible region of $W_s$ is empty if $W_o < \frac{V \mu_s - w_s}{2\mu_s }$, so we consider $W_o \geq \frac{V \mu_s - w_s}{2\mu_s }$.
Since $\frac{\partial^2 \pi(\bar \lambda_s(W_s), W_s,\lambda_o,W_o) }{\partial W_s^2} = \frac{\Lambda(V \mu_s - w_s)^2 - 4\mu_s w_s}{2\mu_s^2 W_s^3} = \frac{\pi_1(\Lambda) ( \pi_1(\Lambda) +  4\sqrt{ \frac{w_s\Lambda}{\mu_s}})}{2\mu_s \Lambda W_s^3}\geq 0$, 
we have $\pi(\bar \lambda_s(W_s), W_s,\lambda_o,W_o)$ is convex in $W_s$ and the best $W_s $ is either $\frac{V\mu_s  -w_s}{2\mu_s}$,   $\frac{\Lambda ( V\mu_s -w_s)}{2\mu_s \lambda_o}$, or $W_o$.
Suppose $W_s = \frac{\Lambda ( V\mu_s -w_s)}{2\mu_s \lambda_o} $ if $W_o \geq \frac{\Lambda ( V\mu_s -w_s)}{2\mu_s \lambda_o} $, then  $ \lambda_s = 0$ is optimal, so that $\pi_{ol,p} \leq \pi^O$. 
Suppose $W_s = W_o$ if $W_o \leq \frac{\Lambda ( V\mu_s -w_s)}{2\mu_s \lambda_o} $, that is, standard and on-demand services have identical lead time. In this situation, for any fixed $(\lambda_o,W_o)$, 
\[  
\pi_{ol,p} =\pi( \lambda_s, W_s, \lambda_o,W_o)\big|_{ \{ \lambda_s = \bar \lambda_s(W_s), W_s = W_o\} } =  \pi_2(\lambda_o ,W_o) + \frac{(V\Lambda - \frac{\Lambda w_s}{\mu_s} - 2\lambda_o W_o )^2 - \frac{4\Lambda w_s}{\mu_s }}{4\Lambda W_o} \equiv A_e.
\]
Suppose $W_s = \frac{V\mu_s -w_s}{2\mu_s} $. That is, the best $\lambda_s$ is $\Lambda - \lambda_o$. Then for any fixed $(\lambda_o,W_o)$, the profit is 
\[\pi_{ol,p} =  \pi( \lambda_s, W_s, \lambda_o,W_o)\big|_{ \{ \lambda_s =\Lambda - \lambda_o, W_s = \frac{V\mu_s -w_s}{2\mu_s}\} } = \pi_2(\lambda_o,W_o) + \frac{(V\mu_s - w_s)}{2\Lambda \mu_s}(\Lambda - \lambda_o)^2 - \frac{2w_s }{V\mu_s -w_s} \equiv A_1. \]
Since  $(\lambda_s,W_s)=(\Lambda -\lambda_o,\frac{V\mu_s -w_s}{2\mu_s} )$ is a feasible solution of $\pi_{ol,f}$, we have $A_1 \leq  \pi_{ol,f}$. Therefore, to show $\pi_{ol,p} \leq \max \{\pi^O,\pi_{ol,f}\}$, it suffices to show 
\begin{align}
A_e \leq \max \{\pi^O,A_1\}  .\label{20240725-1}
\end{align}
There are two cases: $A_1 \geq \pi_2(\lambda_o,W_o)$, and $A_1\leq \pi_2(\lambda_o,W_o)$. If $A_1\geq \pi_2(\lambda_o,W_o)$, then we have   
\begin{align} 
A_1 - A_e &= \frac{2\mu_s W_o +w_s - V\mu_s}{4\Lambda \mu_s^2 W_o (V\mu_s -w_s)} \left[ -4 \Lambda \mu_s w_s - 2\mu_s \lambda_o^2  W_o (V\mu_s -w_s) + \Lambda^2 (V\mu_s -w_s)^2\right] \nonumber\\
& \geq \frac{2\mu_s W_o +w_s - V\mu_s}{4\Lambda \mu_s^2 W_o (V\mu_s -w_s)} \left[ -4 \Lambda \mu_s w_s - \lambda_o \Lambda (V\mu_s -w_s)^2 + \Lambda^2 (V\mu_s -w_s)^2 \right] \nonumber \\
& = \frac{2\mu_s W_o +w_s - V\mu_s}{4 \mu_s^2 W_o (V\mu_s -w_s)} \left[ (\Lambda - \lambda_o) (V\mu_s -w_s)^2 -4 \mu_s w_s \right], \label{65}
\end{align} 
where the inequality is from $\frac{(V \mu_s - w_s)}{2\mu_s} \leq  W_o$ and $ W_o  \leq \frac{\Lambda (V \mu_s - w_s)}{2 \mu_s \lambda_o}$. Note that 
$(\Lambda - \lambda_o) (V\mu_s -w_s)^2 -4 \mu_s w_s  = \frac{1}{\Lambda - \lambda_o} \left[ (\Lambda - \lambda_o)^2 (V\mu_s -w_s)^2 -4 \mu_s w_s(\Lambda - \lambda_o)  \right] \geq \frac{1}{\Lambda - \lambda_o} \left[ (\Lambda - \lambda_o)^2 (V\mu_s -w_s)^2 -4 \mu_s w_s \Lambda \right]\geq 0$, where the last inequality is due to $A_1 \geq \pi_2(\lambda_o,W_o)$. So $\eqref{65}\geq 0$. If $A_1\leq \pi_2(\lambda_o,W_o)$, i.e.,  $\frac{(V\mu_s - w_s)}{2\Lambda \mu_s}(\Lambda - \lambda_o)^2 - \frac{2w_s }{V\mu_s -w_s} \leq 0$, then  we have $\frac{V\mu_s - w_s}{\mu_s }\lambda_o \geq V\Lambda - \frac{w_s \Lambda}{\mu_s} - 2\sqrt{\frac{w_s \Lambda}{\mu_s}}$. Together with $W_o \geq \frac{(V \mu_s - w_s)}{2\mu_s}$, we have $2\lambda_o W_o \geq V\Lambda - \frac{w_s \Lambda}{\mu_s} - 2\sqrt{\frac{w_s \Lambda}{\mu_s}}$, which implies $A_e \leq \pi_2(\lambda_o,W_o)$. So $A_e$ is worse than $\pi^O$. As such, $\pi_{ol,p} \leq \max \{\pi^O,\pi_{ol,f}\}$. 


Case 2: $W_o \leq W_s$. 
Then the objective function in \eqref{model-H}
can be rewritten as $\pi(\lambda_s,W_s, \lambda_o,W_o ) = \pi_2(\lambda_o,W_o) + V \lambda_s - \frac{ \lambda_s^2 W_s + 2 \lambda_s \lambda_o W_o  }{\Lambda} - \frac{w_s}{\mu_s} (\frac{1}{W_s} + \lambda_s)$. 
For fixed $(W_s,\lambda_o,W_o)$, the objective function is concave in $\lambda_s$. 
Solve the first order condition, i.e., $\frac{\partial \pi(\lambda_s,W_s,\lambda_o,W_o)}{\partial \lambda_s}|_{\bar \lambda_s(W_s,\lambda_o,W_o)}=0$, we derive  $\bar \lambda_s(W_s,\lambda_o,W_o)= \frac{\Lambda (V \mu_s - w_s)}{2 \mu_s W_s} - \frac{\lambda_o W_o}{W_s}$. Note that the feasible region of $\lambda_s$ is $0 < \lambda_s \leq \Lambda - \lambda_o$, so the best $\lambda_s$ has three possibilities.
\begin{itemize}
    \item $\lambda_s= 0$ if $V\Lambda - \frac{\Lambda w_s}{\mu_s } - 2 \lambda_o W_o < 0$, i.e., $(W_s, \lambda_o,W_o) \in \mathcal{D}_4$.   Then the optimal profit in this subcase is $\pi^O$.
    \item $\lambda_s = \bar \lambda_s (W_s)$ if $V\Lambda - \frac{\Lambda w_s}{\mu_s } - 2 \lambda_o W_o \geq 0$ and $W_s \geq \frac{V\Lambda - \frac{\Lambda w_s}{\mu_s } - 2 \lambda_o W_o }{2(\Lambda - \lambda_o) } $, i.e., $(W_s, \lambda_o,W_o) \in \mathcal{D}_5$. In this case, the market is partially covered. Let 
\bea 
\pi_{sl,p}=\left\{  \begin{array}{rl}
\max_{} & \pi( \lambda_s, W_s, \lambda_o,W_o)\big|_{ \{ \lambda_s = \bar \lambda_s(W_s)\} }  = \pi_2(\lambda_o ,W_o) + \frac{(V\Lambda  - \frac{\Lambda  w_s}{\mu_s}- 2\lambda_oW_o)^2 - \frac{4\Lambda w_s}{\mu_s} }{4\Lambda W_s} \\
\mbox{s.t.} & 0 \leq  \lambda_o \leq \Lambda, W_s \geq \max\{ W_o, \frac{V\Lambda - \frac{\Lambda w_s}{\mu_s } - 2 \lambda_o W_o }{2(\Lambda - \lambda_o) } \}. 
\end{array}\right.
\eea
    \item $\lambda_s = \Lambda - \lambda_o$ if $V\Lambda - \frac{\Lambda w_s}{\mu_s } - 2 \lambda_o W_o \geq 0$ and $W_s \leq \frac{V\Lambda - \frac{\Lambda w_s}{\mu_s } - 2 \lambda_o W_o }{2(\Lambda - \lambda_o) } $, i.e., $(W_s, \lambda_o,W_o) \in \mathcal{D}_6$. In this case, the market is fully covered. Let 
\bea \label{model-slf}
\pi_{sl,f}= \left\{  \begin{array}{rl}
\max_{} & \pi( \lambda_s, W_s, \lambda_o,W_o)\big|_{ \{ \lambda_s = \Lambda - \lambda_o \} } \\
& = \pi_2(\lambda_o ,W_o) + V(\Lambda - \lambda_o) - \frac{(\Lambda -\lambda_o)^2 W_s}{\Lambda} - \frac{2\lambda_o (\Lambda - \lambda_o) W_o}{\Lambda} - \frac{w_s}{\mu_s}(\frac{1}{W_s} + \Lambda - \lambda_o) \\
\mbox{s.t.} & 0 \leq  \lambda_o \leq \Lambda, W_o \leq W_s \leq \frac{V\Lambda  - \frac{\Lambda  w_s}{\mu_s} - 2\lambda_oW_o}{2(\Lambda - \lambda_o)}. 
\end{array}\right.
\eea
\end{itemize}

Next we show that 
\begin{align}
\pi_{sl,p} \leq \max\{ \pi^O, \pi_{sl,f}, \pi_{ol,f} \}.\label{pislp}
\end{align}
Note that $\pi(\bar \lambda_s(W_s),W_s,\lambda_o,W_o)$ is monotonic in $W_s$.  
Together with constraint $W_s \geq \max\{ W_o, \frac{V\Lambda - \frac{\Lambda w_s}{\mu_s } - 2 \lambda_o W_o }{2(\Lambda - \lambda_o)} \}$, the best $W_s$ is either $W_o$, $\frac{V\Lambda - \frac{\Lambda w_s}{\mu_s } - 2 \lambda_o W_o}{2(\Lambda - \lambda_o)}$ or infinite. There are three situations depending on the value of $W_s$.   
Suppose $\pi_1(\Lambda)   \leq 2 \lambda_o W_o$, i.e., $(V\Lambda  - \frac{\Lambda  w_s}{\mu_s} - 2\lambda_oW_o)^2 - \frac{4\Lambda w_s}{\mu_s} < 0$, then $W_s  $ approaching  infinity is optimal. That is, System O is optimal. So we have $\pi_{sl,p} \leq \pi^O$.  
Suppose $\pi_1(\Lambda) \geq 2 \lambda_o W_o$, then  $W_s = W_o$ if $W_o \geq \frac{V\Lambda - \frac{\Lambda w_s}{\mu_s } - 2 \lambda_o W_o}{2(\Lambda - \lambda_o)}$. That is,  standard and on-demand services have identical leadtime. In this situation, $\pi_{sl,p} = A_e $. It follows from \eqref{20240725-1} that $A_e \leq \max\{ \pi^O, \pi_{ol,f} \}$, which implies   $\pi_{sl,p} \leq \max\{ \pi^O, \pi_{ol,f} \}$.  
Suppose $\pi_1(\Lambda) \geq 2 \lambda_o W_o$, then  $W_s = \frac{V\Lambda - \frac{\Lambda w_s}{\mu_s } - 2 \lambda_o W_o}{2(\Lambda - \lambda_o)}$ if $W_o \leq \frac{V\Lambda - \frac{\Lambda w_s}{\mu_s } - 2 \lambda_o W_o}{2(\Lambda - \lambda_o)}$.   That implies $\lambda_s = \Lambda - \lambda_o$. Hence 
\[ \pi_{sl,p}  = \pi( \lambda_s, W_s, \lambda_o,W_o)\big|_{ \{ \lambda_s = \Lambda -\lambda_o, W_s = \frac{V\Lambda - \frac{\Lambda w_s}{\mu_s } - 2 \lambda_o W_o}{2(\Lambda - \lambda_o)}\} } \leq \pi_{sl,f}, \]
where the inequality is because $(\lambda_s,W_s)=  ( \Lambda -\lambda_o,    \frac{V\Lambda - \frac{\Lambda w_s}{\mu_s } - 2 \lambda_o W_o}{2(\Lambda - \lambda_o)})$ is feasible solution for $\pi_{sl,f}$. As such, we have $\pi_{sl,p} \leq \max \{ \pi^O, \pi_{sl,f}, \pi_{ol,f} \}$.

According to above two cases and \eqref{piolp} and \eqref{pislp}, we have
\[ \pi^* = \max\{\pi^S,\pi^O, \pi^H_{\lambda_s, W_s,\lambda_o,W_o} \} = \max\{\pi^S,\pi^O, \pi_{ol,f}, \pi_{sl,f} \}.  \]
Note that the best $\lambda_s$ in both $\pi_{ol,f}$ and $\pi_{sl,f}$ is $\Lambda - \lambda_o$, so we can rewrite $ \max\{\pi_{ol,f}, \pi_{sl,f} \}$ as 
\[\max\{\pi_{ol,f}, \pi_{sl,f} \}= \max_{(W_s,\lambda_o, W_o) \in \mathcal{D}_3 \cup \mathcal{D}_6} \  \pi(\lambda_s, W_s, \lambda_o,W_o)\big|_{ \{\lambda_s= \Lambda - \lambda_o\}}  
\equiv\pi^H_{W_s, \lambda_o,W_o}  ,\]
and then $\pi^* =  \max\{\pi^S, \pi^O, \pi^H_{W_s, \lambda_o,W_o}\}$. 
\qed

\noindent
{\bf Proof of Lemma \ref{lem-solve-Ws}.} 
In this lemma, we work on $\pi^H_{W_s, \lambda_o,W_o}$ to derive the best $W_s$ for any fixed $(\lambda_o,W_o)$. Following from Lemma \ref{lem-solve-ls}, the feasible region of $\pi^H_{W_s, \lambda_o,W_o}$ is $\mathcal{D}_3 \cup \mathcal{D}_6$. We will proceed in three steps. In step 1, we solve the best $W_s$ when $(W_s, \lambda_o, W_o) \in \mathcal{D}_3$. In step 2, we solve the best $W_s$ when $(W_s, \lambda_o, W_o) \in \mathcal{D}_6$. In step 3,  we compare and combine the results from the two steps. 

Step 1: $(W_s, \lambda_o, W_o) \in \mathcal{D}_3$. In this step, we solve the best $W_s$ in  \eqref{model-olf}. 
Let $\mathcal{R} = \{ (\lambda_o, W_o )| W_o > 0, 0 < \lambda_o \leq \Lambda \}$, and 
\beq
\mathcal{R}_{ol,1} &=& \{ (\lambda_o, W_o )| 0 \leq  \lambda_o \leq \sqrt{\Lambda^2 - \frac{4\Lambda w_s\mu_s}{(V\mu_s -w_s)^2}},  W_o \geq \sqrt{\frac{\Lambda w_s}{\mu_s (\Lambda + \lambda_o)(\Lambda - \lambda_o)}}  \}, \\
\mathcal{R}_{ol,2} &=& \{ (\lambda_o, W_o )|  \sqrt{\Lambda^2 - \frac{4\Lambda w_s\mu_s}{(V\mu_s -w_s)^2}} < \lambda_o \leq \Lambda,  W_o > \frac{ V\mu_s-w_s}{2\mu_s} \}, \\
\mathcal{R}_{ol,3} &=& \{ (\lambda_o, W_o )| 0 \leq \lambda_o \leq \Lambda,0< W_o \leq \min\{ \frac{ V\mu_s-w_s}{2\mu_s} , \sqrt{\frac{\Lambda w_s}{\mu_s (\Lambda + \lambda_o)(\Lambda - \lambda_o)}} \}  \}.
\eeq 
Note that $\mathcal{R} = \cup_{i=1}^3 \mathcal{R}_{ol,i}$, see Figure \ref{fig-Ws} for graphical illustration. Define the objective function in \eqref{model-olf} as $\pi_{ol,f}(W_s,\lambda_o,W_o)$, i.e., 
\bea \label{model-olf-1}
\pi_{ol,f} && =\left\{  \begin{array}{rl}
\max_{} & \pi_{ol,f}(W_s,\lambda_o,W_o) \equiv \pi_2(\lambda_o ,W_o) + V(\Lambda - \lambda_o) - \frac{(\Lambda -\lambda_o)^2 W_s}{\Lambda} - \frac{2\lambda_o (\Lambda - \lambda_o) W_s}{\Lambda} - \frac{w_s}{\mu_s}(\frac{1}{W_s} + \Lambda - \lambda_o). \\
\mbox{s.t.} & 0\leq  \lambda_o \leq \Lambda, \ W_o \geq W_s, \ W_s \leq \frac{ V\mu_s-w_s}{2\mu_s }. 
\end{array}\right. .
\eea  
The objective function $\pi_{ol,f}(W_s,\lambda_o,W_o)$ is concave in $W_s$ since $\frac{\partial^2 \pi_{ol,f}(W_s,\lambda_o,W_o)}{\partial W_s^2} = - \frac{2w_s}{\mu_s W_s^3} \leq 0$. The first order derivative, i.e., $\frac{\partial \pi_{ol,f}(W_s,\lambda_o,W_o)}{\partial W_s}$, equals to zero when $W_s = \sqrt{\frac{\Lambda w_s}{\mu_s (\Lambda + \lambda_o)(\Lambda - \lambda_o)}}$. So the best $W_s$ has three possibilities. 
\begin{itemize}
    \item   $W_s = \sqrt{\frac{\Lambda w_s}{\mu_s (\Lambda + \lambda_o)(\Lambda - \lambda_o)}}$ if  $\sqrt{\frac{\Lambda w_s}{\mu_s (\Lambda + \lambda_o)(\Lambda - \lambda_o)}} \leq \min\{ \frac{ V\mu_s-w_s}{2\mu_s },W_o\}$, i.e., $(\lambda_o,W_o) \in \mathcal{R}_{ol,1}$. We use $\pi^{ol,f}_{\mathcal{R}_{ol,1}}$ to denote the corresponding largest profit. 
    \item     $W_s = \frac{ V\mu_s-w_s}{2\mu_s}$ if $W_o > \frac{ V\mu_s-w_s}{2\mu_s}$ and $\frac{ V\mu_s-w_s}{2\mu_s} \leq \sqrt{\frac{\Lambda w_s}{\mu_s (\Lambda + \lambda_o)(\Lambda - \lambda_o)}}$, i.e., $(\lambda_o,W_o) \in \mathcal{R}_{ol,2}$. We use $\pi^{ol,f}_{\mathcal{R}_{ol,2}}$ to denote the corresponding largest profit. 
    \item   $W_s = W_o$ if $W_o \leq \frac{ V\mu_s-w_s}{2\mu_s} $ and $W_o \leq \sqrt{\frac{\Lambda w_s}{\mu_s (\Lambda + \lambda_o)(\Lambda - \lambda_o)}}$,  i.e., $(\lambda_o,W_o) \in \mathcal{R}_{ol,3}$. We use $\pi^{ol,f}_{\mathcal{R}_{ol,3}}$ to denote the corresponding largest profit. 
\end{itemize}

\begin{figure}[!hbt]
	\vspace{-0.1in}
	\FIGURE
	{
            \includegraphics[width=0.5\textwidth]{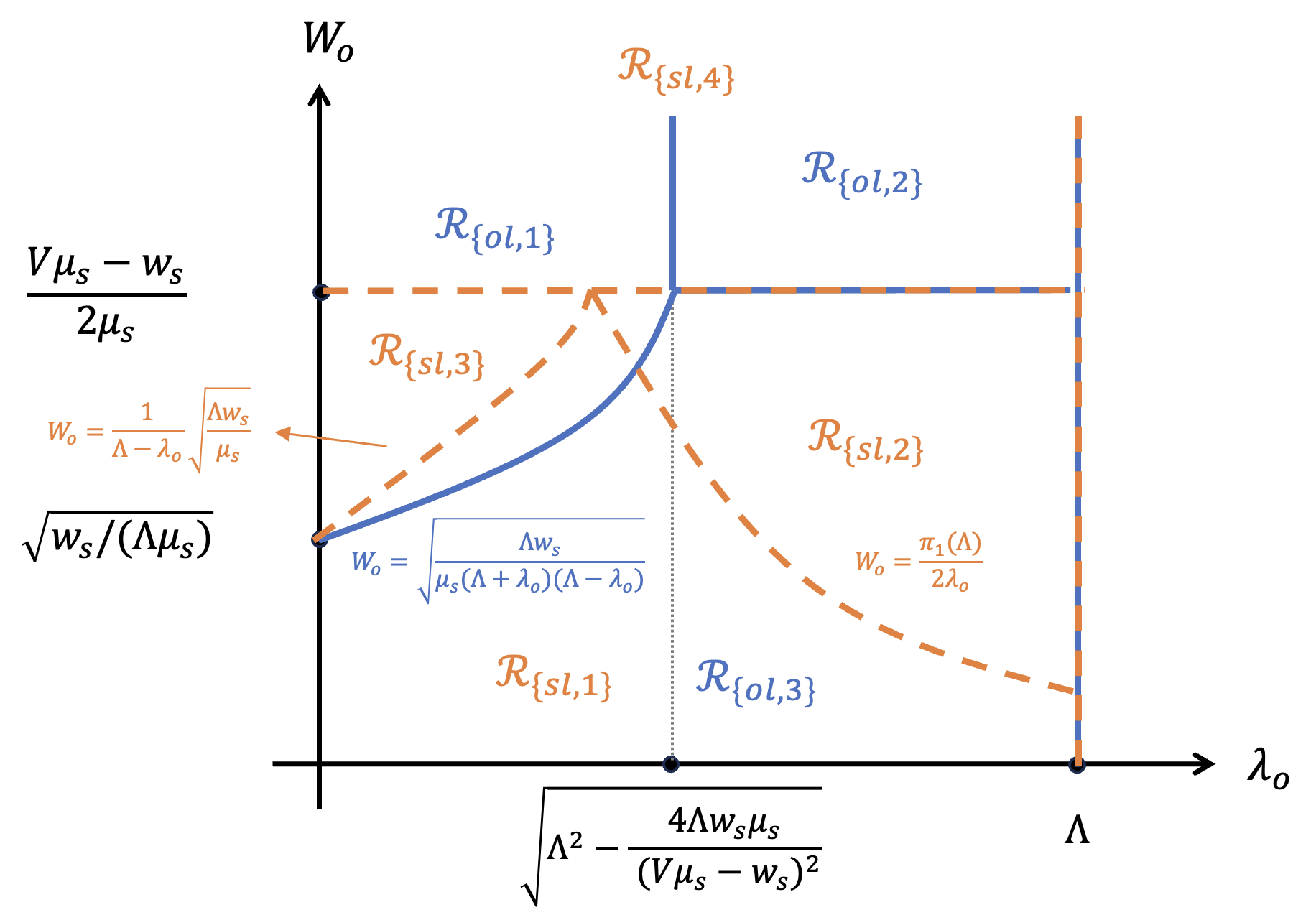} 
	}
	{\bf Illustration of $\mathcal{R}$ in $(\lambda_o,W_o)$ space \label{fig-Ws}}
	{The blue solid lines demarcate the regions $\mathcal{R}_{ol,i}$ (for $i=1,2,3$), and the orange dashed lines delineate the regions $\mathcal{R}_{sl,i}$ (for $i=1,2,3,4$).}
	\vspace{-0.2in}
\end{figure}

Step 2: $(W_s, \lambda_o, W_o) \in \mathcal{D}_6$. In this step, we solve  \eqref{model-slf}. 
Let 
\bea
\mathcal{R}_{sl,1} &=& \{ (\lambda_o, W_o )| 0\leq \lambda_o \leq \Lambda, 0 < W_o \leq \min\{ \frac{V\mu_s - w_s}{2\mu_s},  \frac{1}{\Lambda - \lambda_o} \sqrt{\frac{\Lambda w_s}{\mu_s}}, \frac{\pi_1(\Lambda)}{2\lambda_o} \}  \}, \\
\mathcal{R}_{sl,2} &=& \{ (\lambda_o, W_o )|  0 \leq  \lambda_o \leq \Lambda, \frac{\pi_1(\Lambda)}{2\lambda_o} < W_o \leq \frac{V\mu_s - w_s}{2\mu_s} \}, \\
\mathcal{R}_{sl,3} &=& \{ (\lambda_o, W_o )|  0 \leq   \lambda_o \leq \Lambda, \frac{1}{\Lambda - \lambda_o} \sqrt{\frac{\Lambda w_s}{\mu_s}} <  W_o \leq \frac{V\mu_s - w_s}{2\mu_s} \}, \\
\mathcal{R}_{sl,4} &=& \{ (\lambda_o, W_o )|  0 \leq  \lambda_o \leq \Lambda, W_o > \frac{V\mu_s - w_s}{2\mu_s}  \}. 
\eea 
Note that $\mathcal{R} = \cup_{i=1}^4 \mathcal{R}_{sl,i}$, see Figure \ref{fig-Ws} for graphical illustration.  Define the objective function in  \eqref{model-slf} as $\pi_{sl,f}(W_s,\lambda_o,W_o)$, i.e., 
\bea \label{model-slf-1}
\pi_{sl,f}=\left\{  \begin{array}{rl}
\max_{} & \pi_{sl,f}(W_s,\lambda_o,W_o) \equiv \pi_2(\lambda_o ,W_o) + V(\Lambda - \lambda_o) - \frac{(\Lambda -\lambda_o)^2 W_s}{\Lambda} - \frac{2\lambda_o (\Lambda - \lambda_o) W_o}{\Lambda} - \frac{w_s}{\mu_s}(\frac{1}{W_s} + \Lambda - \lambda_o) \\
\mbox{s.t.} & 0 \leq  \lambda_o \leq \Lambda, \  0 < W_o \leq W_s \leq \frac{V\Lambda  - \frac{\Lambda  w_s}{\mu_s} - 2\lambda_oW_o}{2(\Lambda - \lambda_o)}. 
\end{array}\right.
\eea
Note that the feasible region of $W_s$ is non-empty if and only if $W_o \leq \frac{V\mu_s - w_s}{2\mu_s}$. So \eqref{model-slf-1} is infeasible on $\mathcal{R}_{sl,4}$. We only discuss the best $W_s$ on $\mathcal{R}_{sl,i}$ for $i=1,2,3$ in the following.  The objective function $\pi_{sl,f}(W_s,\lambda_o,W_o)$ is concave in $W_s$ since $\frac{\partial^2 \pi_{sl,f}(W_s,\lambda_o,W_o)}{\partial W_s^2} = - \frac{2w_s}{\mu_s W_s^3} \leq 0$. The first order derivative, i.e., $\frac{\partial \pi_{sl,f}(W_s,\lambda_o,W_o)}{\partial W_s}$, equals to zero when $W_s = \frac{1}{\Lambda - \lambda_o} \sqrt{\frac{\Lambda w_s}{\mu_s}}$. So the best $W_s$ has three possibilities. 
\begin{itemize}
 \item $W_s = \frac{1}{\Lambda - \lambda_o} \sqrt{\frac{\Lambda w_s}{\mu_s}}$ if $W_o < \frac{1}{\Lambda - \lambda_o} \sqrt{\frac{\Lambda w_s}{\mu_s}} < \frac{V\Lambda  - \frac{\Lambda  w_s}{\mu_s} - 2\lambda_oW_o}{2(\Lambda - \lambda_o)}$, i.e., $(\lambda_o,W_o) \in \mathcal{R}_{sl,1}$. We use $\pi^{sl,f}_{\mathcal{R}_{sl,1}}$ to denote the corresponding largest profit. 
 \item  $W_s = \frac{V\Lambda  - \frac{\Lambda  w_s}{\mu_s} - 2\lambda_oW_o}{2(\Lambda - \lambda_o)}$ if $\frac{1}{\Lambda - \lambda_o} \sqrt{\frac{\Lambda w_s}{\mu_s}} > \frac{V\Lambda  - \frac{\Lambda  w_s}{\mu_s} - 2\lambda_oW_o}{2(\Lambda - \lambda_o)}$, i.e., $(\lambda_o,W_o) \in \mathcal{R}_{sl,2}$. We use $\pi^{sl,f}_{\mathcal{R}_{sl,2}}$ to denote the corresponding largest profit. 
\item   $W_s = W_o$ if $W_o \geq \frac{1}{\Lambda - \lambda_o} \sqrt{\frac{\Lambda w_s}{\mu_s}}$, i.e., $(\lambda_o,W_o) \in \mathcal{R}_{sl,3}$. We use $\pi^{sl,f}_{\mathcal{R}_{sl,3}}$ to denote the corresponding largest profit. 
\end{itemize}

Step 3. We compare and combine the results from the previous two steps. Following from step 1 and step 2, we have $\pi^H_{W_s, \lambda_o,W_o} = \max\{ \pi^{ol,f}_{\mathcal{R}_{ol,1}}, \pi^{ol,f}_{\mathcal{R}_{ol,2}},\pi^{ol,f}_{\mathcal{R}_{ol,3}} ,\pi^{sl,f}_{\mathcal{R}_{sl,1}}, \pi^{sl,f}_{\mathcal{R}_{sl,2}},\pi^{sl,f}_{\mathcal{R}_{sl,3}} \}$. 

First, we show that the best $W_s$ can not be $W_o$, i.e., $\pi^H_{W_s, \lambda_o,W_o} \neq \pi^{ol,f}_{\mathcal{R}_{ol,3}} $ and $\pi^H_{W_s, \lambda_o,W_o} \neq \pi^{sl,f}_{\mathcal{R}_{sl,3}} $. 
Note that $\pi_{ol,f}(W_s,\lambda_o,W_o)|_{\{W_s = W_o\}} = \pi_{sl,f}(W_s,\lambda_o,W_o)|_{\{W_s = W_o\}}$, so that the best $W_s= W_o$ if and only if $W_s = W_o$ is the optimal solution for both problems, i.e., \eqref{model-olf-1} and \eqref{model-slf-1}. We consider two regions: $W_o >\frac{V\mu_s - w_s}{2\mu_s}$ and $W_o \leq \frac{V\mu_s - w_s}{2\mu_s}$. Suppose $W_o >\frac{V\mu_s - w_s}{2\mu_s}$, then $W_s = W_o$ is infeasible in both \eqref{model-olf-1} and \eqref{model-slf-1}. Suppose $W_o \leq \frac{V\mu_s - w_s}{2\mu_s}$, then for \eqref{model-olf-1}, the best $W_s = W_o$ if and only if $W_o \leq \sqrt{\frac{\Lambda w_s}{\mu_s (\Lambda + \lambda_o)(\Lambda - \lambda_o)}}$; for \eqref{model-slf-1}, the best $W_s = W_o$ if and only if $W_o \geq \frac{1}{\Lambda - \lambda_o} \sqrt{\frac{\Lambda w_s}{\mu_s}}$. Since $ \sqrt{\frac{\Lambda w_s}{\mu_s (\Lambda + \lambda_o)(\Lambda - \lambda_o)}} < \frac{1}{\Lambda - \lambda_o} \sqrt{\frac{\Lambda w_s}{\mu_s}}$ for any $\lambda_o \in (0,\Lambda)$, $W_s = W_o$ can not be the optimal solution for both problems. Therefore, we have 
\bea 
\pi^H_{W_s, \lambda_o,W_o} = \max\{ \pi^{ol,f}_{\mathcal{R}_{ol,1}}, \pi^{ol,f}_{\mathcal{R}_{ol,2}},\pi^{sl,f}_{\mathcal{R}_{sl,1}}, \pi^{sl,f}_{\mathcal{R}_{sl,2}} \}. \label{eq-Ws-step3-1}
\eea 


Second, we show that $\pi^{ol,f}_{\mathcal{R}_{ol,2}}\leq \pi^O$ and $\pi^{sl,f}_{\mathcal{R}_{sl,2}}\leq \pi^O$. 
If $(\lambda_o,W_o)\in \mathcal{R}_{ol,2}$, then for problem \eqref{model-olf-1}, the best $W_s$ is $\frac{ V\mu_s-w_s}{2\mu_s}$ and the corresponding profit is 
\[ \pi_{ol,f}( W_s, \lambda_o,W_o)|_{\{W_s = \frac{ V\mu_s-w_s}{2\mu_s}\} } = \pi_2(\lambda_o,W_o)  + \frac{(\Lambda - \lambda_o)^2 (V - w_s/ \mu_s)}{2\Lambda} - \frac{2w_s}{V\mu_s -w_s}. \]
Note that  $\frac{(\Lambda - \lambda_o)^2 (V - w_s/ \mu_s)}{2\Lambda} - \frac{2w_s}{V\mu_s -w_s} \geq 0$ if and only if $\lambda_o \leq \Lambda - \sqrt{\frac{4w_s\mu_s \Lambda}{(V\mu_s -w_s)^2}}$, and $\mathcal{R}_{ol,2}$ requires $\lambda_o \geq \sqrt{\Lambda^2 - \frac{4w_s\mu_s \Lambda}{(V\mu_s -w_s)^2}}$. 
Since $\Lambda - \sqrt{\frac{4w_s\mu_s \Lambda}{(V\mu_s -w_s)^2}} \leq \sqrt{\Lambda^2 - \frac{4w_s\mu_s \Lambda}{(V\mu_s -w_s)^2}}$ when $\pi_1(\Lambda) \geq 0$, we have $\pi_{ol,f}( W_s, \lambda_o,W_o)|_{\{W_s = \frac{ V\mu_s-w_s}{2\mu_s}\} } \leq \pi_2(\lambda_o,W_o)$ for any $(\lambda_o,W_o)\in \mathcal{R}_{ol,2}$. This implies $\pi^{ol,f}_{\mathcal{R}_{ol,2}}\leq \pi^O$.

If $(\lambda_o,W_o)\in \mathcal{R}_{sl,2}$, then for problem \eqref{model-slf-1}, the best $W_s$ is $\frac{V\Lambda  - \frac{\Lambda  w_s}{\mu_s} - 2\lambda_oW_o}{2(\Lambda - \lambda_o)}$ and the corresponding profit is 
\[ \pi_{sl,f}( W_s, \lambda_o,W_o)|_{\{ W_s= \frac{V\Lambda  - \frac{\Lambda  w_s}{\mu_s} - 2\lambda_oW_o}{2(\Lambda - \lambda_o)}\} } = \pi_2(\lambda_o,W_o) + \frac{(\Lambda - \lambda_o)(V\Lambda - \frac{\Lambda w_s}{\mu_s} - 2\lambda_o W_o)}{2\Lambda } - \frac{2w_s(\Lambda - \lambda_o) }{V\Lambda - \frac{\Lambda w_s}{\mu_s} - 2\lambda_o W_o} .\]
Note that $\frac{(\Lambda - \lambda_o)(V\Lambda - \frac{\Lambda w_s}{\mu_s} - 2\lambda_o W_o)}{2\Lambda } - \frac{2w_s(\Lambda - \lambda_o) }{V\Lambda - \frac{\Lambda w_s}{\mu_s} - 2\lambda_o W_o} \geq 0$ if and only if $\lambda_o W_o \leq \pi_1(\Lambda)$, and $\mathcal{R}_{sl,2}$ requires $\lambda_o W_o \geq \pi_1(\Lambda)$. Hence we have $\pi_{sl,f}( W_s, \lambda_o,W_o)|_{ \{ W_s= \frac{V\Lambda  - \frac{\Lambda  w_s}{\mu_s} - 2\lambda_oW_o}{2(\Lambda - \lambda_o)}\} } \leq \pi_2(\lambda_o,W_o)$ for any $(\lambda_o,W_o)\in \mathcal{R}_{sl,2}$. This implies $\pi^{sl,f}_{\mathcal{R}_{sl,2}}\leq \pi^O$. Therefore, we have
\bea 
\max\{ \pi^{ol,f}_{\mathcal{R}_{ol,2}}, \pi^{sl,f}_{\mathcal{R}_{sl,2}} \} \leq \pi^O. \label{eq-Ws-step3-2}
\eea 

Third, it follows from $\pi^* =  \max\{\pi^S, \pi^O, \pi^H_{W_s, \lambda_o,W_o}\}$, \eqref{eq-Ws-step3-1} and \eqref{eq-Ws-step3-2} that 
\bea 
\pi^* = \max\{\pi^S, \pi^O, \pi^{ol,f}_{\mathcal{R}_{ol,1}}, \pi^{sl,f}_{\mathcal{R}_{sl,1}} \}. \label{eq-Ws-step3-3}
\eea 
Note that  $\mathcal{R}_{ol,1}\cap \mathcal{R}_{sl,1}\neq  \emptyset$,  so we   compare which is better. 
Suppose $(\lambda_o,W_o) \in \mathcal{R}_{ol,1} \cap \mathcal{R}_{sl,1}$. Let $\pi_{ol,1} = \pi_{ol,f}( W_s, \lambda_o,W_o)|_{ \{W_s=\sqrt{\frac{\Lambda w_s}{\mu_s (\Lambda + \lambda_o)(\Lambda - \lambda_o)}}\} }$ and $\pi_{sl,1} = \pi_{sl,f}(W_s, \lambda_o,W_o) |_{\{W_s=\frac{1}{\Lambda - \lambda_o} \sqrt{\frac{\Lambda w_s}{\mu_s}}\}}$. Note that $
\pi_{ol,1} - \pi_{sl,1} = \frac{2(\Lambda - \lambda_o)}{\Lambda} ( \lambda_o W_o + \sqrt{\frac{\Lambda w_s}{\mu_s}} - \sqrt{\frac{\Lambda w_s (\Lambda +\lambda_o)}{\mu_s (\Lambda - \lambda_o)}} ) \geq 0$ if and only if $W_o \geq \frac{1}{\lambda_o}( \sqrt{\frac{\Lambda w_s (\Lambda +\lambda_o)}{\mu_s (\Lambda - \lambda_o)}} - \sqrt{\frac{\Lambda w_s}{\mu_s}})$. 
For any $\lambda_o \in (0,\Lambda)$, we have $\sqrt{\frac{\Lambda w_s}{\mu_s (\Lambda + \lambda_o)(\Lambda - \lambda_o)}} < \frac{1}{\lambda_o}( \sqrt{\frac{\Lambda w_s (\Lambda +\lambda_o)}{\mu_s (\Lambda - \lambda_o)}} - \sqrt{\frac{\Lambda w_s}{\mu_s}}) <  \frac{1}{\Lambda - \lambda_o} \sqrt{\frac{\Lambda w_s}{\mu_s}}$, so the curve $ \frac{1}{\lambda_o}( \sqrt{\frac{\Lambda w_s (\Lambda +\lambda_o)}{\mu_s (\Lambda - \lambda_o)}} - \sqrt{\frac{\Lambda w_s}{\mu_s}})$
divides $\mathcal{R}_{ol,1}\cap \mathcal{R}_{sl,1}$ into two regions.  It is straightforward to show that $\min\{\frac{1}{\lambda_o}( \sqrt{\frac{\Lambda w_s (\Lambda +\lambda_o)}{\mu_s (\Lambda - \lambda_o)}} - \sqrt{\frac{\Lambda w_s}{\mu_s}}), \frac{\pi_1(\Lambda)}{2\lambda_o} \} \leq \frac{V\mu_s - w_s}{2\mu_s} $. Define 
\beq
\mathcal{R}_1  &=&  \{ (\lambda_o, W_o )|  0\leq \lambda_o \leq \sqrt{\Lambda^2 - \frac{4\Lambda w_s\mu_s}{(V\mu_s -w_s)^2}} ,  W_o \geq \max\{\sqrt{\frac{\Lambda w_s}{\mu_s (\Lambda + \lambda_o)(\Lambda - \lambda_o)}} ,  \\
&&\quad \quad \quad \quad \quad \quad \quad \quad \min\{    
\frac{1}{\lambda_o}( \sqrt{\frac{\Lambda w_s (\Lambda +\lambda_o)}{\mu_s (\Lambda - \lambda_o)}} - \sqrt{\frac{\Lambda w_s}{\mu_s}}), \frac{\pi_1(\Lambda)}{2\lambda_o} \} \}\},\\
\mathcal{R}_2  &=& \{ (\lambda_o, W_o )|  0\leq \lambda_o \leq \Lambda, 0 < W_o \leq \min\{    
\frac{1}{\lambda_o}( \sqrt{\frac{\Lambda w_s (\Lambda +\lambda_o)}{\mu_s (\Lambda - \lambda_o)}} - \sqrt{\frac{\Lambda w_s}{\mu_s}}), \frac{\pi_1(\Lambda)}{2\lambda_o} \} \}.
\eeq
Then $W_o\geq W_s$ is a better solution than $W_o\leq W_s$ on $\mathcal{R}_{1}$, $W_s\geq W_o$ is a better solution than $W_s\leq W_o$ on $\mathcal{R}_{2}$, and $  \mathcal{R}_{ol,1}\cup \mathcal{R}_{sl,1}= \mathcal{R}_{1}\cup \mathcal{R}_{2}$. See Figure \ref{fig-Ws-1} for illustration. Then we can rewrite $\max\{ \pi^{ol,f}_{\mathcal{R}_{ol,1}}, \pi^{sl,f}_{\mathcal{R}_{sl,1}} \}$ as $\pi^H_{\lambda_o,W_o}$ which is defined in 
\eqref{pilaWo}.  
Together with \eqref{eq-Ws-step3-3}, we have $\pi^* = \{\pi^S, \pi^O, \pi^H_{\lambda_o,W_o} \} $.
\begin{figure}[!hbt]
	\vspace{-0.1in}
	\FIGURE
	{
            \includegraphics[width=0.7\textwidth]{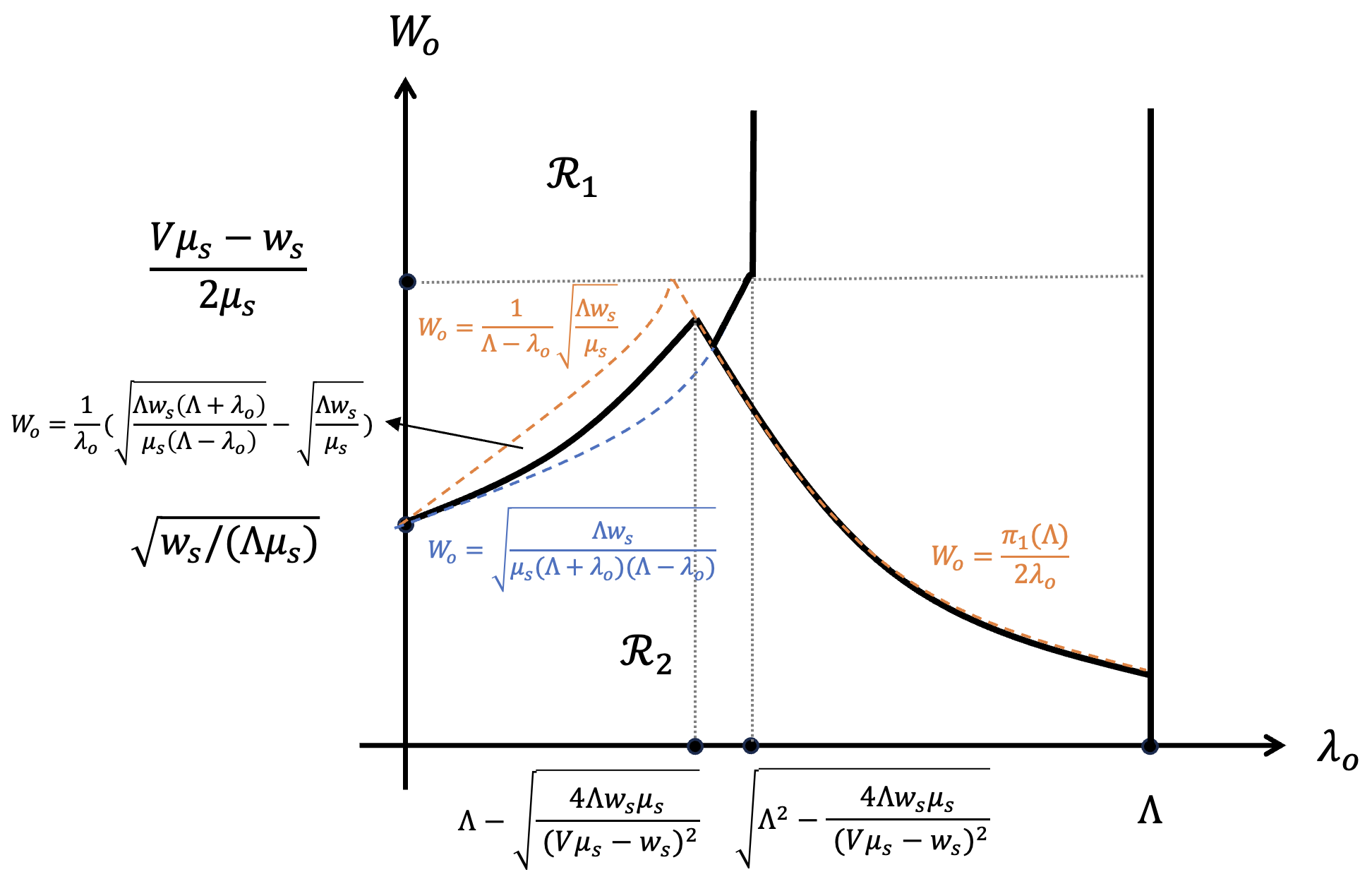} 
	}
	{\bf Illustration of $\mathcal{R}$ in $(\lambda_o,W_o)$ space \label{fig-Ws-1}}
	{}
	\vspace{-0.2in}
\end{figure}
\qed

\noindent
{\bf Proof of Lemma \ref{lem-solve-loWo}.}
Let 
\bee 
  \pi_{ol}(\lambda_o,W_o)& =&  \pi(\lambda_s, W_s, \lambda_o,W_o) \big|_{ \{ \lambda_s = \Lambda - \lambda_o, W_s = \sqrt{\frac{\Lambda w_s}{\mu_s (\Lambda + \lambda_o)(\Lambda - \lambda_o)}}\}}\nonumber \\
&= &\pi_2(\lambda_o,W_o) + (\Lambda - \lambda_o) (V - \frac{w_s }{\mu_s} ) - 2 \sqrt{\frac{w_s (\Lambda + \lambda_o)(\Lambda - \lambda_o)}{\Lambda \mu_s }}, \label{2024-0819-1}\\
 \pi_{sl}(\lambda_o,W_o)& = & \pi(\lambda_s, W_s, \lambda_o,W_o)\big|_{ \{\lambda_s = \Lambda - \lambda_o, W_s = \frac{1}{\Lambda - \lambda_o} \sqrt{\frac{\Lambda w_s}{\mu_s}}\} }  \nonumber\\
&=& \pi_2(\lambda_o,W_o) + (\Lambda - \lambda_o)(V - \frac{w_s}{\mu_s}- 2\sqrt{\frac{w_s}{\Lambda \mu_s}} - \frac{2\lambda_o W_o}{\Lambda}). \label{2024-0819-2}
\eee 
Together with \eqref{pilaWo}, we have 
\begin{align}
\pi^H_{\lambda_o,W_o} = \max \{   \max_{(\lambda_o,W_o)\in \mathcal{R}_1}   \pi_{ol}(\lambda_o,W_o),  \max_{(\lambda_o,W_o)\in \mathcal{R}_2}   \pi_{sl}(\lambda_o,W_o)
\}.\label{eq-solve-Wo}
\end{align}

For fixed $\lambda_o$, $\pi_{ol}(\lambda_o,W_o) $ is concave in $W_o$ since $\frac{\partial^2 \pi_{ol}(\lambda_o, W_o)}{\partial W_o^2} = - \frac{2 (3 + 2\lambda_o W_o)}{K\mu_o^2 W_o^4} < 0$. The first order derivative, i.e., $\frac{\partial \pi_{ol}(\lambda_o, W_o)}{\partial W_o} = - \frac{\lambda_o^2}{\Lambda} + \frac{2(1+ \lambda_o W_o)}{K \mu_o^2 W_o^3}$, equals to zero when $ W_o = \hat W_o (\lambda_o)$, where $ \hat W_o(\lambda_o)$ is given by $\frac{\lambda_o^2}{\Lambda} = \frac{2(1+ \lambda_o W_o)}{K \mu_o^2 W_o^3}$. Therefore, if $(\lambda_o,  \hat W_o(\lambda_o)) \in \mathcal{R}_1$, then $W_{o,ol}(\lambda_o) =  \hat W_o (\lambda_o)$; otherwise, $\pi_{ol}(\lambda_o,W_o)$ decreases in $W_o$ and  $W_{o,ol}(\lambda_o)$ is the smallest feasible $W_o$ that satisfies $(\lambda_o, W_o) \in \mathcal{R}_1$.  

For fixed $\lambda_o$, $\pi_{sl}(\lambda_o,W_o) $ is also concave in $W_o$. The first order derivative, i.e., $\frac{\partial \pi_{sl}(\lambda_o, W_o)}{\partial W_o} = \frac{ 2 ( 1+\lambda_o W_o)}{K \mu_o^2 W_o^3} - \frac{\lambda_o (2\Lambda - \lambda_o)}{\Lambda}$, equals to zero when $W_o = \check  W_o(\lambda_o)$, where $\check W_o(\lambda_o)$ is given by $\frac{\lambda_o(2\Lambda - \lambda_o)}{\Lambda} = \frac{2(1+ \lambda_o W_o)}{K \mu_o^2 W_o^3}$. Therefore, if $(\lambda_o, \check W_o(\lambda_o)) \in \mathcal{R}_2$, then  $W_{o,sl}(\lambda_o) = \check  W_o (\lambda_o)$; otherwise, $\pi_{sl}(\lambda_o,W_o)$ increases in $W_o$ and $W_{o,sl}(\lambda_o)$ is the largest feasible $W_o$ that satisfies $(\lambda_o, W_o) \in \mathcal{R}_2$. This completes the proof. 
\qed 

\noindent
{\bf Proof of Lemma \ref{lem-pi-mono}.}
In this proof we only show the monotonicity  of $\pi^H$ using \eqref{eq-solve-Wo} because  $\pi^H_{\lambda_o,W_o}=\pi^H$. It is straightforward to show the monotonicity of $\pi^O$ and $\pi^S$ and we omit the details. Let $\pi(\lambda_o,W_o) =  \pi_{ol}(\lambda_o,W_o) \mathbb{I}_{\{(\lambda_o, W_o) \in \mathcal{R}_1\}} + \pi_{sl}(\lambda_o,W_o)\mathbb{I}_{\{(\lambda_o, W_o) \in \mathcal{R}_2\}}$ and $\Omega = \mathcal{R}_1 \cup \mathcal{R}_2$. 
Then following from \eqref{eq-solve-Wo}, we have $\pi^H_{\lambda_o,W_o} =  \max_{ (\lambda_o,W_o) \in \Omega } \pi(\lambda_o,W_o) $. 

Let $K_1 < K_2$, $(\lambda_{o1}, W_{o1}) = \arg \max_{(\lambda_o,W_o)\in \Omega}  \pi( \lambda_o,W_o| K=K_1)$, and $( \lambda_{o2}, W_{o2}) = \arg \max_{(\lambda_o,W_o) \in \Omega} \pi(\lambda_o,W_o | K=K_2)$. Then we have 
\[ \pi(\lambda_{o1}, W_{o1} | K = K_1 ) \leq  \pi(\lambda_{o1}, W_{o1} | K = K_2) \leq \pi(\lambda_{o2}, W_{o2} | K = K_2),\]
where the first inequality is because $\pi(\lambda_o,W_o)$ increases in $K$ for any fixed $(\lambda_o,W_o)$, and the second inequality is from optimality. Hence, $\pi^H_{\lambda_o,W_o}$ increases in $K$.  

Let $w_{s1} < w_{s2}$, $(\lambda_{o1}, W_{o1}) = \arg \max_{(\lambda_o,W_o)\in \Omega_1}  \pi(\lambda_o,W_o| w_s = w_{s1})$, and $(\lambda_{o2}, W_{o2}) = \arg \max_{(\lambda_o,W_o) \in \Omega_2} \pi(\lambda_o,W_o | w_s= w_{s2})$, where $\Omega_1$ and $\Omega_2$ are the feasible regions under $w_s = w_{s1}$ and $w_s = w_{s2}$, respectively. 
Then we have 
\[ \pi(\lambda_{o1}, W_{o1} | w_s= w_{s1} ) \geq  \pi(\lambda_{o2}, W_{o2} | w_s= w_{s1}) \geq \pi(\lambda_{o2}, W_{o2} |w_s= w_{s2}),\]
where the first inequality is from optimality and $\Omega_2 \subseteq \Omega_1$, and the second inequality is because $\pi(\lambda_o,W_o)$ decreases in $w_s$ for any fixed $(\lambda_o,W_o)$. Hence, $\pi^H_{\lambda_o,W_o}$ decreases in $w_s$.  \qed 

\noindent
{\bf Proof of Theorem \ref{prop-relation}.}
It follows from Theorem \ref{cor-joint-compare} that $CS^S = \frac{\Lambda W_s^S}{2}$.  
Following from Theorem \ref{opt-two}, when $(w_s ,K) \in \mathcal{H}_1$, we have $\lambda_s = \Lambda - \lambda_o$, $W_s = \sqrt{\frac{\Lambda w_s}{\mu_s (\Lambda + \lambda_o)(\Lambda - \lambda_o)}}$, $W_o = W_{o,ol}(\lambda_o)$.
Let  $\theta_1 = \frac{p_s - p_o}{W_o - W_s}$.
\bea 
CS &=&  \Lambda \int_0^{\theta_1} (V- p_o - \theta W_o) d\theta + \Lambda \int_{\theta_1}^1 (V- p_s - \theta W_s) d\theta \nonumber \\
&=& \Lambda \left[ -\frac{W_o}{2} \theta_1^2 + (V-p_o) \theta_1 - \frac{W_s}{2} + V - p_s + \frac{W_s}{2} \theta_1^2 - (V-p_s) \theta_1 \right] \nonumber \\
&=& \Lambda \left[  - \frac{W_o - W_s}{2} \theta_1^2 + (p_s -p_o) \theta_1 + V- p_s - \frac{W_s}{2}\right]\nonumber \\
&=& \Lambda \left[  V- p_s - \frac{W_s}{2} + \frac{(p_s-p_o)^2}{2(W_o - W_s)}\right] . \label{cs-ol}
\eea 

It follows from \eqref{eq-price-ol} and $\lambda_s = \Lambda - \lambda_o$ that  $p_s = V - \frac{\lambda_o W_s + \lambda_s W_s}{\Lambda} = V - W_s$. Plug $p_s = V-W_s$ into \eqref{cs-ol}, we have 
\[ CS = \Lambda \left[  \frac{W_s}{2} + \frac{(p_s-p_o)^2}{2(W_o - W_s)}\right] > \frac{\Lambda W_s }{2} > \frac{\Lambda W_s^S}{2} = CS^S,\]
where the last inequality is from Theorem \ref{pro:compareTS}. 
\qed 
             
\

\noindent
{\bf Proof of Theorem \ref{lem-flexible}.}
As we discussed above, since a dedicated system is a special case of a flexible system with $q_s=0$ and $q_o=0$,  $\pi^* \le \pi_f^*$. 

In what follows, we show $\pi^* \ge \pi_f^*$ by proving that the SP profit under any $(p_s,p_o,k_s,w_o, q_s,q_o)$ in the flexible system could be dominated by a feasible solution in the dedicated system.  Specifically, we show that for any  $(p_s,p_o,k_s,w_o, q_s,q_o)$ in the flexible system with $q_s\ge 0$ and $q_o \geq 0$, there exist feasible solutions $(p_s',p_o',k_s',w_o',0,0)$ corresponding to the dedicated system satisfying $\pi_f(p_s',p_o',k_s',w_o',0,0)  \geq \pi_{f}(p_s,p_o,k_s,w_o,q_s,q_o)$. The key idea is, for any feasible solution $(p_s,p_o,k_s,w_o, q_s,q_o)$ in the flexible system, we construct corresponding prices $(p_s',p_o')$ in the dedicated system 
such that these two systems result in identical effective arrival rates and expected waiting times for each service. To this end, for the convenience of analysis, we perform a one-to-one variable transformation, from $(p_s,p_o,k_s,w_o,q_s,q_o)$ to $(p_s,p_o,W_s,W_o,q_s,q_o)$, and re-write the optimization problem \eqref{model-flexible} as:
\begin{align} 
&\pi_{f}^* =  \left\{\begin{array}{rl}
 \max\limits_{p_s,p_o,W_s,W_o, q_s,q_o} & \pi_{f}(p_s,p_o,W_s,W_o,q_s,q_o)= p_s \frac{\lambda_s - q_o(\lambda_s+\lambda_o)}{1-q_s-q_o} - \frac{w_s}{\mu_s} (\frac{1}{W_s} + \lambda_s) + p_o \frac{\lambda_o - q_s( \lambda_s + \lambda_o)}{1-q_s-q_o}  - \frac{1}{K \mu_o^2} (\frac{1}{W_o} + \lambda_o)^2    \\
\mbox{s.t.} 
& \eqref{eq-lambda-flexible-1}, \eqref{eq-lambda-flexible-2}.   
\end{array}\right. \nonumber
\end{align}

Next, we show how the $(p_s',p_o')$ can be constructed in the dedicated system such that the waiting times and effective arrival rates under the dedicated system are exactly the same as those under flexible system, and $\pi_f(p_s',p_o',W_s,W_o,0,0)  \geq \pi_{f}(p_s,p_o,W_s,W_o,q_s,q_o)$, or equivalently,   
\begin{align}
 p_s' \lambda_s'+ p_o' \lambda_o'  \geq
 p_s \frac{\lambda_s - q_o(\lambda_s+\lambda_o)}{1-q_s-q_o} + p_o \frac{\lambda_o - q_s( \lambda_s + \lambda_o)}{1-q_s-q_o} .\label{eq-equiva-flexible}
\end{align}
We illustrate the construction of such prices $(p_s',p_o')$ in the following by considering two cases: $W_o > W_s$ and $W_o \leq W_s$.  

{\bf Case 1:} $W_o > W_s$.  In this case, the effective arrival rates \eqref{eq-lambda-flexible-1} - \eqref{eq-lambda-flexible-2} could be simplified as
\bea
{\sf Pr} (U_s \ge 0, U_s \ge U_o)  &=&  {\sf Pr} ( \frac{p_s - p_o}{(1-q_s-q_o)(W_o - W_s) } < \theta \leq  \frac{V-p_s}{(1-q_s)W_s+q_s W_o}),   \label{eq-lambda-flexible-3} \\ 
{\sf Pr} (U_o \ge 0, U_o \ge U_s)  &=&  {\sf Pr} ( \theta \leq \frac{V-p_o}{q_o W_s+ (1-q_o) W_o} , \theta \leq \frac{p_s - p_o}{(1-q_s-q_o) (W_o - W_s)} ). \label{eq-lambda-flexible-4}  
\eea
Note that both $\frac{V-p_o}{q_o W_s+ (1-q_o) W_o}> \frac{p_s - p_o}{(1-q_s-q_o) (W_o - W_s)}$ and $\frac{V-p_s}{(1-q_s)W_s+q_s W_o} > \frac{p_s - p_o}{(1-q_s-q_o) (W_o - W_s)}$ are equivalent to $W_o > \frac{[ (1-q_s -q_o )V + q_o p_s - (1-q_s)p_o ]W_s}{(1-q_s-q_o) V- (1-q_o)p_s + q p_o}$.

{\bf Subcase 1-1:}   $p_o \geq p_s$. Following from \eqref{eq-lambda-flexible-3}-\eqref{eq-lambda-flexible-4} we have ${\sf Pr} (U_s \ge 0, U_s \ge U_o)  = {\sf Pr} (\theta \leq  \frac{V-p_s}{(1-q_s)W_s+q_s W_o})$ and ${\sf Pr} (U_o \ge 0, U_o \ge U_s) =  0$. Together with  \eqref{eq-lambda-flexible-1}-\eqref{eq-lambda-flexible-2} we have
\beq
&& \lambda_s = (1-q_s) \Lambda \frac{V-p_s}{(1-q_s)W_s+q_s W_o}, \ \lambda_o =q_s \Lambda  \frac{V-p_s}{(1-q_s)W_s+q_s W_o},  \\ 
&&\lambda_s' = \Lambda \frac{V-p_s'}{W_s} - \Lambda \frac{p_s' - p_o'}{W_o - W_s} , \ \lambda_o' =\Lambda  \frac{p_s' - p_o'}{W_o - W_s} . 
\eeq
Setting $\lambda_s'=\lambda_s$ and $\lambda_o'=\lambda_o$ yields
 \[p_s' = \frac{q_s V(W_o - W_s) + p_s W_s}{q_s (W_o - W_s) + W_s} \text{ and }p_o' = p_s ,\]
which obviously satisfies \eqref{eq-equiva-flexible} since 
$\lambda_o - q_s(\lambda_s + \lambda_o) = 0$ and $p_s' \geq p_s$. 

{\bf Subcase 1-2:}  $p_o < p_s$ and $W_o \leq \frac{[ (1-q_s -q_o )V + q_o p_s - (1-q_s)p_o ]W_s}{(1-q_s-q_o) V- (1-q_o)p_s + q p_o}$. 
Following from \eqref{eq-lambda-flexible-3}-\eqref{eq-lambda-flexible-4} we have ${\sf Pr} (U_s \ge 0, U_s \ge U_o)  = 0$ and $ {\sf Pr} (U_o \ge 0, U_o \ge U_s) = {\sf Pr} ( \theta \leq \frac{V-p_o}{q_o W_s+ (1-q_o) W_o})$. Together with  \eqref{eq-lambda-flexible-1}-\eqref{eq-lambda-flexible-2} we have
\beq 
&& \lambda_s = q_o \Lambda \frac{V-p_o}{q_o W_s+ (1-q_o) W_o}, \ \lambda_o = (1-q_o) \Lambda \frac{V-p_o}{q_o W_s+ (1-q_o) W_o}, \\  
&&\lambda_s' = \Lambda \frac{V-p_s'}{W_s} - \Lambda \frac{p_s' - p_o'}{W_o - W_s} , \ \lambda_o' =\Lambda  \frac{p_s' - p_o'}{W_o - W_s} . 
\eeq 
Setting $\lambda_s'=\lambda_s$ and $\lambda_o'=\lambda_o$ yields
 \[p_s' = \frac{ (1-q_o) V(W_o - W_s) + p_o W_s}{W_o - q_o (W_o - W_s) } \text{ and }p_o' = p_o ,\]
which obviously satisfies \eqref{eq-equiva-flexible} since 
$\lambda_s - q_o(\lambda_s + \lambda_o) = 0$ and $p_s' \geq p_o$. 

{\bf Subcase 1-3:}   $p_o < p_s$ and $W_o > \frac{[ (1-q_s -q_o )V + q_o p_s - (1-q_s)p_o ]W_s}{(1-q_s-q_o) V- (1-q_o)p_s + q p_o}$. Following from \eqref{eq-lambda-flexible-3}-\eqref{eq-lambda-flexible-4} we have ${\sf Pr} (U_s \ge 0, U_s \ge U_o) =  {\sf Pr} ( \frac{p_s - p_o}{(1-q_s-q_o)(W_o - W_s) } < \theta \leq  \frac{V-p_s}{(1-q_s)W_s+q_s W_o})$ and ${\sf Pr} (U_o \ge 0, U_o \ge U_s)  = {\sf Pr} ( \theta \leq \frac{p_s - p_o}{(1-q_s-q_o) (W_o - W_s)} )$. Together with  \eqref{eq-lambda-flexible-1}-\eqref{eq-lambda-flexible-2} we have
\beq
&& \lambda_s = (1-q_s) \Lambda \frac{V-p_s}{(1-q_s)W_s+q_s W_o} - \Lambda \frac{p_s - p_o}{(W_o - W_s) }, \  \lambda_o = q_s \Lambda \frac{V-p_s}{(1-q_s)W_s+q_s W_o} + \Lambda \frac{p_s - p_o}{(W_o - W_s) } , \   \\
&&\lambda_s' = \Lambda \frac{V-p_s'}{W_s} - \Lambda \frac{p_s' - p_o'}{W_o - W_s} , \ \lambda_o' =\Lambda  \frac{p_s' - p_o'}{W_o - W_s}. 
\eeq 
Setting $\lambda_s'=\lambda_s$ and $\lambda_o'=\lambda_o$ yields
\[ p_s' = \frac{q_s V(W_o - W_s) + p_s W_s}{q_s (W_o - W_s) + W_s} , \ p_o' = p_o, \]
which satisfies \eqref{eq-equiva-flexible} since $p_s' \geq p_s$, $p_o' = p_o$ and $\lambda_s' + \lambda_o' = \frac{\lambda_s - q_o(\lambda_s+\lambda_o)}{1-q_s-q_o} + \frac{\lambda_o - q_s( \lambda_s + \lambda_o)}{1-q_s-q_o}$.

{\bf Case 2:} $W_o \leq W_s$. The proofs are similar as Case 1 and we omit the details. 
\qed

\
\end{APPENDICES}

\end{document}